\documentclass{memo-l}

\usepackage{}
\usepackage{bm}
\def\nz{\ifmmode {I\hskip -3pt N} \else {\hbox {$I\hskip -3pt N$}}\fi}
\def\zz{\ifmmode {Z\hskip -4.8pt Z} \else
       {\hbox {$Z\hskip -4.8pt Z$}}\fi}
\def\qz{\ifmmode {Q\hskip -5.0pt\vrule height6.0pt depth 0pt
       \hskip 6pt} \'else {\hbox
       {$Q\hskip -5.0pt\vrule height6.0pt depth 0pt\hskip 6pt$}}\fi}
\def\rz{\ifmmode {I\hskip -3pt R} \else {\hbox {$I\hskip -3pt R$}}\fi}
\def\cz{\ifmmode {C\hskip -4.8pt\vrule height5.8pt\hskip 6.3pt} \else
       {\hbox {$C\hskip -4.8pt\vrule height5.8pt\hskip 6.3pt$}}\fi}

\def\dist {{\rm \; dist \;}}
\def\i{{\mbox{\raisebox{.5ex}{$\chi$}}}}

\def\supp {{\rm \;Supp\;}}

\def\im{{\rm \;Im\;}}
\def\re{{\rm \;Re\;}}
\def\R{{\rz}} 
\def\Z{{\zz}}
 
\def\N{{\nz}} 

\def\C{{\cz}}
\def\la{\langle}
\def\be{\begin{equation}}
\def\ee{\end{equation}}
\def\ra{\rangle}

\def\Max{{\rm Max}}

\makeatletter
\@addtoreset{equation}{chapter}

\makeatother

\newtheorem{theorem}{Theorem}[chapter]
\newtheorem{lemma}[theorem]{Lemma}
\newtheorem{proposition}[theorem]{Proposition}
\newtheorem{definition}[theorem]{Definition}
\newtheorem{remark}[theorem]{Remark}
\newtheorem{corollary}[theorem]{Corollary}

\makeindex

\begin{document}

\frontmatter

\title{Twisted Pseudodifferential Calculus and Application to  the Quantum Evolution of  Molecules}


\author{Andr\'e Martinez}
\address{Universit\`{a} di Bologna, Dipartimento di Matematica, Piazza di Porta San Donato, 5, I-40127 Bologna, 
Italy}
\curraddr{}
\email{martinez@dm.unibo.it}
\thanks{Investigation supported by University of Bologna. Funds for selected
research topics.}

\author{Vania Sordoni}
\address{Universit\`{a} di Bologna, Dipartimento di Matematica, Piazza di Porta San Donato, 5, I-40127 Bologna, 
Italy}
\curraddr{}
\email{sordoni@dm.unibo.it}
\thanks{Investigation supported by University of Bologna. Funds for selected
research topics.}

\date{April 14, 2006}



\begin{abstract}
We construct an abstract pseudodifferential calculus with operator-valued symbol, suitable for the treatment of Coulomb-type interactions, and we apply it to the study of the quantum evolution of molecules in the Born-Oppenheimer approximation, in the case of the electronic Hamiltonian admitting a local gap in its spectrum. In particular, we show that the molecular evolution can be reduced to the one of a  system of smooth semiclassical operators, the symbol of which can be computed explicitely. In addition, we study the propagation of certain wave packets up to long time values of Ehrenfest order. (This work has been accepted for
publication as part of the Memoirs of the American Mathematical Society and
will be published in a future volume.)
\end{abstract}

\maketitle


\setcounter{page}{4}

\tableofcontents


\mainmatter

\chapter{Introduction}
\setcounter{equation}{0}%
\setcounter{theorem}{0}%
In quantum physics, the evolution of a molecule is described by the initial-value Schr\"odinger system,
\be
\label{schr}
\left\{
\begin{array}{l}
i\partial_t\varphi =H\varphi ;\\
\varphi\left\vert_{t=0}\right. =\varphi_0,
\end{array}
\right.
\ee
where $\varphi_0$ is the initial state of the molecule and $H$ stands for the molecular Hamiltonian involving all the interactions between the particles constituting the molecule (electron and nuclei). In case the molecule is imbedded in an electromagnetic field, the corresponding potentials enter the expression of $H$, too. Typically, the interaction between two particles of  positions $z$ and $z'$, respectively, is of Coulomb type, that is, of the form $\alpha\vert z-z'\vert^{-1}$ with $\alpha \in \R$ constant.
\vskip 0.2cm
In the case of a free molecule, a first approach for studying the system (\ref{schr}) consists in considering bounded initial states only, that is, initial states that are eigenfunctions of the Hamiltonian after removal of the center of mass motion. More precisely, one can split the Hamiltonian into,
$$
H=H_{\rm CM} + H_{\rm Rel},
$$
where the two operators $H_{CM}$ (corresponding to the kinetic energy of the center of mass) and $H_{Rel}$ (corresponding to the relative motion of electrons and nuclei) commute. As a consequence, the quantum evolution factorizes into,
$$
e^{-itH} = e^{-itH_{\rm CM}}\hskip 1pt e^{-itH_{\rm Rel}},
$$
where the (free) evolution $e^{-itH_{\rm CM}}$ of the center of mass can be explicitly computed (mainly because $H_{\rm CM}$ has constant coefficients), while the relative motion $e^{-itH_{\rm Rel}}$ still contains all the interactions (and thus, all the difficulties of the problem). Then, taking $\varphi_0$ of the form,
\be
\label{boindt}
\varphi_0 = \alpha_0 \otimes \psi_j
\ee
where $\alpha_0$ depends on the position of the center of mass only, and $\psi_j$ is an eigenfunction of $H_{\rm Rel}$ with eigenvalue $E_j$, the solution of (\ref{schr}) is clearly given by,
$$
\varphi (t) = e^{-it E_j}(e^{-itH_{\rm CM}}\alpha_0) \otimes \psi_j.
$$
Therefore, in this case, the only real problem is to know sufficiently well the eigenelements of $H_{Rel}$, in order to be able to produce initial states of the form (\ref{boindt}).
\vskip 0.2cm
In 1927, M. Born and R. Oppenheimer \cite{BoOp} proposed a formal method for constructing such an approximation of eigenvalues and eigenfunctions of $H_{Rel}$. This method was based on the fact that, since the nuclei are much heavier than the electrons, their motion is slower and allows the electrons to adapt almost instantaneously  to it. As a consequence, the motion of the electrons is not really perceived by the nuclei, except as a surrounding electric field created by their total potential energy (that becomes a function of the positions of the nuclei). In that way, the evolution of the molecule reduces to that of the nuclei imbedded  in an effective electric potential created by the electrons. Such a reduction (that is equivalent to a decomposition of the problem into two different position-scales) permits, in a second step, to use semiclassical tools  in order to find the eigenelements of the final effective Hamiltonian. 
\vskip 0.2cm
At this point, it is important to observe that this method was formal only, in the sense that it allowed to produce formal series of functions that were (formally) solutions of the eigenvalue problem for $H_{Rel}$, but without any estimates on the remainder terms,  and no information about the possible closeness of these functions to true eigenfunctions, nor to the possible exhaustivity of such approximated eigenvalues. 
\vskip 0.2cm
Many years later,  a first attempt to justify rigorously (from the mathematical point of view) the Born-Oppenheimer approximation (in short: BOA) was made by J.-M. Combes, P. Duclos and R. Seiler \cite{CDS} for the diatomic molecules, with an accuracy of order $h^2$, where $h:=\sqrt{m/M}$ is the square-root of the ratio of the electron masses to nuclear masses. After that, full asymptotics in $h$ were obtained by G. Hagedorn  \cite{Ha2, Ha3}, both in the case of diatomic molecules with Coulomb interactions, and in the case of smooth interactions. In these two cases, these results gave a positive answer to the first question concerning the justification of the BOA, namely, the existence of satisfactory estimates on the remainder terms of the series. Later, by using completely different methods (mostly inspired by the microlocal treatment of semiclassical spectral problems, developed by B. Helffer and J. Sj\"ostrand in \cite{HeSj11}), and in the case of smooth interactions, the first author \cite{Ma1} extended this positive answer to the two remaining questions, that is, the exhaustivity and the closeness of the formal eigenfunctions to the true ones. Although such a method (based on microlocal analysis) seemed to require a lot of smoothness, it appeared that it could be adapted to the case of Coulomb interactions, too, giving rise to a first complete rigorous justification of the BOA in a work by M. Klein, A. Martinez, R. Seiler and X.P. Wang \cite{KMSW}. The main trick, that made possible such an adaptation, consists in a change of variables in the positions of the electrons, that depends in a convenient way of the position (say, $x$) of the nuclei. This permits to make the singularities of the interactions electron-nucleus independent of $x$, and thus, in some sense, to regularize these interactions with respect to $x$. Afterwards,  the standard  microlocal tools (in particular, the pseudodifferential calculus with operator-valued symbols, introduced in \cite{Ba}) can be applied and give the conclusion.
\vskip 0.2cm
Of course, all these justifications concerned the eigenvalue problem for $H_{Rel}$, not the general problem of evolution described in (\ref{schr}). In the general case, one could think about expanding any arbitrary initial state according to the eigenfunctions of $H_{Rel}$, and then apply the previous constructions to each term. However, this would lead to remainder terms quite difficult to estimate with respect to the small parameter $h$, mainly because one would have to mix two types of approximations that have nothing to do with each other: The semiclassical one, and the eigenfunctions expansion one. In other words, this would correspond to handle both functional and microlocal analysis, trying to optimize both of them at the same time. It is folks that such a method is somehow contradictory, and does not produce good enough estimates. For this reason, several authors have looked for an alternative way of studying (\ref{schr}), by trying to adapt Born-Oppenheimer's ideas directly to the problem of evolution.
\vskip 0.2cm
The first results in this direction are due to G. Hagedorn \cite{Ha4, Ha5, Ha6}, and provide complete asymptotic expansions of the solution of (\ref{schr}), in the case of smooth interactions and when the initial state is a convenient perturbation of a single electronic-level state. More precisely, splitting the Hamiltonian into,
$$
H = K_{\rm n}(hD_x) + H_{\rm el}(x),
$$
where $K_{\rm n}(hD_x)$ stands for the quantum kinetic energy of the nuclei, and $H_{\rm el}(x)$ is the so-called electronic Hamiltonian (that may be viewed as acting on the position variables $y$ of the electrons, and depending on the position $x$ of the nuclei), one assumes that $H_{\rm el}(x)$ admits an isolated eigenvalue $\lambda (x)$ (say, for $x$ in some open set of $\R^3$) with corresponding eigenfunction $\psi (x,y)$, and one takes $\varphi_0$ of the form,
$$
\varphi_0 (x,y) = f(x)\psi (x,y) + \sum_{k\geq 1}h^k\varphi_{0,k}(x,y)= f(x)\psi (x,y) + {\mathcal  O}(h),
$$
where $f(x)$ is a coherent state in the $x$-variables.
Then, it is shown that, if the  $\varphi_{0,k}$'s are conveniently chosen, the solution of (\ref{schr}) (with a rescaled time $t\mapsto t/h$) admits an asymptotic expansion of the type,
$$
\varphi_t(x,y) \sim f_t(x)\psi (x,y) + \sum_{k\geq 1}h^k\varphi_{t,k}(x,y),
$$
where all the terms can be explicitly computed by means of the classical flow of the effective Hamiltonian $H_{\rm eff}(x,\xi ):=K_{\rm n}(\xi) + \lambda (x)$.
\vskip 0.2cm
Such a result is very encouraging, since it provides a case in which the relevant information on the initial state is not anymore connected with the point spectrum of  $H_{\rm rel}$, but rather with the  localization in energy of the electrons and the  localization in phase space of the nuclei. This certainly fits much better with the semiclassical intuition of this problem, in accordance with the fact that the classical flow of $H_{\rm eff}(x,\xi )$ is involved.
\vskip 0.2cm
Nevertheless, from a conceptual point of view, something is missing in the previous result. Namely, one would like to have an even closer relation between the complete quantum evolution $e^{-itH/h}$ and some {\it reduced quantum evolution} of the type $e^{-it\tilde H_{\rm eff}(x,hD_x)/h}$, for some $\tilde H_{\rm eff}$ close to $H_{\rm eff}$. In that way, one would be able to use all the well developed semiclassical (microlocal) machinery on the operator $\tilde H_{\rm eff}(x,hD_x)$, in order to deduce many results on its quantum evolution group $e^{-it\tilde H_{\rm eff}(x,hD_x)/h}$ (e.g., a representation of it as a Fourier integral operator). In the previous result, the presence of a coherent state in the expression of $\varphi_0$ has allowed the author to, somehow, by-pass this step, and to relate directly the complete quantum evolution to its semiclassical approximation (that is, to objects involving the underlying classical evolution). However, a preliminary link between $e^{-itH/h}$ and some $e^{-it\tilde H_{\rm eff}(x,hD_x)/h}$ would have the advantage of allowing more general initial states, and, by the use of more sophisticated results of semiclassical analysis,  should permit to have a better understanding of the phenomena related to this approximation. Moreover, as we will see, this preliminary link is usually valid for very large time intervals of the form $[-h^{-N}, h^{-N}]$ with $N\geq 1$ arbitrary, while it is well known that the second step (that is, the semiclassical approximation of $e^{-it\tilde H_{\rm eff}(x,hD_x)/h}$) has, in best cases, the Ehrenfest-time limitation $\vert t\vert ={\mathcal  O}(\ln \frac1{h})$ (see (\ref{sol}) and Theorem \ref{WP} below).
\vskip 0.2cm
The first results concerning a reduced quantum evolution have been obtained recently (and independently) by H. Spohn and S. Teufel in \cite{SpTe}, and by the present authors in \cite{MaSo}. In both cases, it is assumed that, at time $t=0$, the energy of the electrons is localized in some isolated part of the electronic Hamiltonian $H_{\rm el}(x)$. In \cite{SpTe}, the authors find an approximation of  $e^{-itH/h}$ in terms of $e^{-it H_{\rm eff}(x,hD_x)/h}$, and prove an error estimate in ${\mathcal  O}(h)$ (actually, it seems that such a result was already present in a much older, but unpublished, work by A. Raphaelian \cite{Ra}). In \cite{MaSo} (following a procedure of \cite{NeSo, So}, and later reproduced with further applications in \cite{PST, Te}), a whole perturbation $\tilde H_{\rm eff} \sim H_{\rm eff} + \sum_{k\geq 1}h^kH_k$ of $H_{\rm eff}$ is constructed, allowing an error estimate in ${\mathcal  O}(h^\infty)$ for the quantum evolution.
\vskip 0.2cm
However, these two papers have the defect of assuming all the interactions smooth, and thus of excluding the physically interesting case of Coulomb interactions. Here, our goal is precisely to allow this case. More precisely, we plan to mix the arguments of \cite{MaSo} and those of \cite{KMSW} in order to include Coulom-type (or, more generally, Laplace-compact)  singularities of the potentials.
\vskip 0.2cm
In \cite{KMSW}, the key-point consists in a refinement of the Hunziker distorsion method, that leads to a family of $x$-dependent unitary operators (where, for each operator, the nuclei-position variable $x$ has to stay in some small open set) such that, once conjugated with these operators, the electronic Hamiltonian becomes smooth with respect to $x$. Then, by using local pseudodifferential calculus with operator-valued symbols, and various tricky patching techniques, a constructive Feshbach method (through a Grushin problem) is performed and leads to the required result.
\vskip 0.2cm
When reading \cite{KMSW}, however, one has the impression that all the technical difficulties and tricky arguments actually hide a somewhat simpler concept, that should be related to some global pseudodifferential calculus adapted to the singularities of the interactions. In other words, it seems that interactions such as Coulomb electron-nucleus ones are, indeed, smooth with respect to $x$ for some `exotic' differential structure on the $x$-space, and that such a differential structure could be used to construct a complete pseudodifferential calculus (with operator-valued symbols). 
Such considerations (that are absent  in \cite{KMSW}) have naturally led us to the notion of {\it twisted pseudodifferential operator} that we describe in Capters \ref{twao} and \ref{twpdo}. This new tool permits in particular to handle a certain type of partial differential operators with singular operator-valued coefficients, mainly as if their coefficients were smooth. To our opinion, the advantages are at least two. First of all, it simplifies considerably (making them clearer and closer to the smooth case) the 
arguments leading to the reduction of the quantum evolution of a molecule. Secondly, thanks to its abstract setting, we believe that it can be applied in other situations where singularities appear.
\vskip 0.2cm
Roughly speaking, we say that an operator $P$ on $L^2(\R_x^n; {\mathcal  H})$ (${\mathcal  H} =$ abstract Hilbert space) is a twisted $h$-admissible pseudodifferential operator, if each operator $U_jPU_j^{-1}$ (where, for any $j$,  $U_j =U_j(x)$ is a given  unitary operator  defined for $x$ in some open set $\Omega_j\subset\R^n$) is $h$-admissible (e.g., in the sense of \cite{Ba, GMS}). Then, under few general conditions on the finite family $(U_j,\Omega_j)_j$, we show that these operators enjoy all the nice properties of composition, inversion, functional calculus and symbolic calculus, similar to those present in the smooth case. Thanks to this, the general strategy of \cite{MaSo} can essentially be reproduced, and leads to the required reduction of the quantum evolution. 
More precisely, we prove that, if the initial state $\varphi_0$ is conveniently localized in space, in energy, and on a $L$-levels isolated part of the electronic spectrum ($L\geq 1$), then, during a certain interval of time (that can be estimated), its quantum evolution can be described by that of a selfadjoint $L\times L$ matrix  $A=A(x,hD_x)$ of smooth semiclassical pseudodifferential operators in the nuclei-variables, in the sense that one has,
$$
e^{-itH/h}\varphi_0 = {\mathcal  W}^*e^{-itA/h} {\mathcal  W}\varphi_0 +{\mathcal  O}(\la t\ra h^\infty ),
$$
where $\mathcal  W$ is a bounded operator onto $L^2(\R^n)^{\oplus L}$, such that  ${\mathcal  W}{\mathcal  W}^* =1$ and ${\mathcal  W}^*{\mathcal  W}$ is an orthogonal projection (that projects onto a so-called almost-invariant subspace). We refer to Theorem \ref{MAINTH} for a precise statement, and to Theorem \ref{th2} for an even better result in the case where the spectral gap of the electronic Hamiltonian is global. In the particular case $L=1$, this also permits to give a geometrical description (involving the underlying classical Hamilton flow of $A$) of the time interval in which such a reduction is possible.
Then, to make the paper more complete, we consider the case of coherent initial states (in the same spirit as in \cite{Ha5, Ha6}) and, applying a semiclassical result of M. Combescure and D. Robert \cite{CoRo}, we  justify the expansions given in  \cite{Ha6} up to times of order $\ln\frac1{h}$ (at least when the geometry makes it possible).\\\
\vskip 0.2cm
\noindent
Outline of the paper:
\vskip 0.2cm
In Chapter \ref{sect2}, we introduce our notations and assumptions, and we state our main results concerning the reduction of the quantum evolution in the case where the electronic Hamiltonian admits a local gap in its spectrum. In Chapter \ref{modop}, we modify the electronic operator away from the relevant region in $x$, in order to deal with a globally nicer operator, admitting a global gap in its spectrum. Chapters 4 and 5 are devoted to the settlement of an abstract singular pseudodifferential calculus (bounded in Chapter \ref{twao}, and partial differential in Chapter \ref{twpdo}). In Chapter \ref{qis}, following \cite{MaSo}, we construct a quasi-invariant subspace that permits, in Chapter \ref{decmodop}, to have a global reduction of the evolution associated with the modified operator constructed in Chapter \ref{modop}. In Chapters \ref{mainproof} and \ref{proofcor}, we complete the proofs of our main results, and, in Chapter \ref{expexpl}, we give a simple way of computing the effective Hamiltonian. Then, in Chapter
\ref{wp}, we apply these results to study the evolution of wave packets. Chapter \ref{PolMol} treats, more specifically, the case of polyatomic molecules, by showing how it can be inserted  in our general framework. The remaining three chapters are just appendices: Chapter \ref{stpdo} reviews standard results on pseudodifferential calculus; Chapter \ref{psupp} gives an estimate on the propagation-speed of the support (up to ${\mathcal  O}(h^\infty )$) of the solutions of (\ref{schr}); Chapter \ref{technr} contains two technical results used in the paper.
\chapter{Assumptions and Main Results}
\label{sect2}
\setcounter{equation}{0}%
\setcounter{theorem}{0}%
The purpose of this paper is to  investigate the asymptotic behavior as 
$h\rightarrow 0_+$ of the solutions of  the time-dependent Schr\"odinger equation,
\begin{equation}
\label{eqsch}
ih\frac{\partial \varphi}{\partial t} = P(h)\varphi
\end{equation}
with 
\be
\label{P}
P(h)= \boldsymbol{\omega} +Q(x) +W(x),
\ee
where $Q(x)$ ($x\in \R^n$) is a family of selfadjoint operators on some fix Hilbert space ${\mathcal  H}$ with same dense domain ${\mathcal  D}_Q$, $\boldsymbol{\omega}=\sum_{|\alpha|\leq m} c_{\alpha} (x;h)(hD_x)^\alpha$ is a symmetric semiclassical differential operator of order 0 and degree $m$, with  scalar coefficients depending smoothly on $x$, and  $W(x)$ is a non negative function defined almost everywhere on $\R^n$. 
\vskip 0.3cm
Typically, in the case of a molecular system, $x$ stands for the position of the nuclei,
 $Q(x)$ represents the electronic Hamiltonian that includes the electron-electron and nuclei-electron interactions (all of them of Coulomb-type), $\boldsymbol{\omega}$ is the quantized cinetic energy of the nuclei, and $W(x)$ represents the nuclei-nuclei interactions. Moreover, the parameter $h$ is supposed to be small and, in the case of a molecular system,
$h^{-2}$ actually represents  the quotient of electronic and nuclear  masses.  In more general systems, one can also include a magnetic potential and an exterior electric potential both in $\boldsymbol{\omega}$ and $Q(x)$. We refer to Chapter \ref{PolMol} for more details about this case.
\vskip 0.3cm
We make the following assumptions:\\\\
$(H1)$ For all $\alpha ,\beta \in \Z^n_+$ with $|\alpha|\leq m$, $\partial^\beta c_\alpha (x,h) ={\mathcal  O}(1)$ uniformly for $x\in\R^n$ and $h>0$ small enough. Moreover, setting $\omega (x,\xi ;h):=\sum_{|\alpha| \leq m} c_{\alpha} (x;h)\xi^\alpha$, we assume 
that there exists a constant $C_0\geq 1$ such that, for all $(x,\xi )\in  \R^{2n}$ and $h>0$ small enough,
$$
\re\omega (x,\xi ;h)\geq \frac1{C_0}\la \xi\ra^m -C_0.
$$
In particular, Assumption $(H1)$ implies that $m$ is even and $\boldsymbol{\omega}$ is well defined as a selfadjoint operator on $L^2(\R^n)$ (and, by extension, on $L^2(\R^n;{\mathcal  H})$) with domain $H^m(\R^n)$. Moreover, by the Sharp G\aa  rding Inequality (see, e.g., \cite{Ma2}), it is uniformly semi-bounded from below.
\\\\
$(H2)$ 
$W\geq 0$ is $\la D_x\ra^m$-compact on $L^2(\R^n)$, and there exists $\gamma \in  \R$ such that, for all $x\in {\R}^n $, $Q(x)\geq \gamma$ on ${\mathcal  H}$.
\vskip 0.3cm
Assumptions $(H1)-(H2)$ guarantee that, for $h$ sufficiently small,  $P(h)$ can be realized as
a selfadjoint operator on $L^2(\R^n ; {\mathcal  H})$ with domain ${\mathcal  D}(P)\subset H^m(\R^n; {\mathcal  H})\cap L^2(\R^n; {\mathcal  D}_Q)$, and
verifies  
$P(h)\geq \gamma_0$, with $\gamma_0\in  \R$ independent of $h$.
\vskip 0.3cm
(Of course, in the case of a molecular system, $P(h)$ is essentially selfadjoint, and the domain of its selfadjoint extension is $H^2(\R^n\times Y)$, where $Y$ stands for the space of electron positions.)
\vskip 0.3cm
For $L\geq 1$ and $L'\geq 0$, we denote by $\lambda_1(x)\leq\dots\leq \lambda_{L+L'}(x)$ the first $L+L'$ values given by the Min-Max principle for $Q(x)$ on ${\mathcal  H}$, and we make the following local gap assumption on the spectrum $\sigma(Q(x))$ of $Q(x)$:
\\\\
$(H3)$ There exists a contractible bounded  open set $\Omega \subset \R^n$ and $L\geq 1$ such that, for all $x\in  \Omega$, 
$\lambda_1(x), \dots ,\lambda_{L+L'}(x)$ are discrete eigenvalues of $Q(x)$, and one has,
$$
\inf_{x\in \Omega} \dist \left(\sigma(Q(x))\backslash \{\lambda_{L'+1}(x), \dots ,\lambda_{L'+L}(x)\} , \{\lambda_{L'+1}(x), \dots ,\lambda_{L'+L}(x)\}\right) >0.
$$
Furthermore, the spectral projections $\Pi^-_0(x)$, associated with $ \{\lambda_{1}(x), \dots ,\lambda_{L'}(x)\}$, and $\Pi_0(x)$, associated with $ \{\lambda_{L'+1}(x), \dots ,\lambda_{L'+L}(x)\}$,   both depend continuously on $x\in\Omega$.
\vskip 0.5cm
Then, we assume that $P$ can be ``regularized'' with respect to $x$ in $\Omega$, in the following sense: 
\\\\
$(H4)$ There exists  a finite family of bounded open sets $(\Omega_j)_{j=1}^r$ in $\R^n$, a corresponding family of unitary operators $U_j(x)$ ($j=1,\cdots ,r$, $x\in \Omega_j$), and some fix   selfadjoint operator $Q_0\geq C_0$ on ${\mathcal  H}$ with domain ${\mathcal  D}_Q$, such
that (denoting by $U_j$ the unitary operator on $L^2(\Omega_j;{\mathcal  H})\simeq L^2(\Omega_j)\otimes{\mathcal  H}$ induced by the action of $U_j(x)$ on ${\mathcal  H}$),
\begin{itemize}
\item $\Omega =\cup_{j=1}^r\Omega_j$;
\item For all $j=1,\cdots, r$ and $x\in  \Omega_j$, $U_j(x)$ leaves ${\mathcal  D}_Q$ invariant;
\item For all $j$, the operator 
$U_j \boldsymbol{\omega}U_j^{-1}$
 is a semiclassical differential operator with operator-valued symbol, of the form,
\be
\label{conjomega}
 U_j \boldsymbol{\omega}U_j^{-1} =\boldsymbol{\omega}+h\sum_{|\beta|\leq m-1}\omega_{\beta ,j}(x;h)(hD_x)^\beta,
 \ee
 where $\omega_{\beta ,j}Q_0^{\frac{|\beta|}{m} -1}\in C^\infty (\Omega_j; {\mathcal  L}({\mathcal  H}))$ for any $\gamma\in\Z^n_+$ (here, ${\mathcal  L}({\mathcal  H})$ stands for the Banach space of bounded operators on ${\mathcal  H}$), and the quantity  $\Vert\partial_x^\gamma \omega_{\beta ,j}(x;h) Q_0^{\frac{|\beta|}{m} -1}\Vert_{{\mathcal  L}({\mathcal  H})} $ is bounded uniformly with respect to $h$ small enough and locally uniformly with respect to $x\in\Omega_j$;
\item For all $j$, the operators $U_j(x)Q(x)U_j(x)^{-1}$ and $U_j(x)Q_0U_j(x)^{-1}$ are in $C^\infty (\Omega_j; {\mathcal  L}({\mathcal  D}_Q, {\mathcal  H}))$ (where ${\mathcal  L}({\mathcal  D}_Q, {\mathcal  H})$ stands for the space of bounded operators from ${\mathcal  D}_Q$ to ${\mathcal  H}$);
\item $W\in  C^\infty (\cup_{j=1}^r\Omega_j)$;
\item There exists a dense subspace ${\mathcal  H}_\infty \subset {\mathcal  D}_Q\subset {\mathcal  H}$, such that, for any  $v\in{\mathcal  H}_\infty$ and any $j=1,\cdots ,r$, the application $x\mapsto U_j(x)v$ is in $C^\infty (\Omega_j, {\mathcal  D}_Q)$.
\end{itemize}
\vskip 0.3cm
Note that, for physical molecular systems, a construction of such operators $U_j(x)$'s  is made in \cite{KMSW}, and can be performed around any point of $\R^n$ where $W$ is smooth. Moreover, in that case  one can take $Q_0 =-\Delta_y+1$ (where $y$ stands for the position of the electrons), and  the last point in (H4) can be realized by taking  ${\mathcal  H}_\infty =C_0^\infty (Y)$. Again, we refer the interested reader to Chapter \ref{PolMol}. Let us also observe that, in the case $L'+L=1$, one does not need to assume that $\Omega$ is contractible.
\vskip 0.3cm
 For any $\varphi_0\in L^2(\R^n;{\mathcal  H})$ (possibly $h$-dependent) such that $\Vert\varphi_0\Vert_{L^2(K_0^c;{\mathcal  H})}={\mathcal  O}(h^\infty )$ for some compact set $K_0\subset\subset\R^n$, and for any $\Omega'\subset\subset\R^n$ open neighborhood of $K_0$, we set,
$$
T_{\Omega'}(\varphi_0):= \sup\{T>0\, ;\,  \exists K_T\subset\subset\Omega' ,\,\, \sup_{t\in [0,T]}\Vert e^{-itP/h}\varphi_0\Vert_{L^2(K_T^c;{\mathcal  H})}={\mathcal  O}(h^{\infty})\}.
$$
Then, $T_{\Omega'}(\varphi_0)\leq +\infty$, and, if one also assume that $\Vert (1-f(P))\varphi_0\Vert ={\mathcal  O}(h^\infty)$ for some $f\in C_0^\infty(\R)$, 
Theorem \ref{th:appA} in Appendix B shows that, 
$$
T_{\Omega'}(\varphi_0)\geq \frac{2\dist (K_0,\partial\Omega ')}{\Vert \nabla_\xi\omega (x,hD_x) g(P)\Vert},
$$
for any $g\in C_0^\infty(\R)$ verifying $gf=f$.
\vskip 0.2cm
As a main result, we obtain (denoting by $L^2(\R^n)^{\oplus L}$ the space $(L^2(\R^n))^L$ endowed with its natural Hilbert structure),
\begin{theorem}\sl  
\label{MAINTH}
 Assume (H1)-(H4) and let $\Omega'\subset\subset\Omega$ with  $\Omega'$ open subset of $\R^n$.
Then, for any $g \in  C_0^{\infty}({\R})$, there exists an orthogonal projection $\Pi_g$ on  $L^2(\R^n;{\mathcal  H})$, an operator ${\mathcal  W} : L^2(\R^n ;{\mathcal  H})\rightarrow  L^2(\R^n)^{\oplus L}$, uniformly bounded with respect to $h$, and a selfadjoint $L\times L$ matrix $A$ of  $h$-admissible operators  $H^m(\R^n)\to L^2(\R^n)$, with the following properties:
\begin{itemize}
\item For all $\i \in C_0^\infty (\Omega')$,
$$
\Pi_g \i =\Pi_0 \i+ {\mathcal  O}(h);
$$
\item ${\mathcal  W}{\mathcal  W}^* = 1$ and ${\mathcal  W}^*{\mathcal  W}=\Pi_g$;
\item For $x\in \Omega'$, the  symbol $a(x,\xi ;h )$ of $A$ verifies,
$$
a(x,\xi ;h) = \omega (x,\xi ;h){\bf I}_L+ {\mathcal  M}(x) + W(x){\bf I}_L+hr(x,\xi ;h )
$$
where ${\bf I}_L$ stands for the $L$-dimensional identity matrix,  ${\mathcal  M}(x)$ is a $L\times L$ matrix depending smoothly on $x\in \Omega'$ and admitting $\lambda_{L'+1}(x)$, ..., $\lambda_{L'+L}(x)$ as eigenvalues,
and where
$
\partial^\alpha r (x,\xi ;h ) ={\mathcal  O}(\langle \xi\rangle^{m-1})
$
for any multi-index $\alpha$ and uniformly with respect to $(x,\xi )\in \Omega'\times\R^{n}$ and $h>0$ small enough; 
\item For any $f\in C_0^\infty (\R)$ with $\supp f\subset \{ g=1\}$, and for any $\varphi_0\in L^2(\R^n;{\mathcal  H})$
such that $\Vert \varphi_0\Vert =1$, and,
\be
\label{condloc}
\Vert\varphi_0\Vert_{L^2(K_0^c;{\mathcal  H})}+\Vert (1-\Pi_g )\varphi_0\Vert +\Vert (1-f(P))\varphi_0\Vert ={\mathcal  O}(h^\infty ),
\ee
for some $K_0\subset\subset\Omega'$,
one has,
\begin{equation}
\label{sol}
e^{-itP/h}\varphi_0= {\mathcal  W}^*e^{-itA/h}{\mathcal  W}\varphi_0  +{\mathcal  O}\left( \la t\ra h^\infty  \right) 
\end{equation}
uniformly with respect to $h>0$ small enough and $t\in  [0, T_{\Omega'}(\varphi_0))$.
\end{itemize}
\end{theorem}
\begin{remark}\sl  Actually, much more informations are obtained on the operators $\Pi_g$, ${\mathcal  W}$ and $A$, and we refer to Theorems \ref{th2} and \ref{th3} for more details, and to Chapter \ref{expexpl} for an explicit computation of $A$, up to  ${\mathcal  O}(h^4)$.
\end{remark}
\begin{remark}\sl  Condition (\ref{condloc}) on the initial data may seems rather strong, but in fact, it will become clear from the proof that the operators $\Pi_g$, $f(\tilde P)$ and $\i$ (where $\i \in C_0^\infty (\R^n)$ is supported in $K_0$) essentially commutes two by two (up to ${\mathcal  O}(h))$. Indeed, in the case of a molecular system, they respectively correspond to a localization in space for the nuclei, a localization in energy for the electrons, and a localization in energy for the whole molecule.
\end{remark}
\begin{remark}\sl  Here, we have assumed that both $\Pi_0^-(x)$ and $\Pi_0(x)$ have finite rank, since this corresponds to the main  applications that we have in mind. However, it will become clear from the proof that the case where one or both of them have infinite rank could  be treated in a similar way, with the difference that, if ${\rm Rank}\hskip 1pt \Pi_0(x) =\infty$, then  ${\mathcal  W}^*e^{-itA/h}{\mathcal  W}$ must be replaced by $e^{-it\Pi_gP\Pi_g /h}$  (there will not be any operator $A$ anymore). Moreover, some assumption must be added in order to be able to construct a modified operator as in Chapter \ref{modop} (for instance, that both $\Pi_0^-(x)$ and $\Pi_0(x)$ admit convenient extensions to all $x\in\R^n$ that depend smoothly on  $x$ away from a neighborhood of $K$).
\end{remark}
\begin{remark}\sl
In the next chapter, we modify the operator $Q(x)$ away from the interesting region, in such a way that the new operator $\tilde Q(x)$ admits a {\it global} gap in its spectrum. With such an operator, a much better result can be obtained, that permits to decouple  the evolution in a somewhat more complete and abstract way: see Theorem \ref{th2} (especially (\ref{reduc1})). In particular, even if $\Vert (1-\Pi_g )\varphi_0\Vert$ is not small, Theorem \ref{th2} gives a description of the quantum evolution of $\varphi_0$ in terms of two independent reduced evolutions. 
\end{remark}
As a corollary, in the case $L=1$ we also obtain the following geometric lower bound on $T_{\Omega'}(\varphi_0)$, that relates it to the underlying classical Hamilton flow of the operator $A$:
\begin{corollary}\sl 
\label{CORT}
Assume moreover that $L=1$ and the coefficients $c_{\alpha}=c_\alpha (x;h)$ of $\boldsymbol{\omega}$ verify,
\be
\label{coeffprinc}
c_\alpha (x;h) =c_{\alpha ,0}(x) + \varepsilon (h)\tilde c_\alpha (x;h),
\ee
with $c_{\alpha ,0}$ real-valued and independent of $h$,  $\varepsilon (h)\to 0$ as $h\to 0$, and, for any $\beta$, $\vert\partial^\beta c_{\alpha,0} (x)\vert +\vert\partial^\beta \tilde c_\alpha (x,h)\vert ={\mathcal  O}(1)$ uniformly, and set,
$$
a_0(x,\xi ):= \sum_{|\alpha| \leq m} c_{\alpha,0} (x)\xi^\alpha + \lambda_{L'+1}(x) +W(x) \qquad (x\in \Omega').
$$
Also, denote by $H_{a_0}:= \partial_\xi a_0\partial_x -\partial_xa_0\partial_\xi$ the Hamilton field of  $a_0$.
Then, for any $f\in C_0^\infty (\R)$ with $\supp f\subset \{ g=1\}$, and for any $\varphi_0\in L^2(\R^n;{\mathcal  H})$
such that $\Vert \varphi_0\Vert =1$, and,
$$
\Vert\varphi_0\Vert_{L^2(K_0^c;{\mathcal  H})}+\Vert (1-\Pi_g )\varphi_0\Vert +\Vert (1-f(P))\varphi_0\Vert ={\mathcal  O}(h^\infty ),
$$
one has,
\be
\label{esttemps}
T_{\Omega'}(\varphi_0)\geq \sup\{ T>0\, ;\, \pi_x(\cup_{t\in [0,T]}\hskip 1pt \exp tH_{a_0}(K(f)))\subset \Omega'\},
\ee
where  $\pi_x$ stands for the projection $(x,\xi )\mapsto x$, and $K(f)$ is the compact subset of $\R^{2n}$ defined by,
$$
K(f):=\{ (x,\xi)\,; \, x\in K_0,\, \omega(x,\xi ) +\gamma \leq C_f\}
$$
with $\gamma = \inf_{x\in\Omega'}\inf\sigma (Q(x))$ and $C_f:=\Max |\supp f|$.
\end{corollary}
\begin{remark}\sl 
Thanks to (H1) and (H2), it is easy to see that $\exp tH_{a_0}(x,\xi )$ is well defined for all $(t,x,\xi )\in\R\times\R^{2n}$.
\end{remark}
\begin{remark}\sl 
\label{esttemps2}
Actually, as it will be seen in the proof, in (\ref{esttemps}) one can replace the set $K(f)$ by $\cup_{j=1}^r FS(U_j\Pi_g\varphi_0)$, where $FS$ stands for the Frequency Set of locally $L^2$ functions introduced in \cite{GuSt} (we refer to Chapter \ref{proofcor} for more details).
\end{remark}
\begin{remark}\sl 
Our proof would permit to state a similar result in the case $L>1$, but under the additional assumption that the set $\{\lambda_{L'+1}(x),\dots,\lambda_{L'+L}(x)\}$ can be written as $\{ E_1(x),\dots ,E_{L''}(x)\}$, where the (possibly degenerate) eigenvalues $E_j(x)$ are such that $E_j(x) \not= E_{j'}(x)$ for $j\not= j'$ and $x\in\Omega$. In the general case where crossings may occur, such a type of result relies on the microlocal propagation of the Frequency Set for solutions of semiclassical matrix evolution problems (for which not much is known, in general).
\end{remark}
\begin{remark}\sl 
The proof also provides a very explicit and somehow optimal bound on $T_{\Omega'}(\varphi_0)$ in the case where $\varphi_0$ is a coherent state with respect to the $x$-variables: see Theorem \ref{WP} and (\ref{defthamilt}).
\end{remark}
\chapter{A Modified Operator}
\label{modop}
\setcounter{equation}{0}%
\setcounter{theorem}{0}%
In this chapter, we consider an arbitrary compact subset $K\subset\subset \Omega$ and an open neighborhood $\Omega_{K} \subset\subset \Omega$ of $K$. We also denote by $\Omega_0$  an open subset of $\R^n$, with closure disjoint from $\overline {\Omega_{K}}$,  and such that
$(\Omega_j)_{j=0}^r$ covers all of  $\R^n$, and we set $U_0:=\bf 1$.
The purpose of this chapter is to modify $Q(x)$ for $x$ outside a neighborhood of $K$, in order to make it regular with respect to $x$ there, and to deal with a global gap instead of a local one.
\vskip 0.3cm
Due to the contractibility of $\Omega$, we know that there exist $L'+L$ continuous functions $u_1,\dots, u_{L'+L}$ in $C(  \Omega ;{\mathcal  H})$, such that  the families $(u_1(x),\dots, u_{L'}(x))$ and $(u_{L'+1}(x),\dots, u_{L'+L}(x))$ span ${\rm Ran}\hskip 1pt  \Pi_0^-(x)$ and ${\rm Ran}\hskip 1pt  \Pi_0(x)$ respectively,  for all $x\in\Omega$ (see, e.g., \cite{KMSW}).
\vskip 0.3cm
Then, following Lemma 1.1 in \cite{KMSW}, we first prove,
\begin{lemma}\sl 
\label{sectionscinf}
For all $x\in {\R}^n$, there exist  $\tilde u_1(x),\dots, \tilde u_{L'+L}(x)$ in ${\mathcal  D}_Q$, such that the family $(\tilde u_1(x), \dots ,\tilde u_{L'+L}(x))$ is orthonormal in ${\mathcal  H}$ for all $x\in \R^n$,  the families $(\tilde u_1(x), \dots ,\tilde u_{L'}(x))$ and $(\tilde u_{L'+1}(x), \dots ,\tilde u_{L'+L}(x))$ span ${\rm Ran}\hskip 1pt  \Pi_0^-(x)$ and ${\rm Ran}\hskip 1pt  \Pi_0(x)$, respectively, when $x\in  \Omega_{K}$, and, for all $j=0,1,\cdots, r$ and $k=1,\dots, L'+L$,  
$$
U_j(x)\tilde u_k(x)\in  C^{\infty}(\Omega_j; {\mathcal  D}_Q).
$$
\end{lemma}
{\em Proof -- } Let $\zeta_1, \zeta_2\in  C^\infty (\R^n ; [0,1])$, such that $\supp\zeta_1\subset\Omega_0^c$, $\zeta_1 =1$ on  $\Omega_{K}$ and $\zeta_1^2 +\zeta_2^2 =1$ everywhere. Since $u_1(x), \dots ,u_{L'+L}(x)$ depend continuously on $x$ in $\Omega$, for any $\varepsilon >0$ one can find a finite number of points $x_1, \cdots, x_N\in  \supp\zeta_1$ and a partition of unity ${\i}_1, \cdots, {\i}_N\in  C_0^\infty (\Omega)$ on $\supp\zeta_1$, such that, for all $k=1,\dots, L'+L$,
$$
\sup_{x\in \supp\zeta_1}\Vert u_k(x) -\sum_{\ell=1}^N {\i}_\ell(x)u_k(x_\ell)\Vert_{{\mathcal  H}}\leq \varepsilon.
$$
On the other hand, using the last assertion of (H4), for any $(k,\ell )$ one can find $v_{k,\ell}$ in ${\mathcal  D}_Q$, such that, $\Vert v_{k,\ell}-u_k(x_\ell)\Vert_{{\mathcal  H}}\leq \varepsilon$ and $U_j(x)v_{k,\ell}\in  C^\infty (\Omega_j ,{\mathcal  D}_Q)$ for all $j=1,\dots ,r$. Moreover,
it follows from $(H3)$ and $(H4)$ that, for all $j=1,\cdots, r$,
$$
U_j(x)\Pi_0^-(x)U_j^*(x) \mbox{ and } U_j(x)\Pi_0(x)U_j^*(x)\in  C^{\infty}(\Omega_j,{\mathcal  L}({\mathcal  H}, {\mathcal  D}_Q)).
$$
Therefore, if we set, 
\begin{eqnarray*}
&& v_k(x):=\Pi_0^-(x)\sum_{\ell=1}^N {\i}_\ell(x)v_{k,\ell}\quad  (k=1,\dots,L');\\
&& v_k(x):=\Pi_0(x)\sum_{\ell=1}^N {\i}_\ell(x)v_{k,\ell}\quad  (k=L'+1,\dots,L'+L),
\end{eqnarray*}
and since $\sum_{\ell=1}^N {\i}_\ell(x) =1$ on  $\supp\zeta_1$, we obtain (also using that $\Pi_0^-(x)u_k(x)=u_k(x)$ for $k\leq L'$, and $\Pi_0(x)u_k(x)=u_k(x)$ for $k\geq L'+1$),
\begin{eqnarray*}
&&\sup_{x\in \supp\zeta_1}\Vert u_k(x) -v_k(x)\Vert_{{\mathcal  H}}\leq 2\varepsilon\\
&& U_j(x)v_k(x) \in  C^\infty (\Omega_j ,{\mathcal  D}_Q)\quad (j=1,\dots ,r).
\end{eqnarray*} 
In particular, by taking $\varepsilon$ small enough, we see that the families $(v_1(x), \dots, v_{L'}(x))$ and $(v_{L'+1}(x), \dots, v_{L'+L}(x))$ span ${\rm Ran}\hskip 1pt  \Pi_0^-(x)$ and ${\rm Ran}\hskip 1pt  \Pi_0(x)$, respectively,  for $x\in {\rm Supp}\hskip 1pt \zeta_1$.
Moreover, by Gram-Schmidt, this families can also be assumed to be orthonormal.
\vskip 0.3cm
Then, using again the last point of (H4),  one can find an orthonormal family $w_1, \dots, w_{L'+L} \in {\mathcal  D}_Q$, such that  $\vert\la w_m,u_k(x_\ell)\ra\vert\leq\varepsilon$ for all $1\leq k, m\leq {L'+L}$, $1\leq \ell\leq  N$, and $U_j(x)w_m\in  C^\infty (\Omega_j ,{\mathcal  D}_Q)$ ($j=1,\dots,r$). Thus, setting,
$$
\tilde w_k(x) := \zeta_1(x)v_k(x) + \zeta_2(x)w_k,
$$
we see that, for all $k,k'\in \{1,\dots,L'+L\}$,
$$
\la\tilde w_k(x), \tilde w_{k'}(x) \ra_{\mathcal  H} = \delta_{k,k'}+{\mathcal  O}(\varepsilon).
$$
As a consequence, taking $\varepsilon >0$ sufficiently small and orthonormalizing the family $(\tilde w_1(x), \dots, \tilde w_{L'+L}(x))$, we obtain a new family $(\tilde u_1(x), \dots, \tilde u_{L'+L}(x))$ that verifies all the properties required in the lemma. \hfill$\bullet$
\vskip 0.3cm
Then, (with the usual convention $\sum_{k=1}^{L'} =0$ if $L'=0$) we set,
\begin{eqnarray*}
&& \tilde \Pi_0^- (x) = \sum_{k=1}^{L'}\la \cdot ,\tilde u_k(x)\ra_{{\mathcal  H}}\tilde u_k(x),\\
&& \tilde \Pi_0 (x) = \sum_{k=L'+1}^{L'+L}\la \cdot ,\tilde u_k(x)\ra_{{\mathcal  H}}\tilde u_k(x)
\end{eqnarray*}
so that $\tilde \Pi_0^- (x)$ and $\tilde \Pi_0 (x)$ are orthogonal projections of  rank $L'$ and $L$ respectively, are orthogonal each other,   coincide with $\Pi_0^- (x)$ and $\Pi_0 (x)$ for $x$ in $\Omega_{K}$, and  verify,
\be
\label{regpi}
U_j(x)\tilde \Pi_0^-(x)U_j(x)^* \mbox{ and } U_j(x)\tilde \Pi_0 (x)U_j(x)^*\in  C^\infty (\Omega_j,{\mathcal  L}({\mathcal  H})),
\ee
for all $j=0, 1, \cdots,r$.
\vskip 0.3cm
Now, with the help of  $\tilde\Pi_0^-(x), \tilde\Pi_0(x)$,  we modify $Q(x)$ outside a neighborhood of $K$ as follows.
\begin{proposition}\sl  
\label{qtilde}
Let $\Omega_{K}'\subset\subset \Omega_{K}$ be an open neighborhood of $K$. Then, 
for all $x\in \R^n$, there exists a selfadjoint  operator $\tilde Q(x)$ on ${\mathcal  H}$, with domain ${\mathcal  D}_Q$, and uniformly semi-bounded from below, such that,
\begin{eqnarray}
&& \tilde Q (x) =Q(x)\quad \mbox{if}\,\,  x\in \Omega_{K}';\\
&&[\tilde Q (x), \tilde\Pi_0^-(x)] = [\tilde Q (x), \tilde\Pi_0(x)] =0\quad \mbox{for all }x\in \R^n,
\end{eqnarray}
and the application  $x\mapsto U_j(x)\tilde Q(x)U_j(x)^{-1}$ is in $C^\infty (\Omega_j; {\mathcal  L}({\mathcal  D}_Q, {\mathcal  H}))$  for all  $j=0,1,\cdots,r$. Moreover, the bottom of the spectrum of $\tilde Q(x)$ consists in $L'+L$  eigenvalues $\tilde\lambda_1(x), \dots, \tilde\lambda_{L'+L}(x)$, and 
$\tilde Q(x)$
admits a global gap in its spectrum, in the sense that,
$$
\inf_{x\in \R^n}{\dist}(\sigma (\tilde Q(x))\backslash \{\tilde\lambda_{L'+1}(x), \dots, \tilde\lambda_{L'+L}(x)\}, \{\tilde\lambda_{L'+1}(x), \dots, \tilde\lambda_{L'+L}(x)\}) >0.
$$
\end{proposition}
{\bf Proof } 
We set $\tilde\Pi_0^+(x)=1-\tilde\Pi_0^-(x)-\tilde\Pi_0(x)$ and we  choose a function  $\zeta\in  C_0^{\infty}(\Omega_{K} ; [0,1])$ such that $\zeta=1$ on $\Omega_{K}'$. Then, with $Q_0$ as in (H4), we set,
$$
\tilde Q(x)=\zeta(x)Q(x)+(1-\zeta(x))\tilde\Pi_0^+(x) Q_0\tilde\Pi_0^+(x) - (1-\zeta (x))\tilde\Pi_0^-(x).
$$
Since $\tilde\Pi_0^-(x)=\Pi_0^-(x)$ and $\tilde\Pi_0(x)=\Pi_0(x)$ on ${\rm Supp}\hskip 1pt \zeta$, we see that $\tilde\Pi_0^-(x)$ and $\tilde\Pi_0(x)$ commute with $\tilde Q(x)$, and it is also clear that $\tilde Q(x)$ is selfadjoint with domain ${\mathcal  D}_Q$.
Moreover, 
\begin{eqnarray*}
&& \tilde\Pi_0^-(x)\tilde Q(x)\tilde\Pi_0^-(x)= \zeta (x)\Pi_0^-(x) Q(x)\Pi_0^-(x) - (1-\zeta (x))\Pi_0^-(x);\\
&& \tilde\Pi_0(x)\tilde Q(x)\tilde\Pi_0(x)= \zeta (x)\Pi_0(x) Q(x)\Pi_0(x),
\end{eqnarray*}
and, setting,
$$
\lambda_{L+L'+1}(x):=\inf \left(\sigma (Q(x))\backslash \{\lambda_1(x), \dots, \lambda_{L+L'}(x)\}\right),
$$
one has,
$$
\tilde \Pi_0^+(x)\tilde Q(x)\tilde \Pi_0^+(x)\geq \left(\zeta (x)\lambda_{L+L'+1}(x) + (1-\zeta (x)\right)\tilde\Pi_0^+(x).
$$
In particular, the bottom of the spectrum of $\tilde Q(x)$ consists in the $L+L'$ eigenvalues $\tilde\lambda_k(x) =\zeta (x)\lambda_k(x) - (1-\zeta (x))$ ($k=1,\dots,L'$), $\tilde\lambda_k(x) =\zeta (x)\lambda_k(x)$ ($k=L'+1,\dots,L'+L$), and, due to (H3), one has,
$$
\inf_{x\in\R^n}\left(\tilde\lambda_{L'+1}(x) - \tilde\lambda_{L'}(x) \right)= \inf_{x\in\R^n}\left(\zeta (x)(\lambda_{L'+1}(x) - \lambda_{L'}(x)+ (1-\zeta (x))\right)>0,
$$
and 
\begin{eqnarray*}
&&\inf_{x\in \Omega}{\dist}(\sigma (\tilde Q(x))\backslash \{\tilde\lambda_1(x), \dots, \tilde\lambda_{L'+L}(x)\}, \{\tilde\lambda_1(x), \dots, \tilde\lambda_{L'+L}(x)\})\\
&&\hskip 3cm \geq
\inf_{x\in \Omega}\vert \zeta (x)(\lambda_{L'+L+1}(x)-\lambda_{L'+L}(x)) + (1-\zeta (x))\vert >0,
\end{eqnarray*}
while, since $\supp\zeta\subset\Omega$,
$$
\inf_{x\in \R^n\backslash\Omega}{\dist}(\sigma (\tilde Q(x))\backslash\{\tilde\lambda_1(x), \dots, \tilde\lambda_{L'+L}(x)\}, \{\tilde\lambda_1(x), \dots, \tilde\lambda_{L'+L}(x)\})\geq 1.
$$
In particular, $\tilde Q(x)$
admits a fix global gap in its spectrum as stated in the proposition.
Finally, using (H4) and (\ref{regpi}), we see that 
$U_j(x)\tilde Q(x)U_j^*(x)$ depends smoothly on $x$ in $\Omega_j$ for all $j=0,1,\cdots,r$.\hfill$\bullet$
\vskip 0.3cm
In the sequel, we also set,
\be
\label{ptilde}
\tilde P =\boldsymbol{\omega}+\boldsymbol{Q}:=\boldsymbol{\omega}+\tilde Q(x)+\zeta (x)W(x),
\ee 
and we denote by $\tilde \Pi_0$ the projection on $L^2({\R}^{n};{\mathcal  H})$ induced by the action of 
$\tilde \Pi_0(x)$ on ${\mathcal  H}$, i.e. the unique projection on $L^2({\R}^{n};{\mathcal  H})$ that verifies 
$$
\tilde\Pi_0(f\otimes g)(x) = f(x)\tilde\Pi_0(x)g \quad (\mbox{a.e. on }  \R^n\ni x)
$$
 for all $f\in  L^2(\R^n)$ and $g\in  {\mathcal  H}$.
 \chapter{Twisted $h$-Admissible Operators}
 \label{twao}
\setcounter{equation}{0}%
\setcounter{theorem}{0}%
In  order to construct  (in the same spirit as in \cite{BrNo, HeSj12, MaSo, NeSo, Sj2, So})  an orthogonal projection $\Pi$ on $L^2({\R}^{n};{\mathcal  H})$ such that 
$\Pi-\Pi_0={\mathcal  O}(h)$
and  $[\tilde P,\Pi]={\mathcal  O}(h^\infty)$ (locally uniformly in energy), we need to generalize the notion of $h$-admissible operator with operator-valued symbol (see, e.g., \cite{Ba, GMS} and the Appendix) by taking into account the possible singularities of $Q(x)$. To avoid complications, in this chapter we also restrict our attention to the case of bounded operators. The case of unbounded ones will be considered in the next chapter, at least from the point of view of {\it differential} operators.
\begin{definition} \sl
\label{regcov}
We call ``regular covering'' of $\R^n$ any finite  family $(\Omega_j)_{j=0,\cdots,r}$ of open subsets of $\R^n$ such that $\cup_{j=0}^r\Omega_j =\R^n$ and such that there exists a family of  functions 
${\i}_{j}\in  C_b^{\infty}(\R^n)$ (the space of smooth functions on $\R^n$ with uniformly bounded derivatives of all order) with  $\sum_{j=0}^r{\i}_j =1$, $0\leq \i_j\leq 1$, and ${\rm dist} \left( \supp ({\i}_{j}), \R^n\backslash \Omega_j\right) >0$ (${j=0,\cdots,r}$). Moreover, if $U_j(x)$ ($x\in \Omega_j, \, 0\leq j\leq r$)   is a family of unitary operators on ${\mathcal  H}$, the family $\left( U_j,\Omega_j\right)_{j=0,\cdots,r}$ (where $U_j$ denotes the unitary operator on $L^2(\Omega_j;{\mathcal  H})\simeq L^2(\Omega_j)\otimes{\mathcal  H}$ induced by the action of $U_j(x)$ on ${\mathcal  H}$) will be called a ``regular unitary covering'' of $L^2(\R^n;{\mathcal  H})$. 
\end{definition}
\begin{remark}\sl 
Despite the terminology that we use, no assumption is made on any possible regularity of $U_j(x)$ with respect to $x$. 
\end{remark}
\begin{remark}\sl  Possibly by shrinking a little bit $\Omega$ around the compact set $K$, one can always assume that the family $(U_j,\Omega_j)_{j=0,1,\cdots,r}$ defined in Chapter \ref{sect2} is a regular unitary covering of $L^2(\R^n;{\mathcal  H})$.
\end{remark}
In the sequels, we denote by $C_d^\infty (\Omega_j)$ the space of functions ${\i}\in  C_b^{\infty}(\R^n)$ such that  ${\rm dist} \left( \supp ({\i}), \R^n\backslash \Omega_j\right) >0$.
\begin{definition}[Twisted $h$-Admissible Operator]\sl
\label{utruc}
Let  ${\mathcal  U}:=\left( U_j,\Omega_j\right)_{j=0,\cdots,r}$ be a regular unitary covering (in the previous sense) of $L^2(\R^n;{\mathcal  H})$.
We say that an operator $A:L^2({\R^n}; {\mathcal  H})\rightarrow L^{2}({\R^n};{\mathcal  H})$ is a ${\mathcal  U}$-twisted  $h$-admissible operator, if there exists  a family of  functions 
${\i}_{j}\in  C_d^{\infty}(\Omega_j)$  such that,  for any $N\geq 1$, $A$ can be
written in the form,
\be
\label{adapt}
A=\sum_{j=0}^rU_j^{-1}{\i}_j A^N_jU_j{\i}_j + R_N,
\ee
where $\Vert R_N\Vert_{{\mathcal  L}(L^2(\R^n;{\mathcal  H}))}=  {\mathcal  O}(h^N)$, and, for any $j=0, ..,r$, 
$A^N_j$ is a bounded $h$-admissible operator on  $L^{2}({\R^n};{\mathcal  H})$ with  symbol $a^N_j(x,\xi)\in 
C_b^{\infty}(T^*{\R}^n; {\mathcal  L}({\mathcal  H}))$,
and, for any $\varphi_\ell\in C_d^\infty (\Omega_\ell)$ ($\ell=0,..,r$), the operator
$$
U_\ell\varphi_\ell U_j^{-1}{\i}_{j}A^N_j{\i}_{j}U_jU_\ell^{-1}\varphi_\ell,
$$
is still  an $h$-admissible operator on  $L^{2}({\R^n};{\mathcal  H})$.\\
\end{definition}
\begin{remark}\sl 
In particular, by the Calder\'on-Vaillancourt theorem, the norm of $A$ on $L^{2}({\R^n};{\mathcal  H})$ is bounded uniformly with respect to $h\in(0,1]$.
\end{remark}
An equivalent definition is given by the following proposition:
\begin{proposition}\sl 
\label{eqdef}
An operator $A: L^{2}({\R^n};{\mathcal  H})\rightarrow L^{2}({\R^n};{\mathcal  H})$ is a ${\mathcal  U}$-twisted $h$-admissible operator if and only if the two following properties are verified:
\begin{enumerate}
\item For any $N\geq 1$ and any  functions ${\i}_1,\cdots,{\i}_N\in  C_b^\infty (\R^n)$, one has,
$$
{\rm ad}_{{\i}_1}\circ\cdots\circ {\rm ad}_{{\i}_N}(A) = {\mathcal  O}(h^N)\;\; : \,\, L^{2}({\R^n};{\mathcal  H})\rightarrow L^{2}({\R^n};{\mathcal  H})
$$
where we have used the notation ${\rm ad}_{{\i}}(A):=[{\i} ,A] = {\i} A- A{\i}$.
\item For any $\varphi_j \in  C_d^\infty (\Omega_j)$,  $U_j \varphi_j A U_j^{-1}\varphi_j$ is a bounded $h$-admissible operator on  $L^{2}({\R^n};{\mathcal  H})$.
\end{enumerate}
\end{proposition}
{\em Proof -- } From Definition \ref{utruc}, it is clear that any ${\mathcal  U}$-twisted $h$-admissible operator  verifies the properties of the Proposition. Conversely, assume $A$ verifies these properties, and denote by $({\i}_j)_{j=0,\cdots,r}\subset C_b^\infty (\R^n)$ a partition of unity on $\R^n$ such that ${\rm dist} \left( \supp ({\i}_{j}), \R^n\backslash \Omega_j\right) >0$. Then, for all $j$ one can construct $\varphi_j , \psi_j \in  C_d^\infty (\Omega_j)$,  such that   $\varphi_j {\i}_j ={\i}_j$ and  $\psi_j \varphi_j =\varphi_j$, and, for any $N\geq 1$, we can write,
\begin{eqnarray*}
A &=& \sum_{j=0}^r {\i}_jA = \sum_{j=0}^r \left( {\i}_j A\varphi_j + {\i}_j \mbox{ad}_{\varphi_j}(A)\right)\\
&=&\sum_{j=0}^r \left( {\i}_j A\varphi_j+ {\i}_j \mbox{ad}_{\varphi_j}(A)\varphi_j+ {\i}_j \mbox{ad}_{\varphi_j}^2(A)\right)\\
&=& \cdots =\sum_{j=0}^r \left( \sum_{k=0}^{N-1} {\i}_j \mbox{ad}_{\varphi_j}^k(A)\varphi_j+ {\i}_j \mbox{ad}_{\varphi_j}^{N}(A)\right)\\
&=& \sum_{j=0}^r \left( \sum_{k=0}^{N-1} \psi_j {\i}_j \mbox{ad}_{\varphi_j}^k(A)\varphi_j\psi_j+ {\i}_j \mbox{ad}_{\varphi_j}^{N}(A)\right).
\end{eqnarray*}
In particular, since $\mbox{ad}_{\varphi_j}^{N}(A) ={\mathcal  O}(h^N)$, and $U_j$ commutes with the multiplication by functions of $x$, we obtain
\be
\label{reconst1}
A= \sum_{j=0}^r U_j^{-1}\psi_j A^N_j U_j\psi_j+ {\mathcal  O}(h^N)
\ee
with
\be
\label{reconst2}
A^N_j:= \sum_{k=0}^{N-1}U_j {\i}_j\mbox{ad}_{\varphi_j}^k(A) U_j^{-1}\varphi_j = \sum_{k=0}^{N-1} {\i}_j\mbox{ad}_{\varphi_j}^k(U_j\varphi_jAU_j^{-1}\varphi_j).
\ee
Therefore, $A^N_j$ is a bounded $h$-admissible operator, and for any $\tilde \psi_l\in C_d^\infty (\Omega_l)$, it verifies,
$$
U_l\tilde\psi_lU_j^{-1}\psi_jA^N_j\psi_jU_j\tilde\psi_lU_l^{-1} = \sum_{k=0}^{N-1} {\i}_j\mbox{ad}_{\varphi_j}^k(U_l\tilde\psi_lAU_l^{-1}\tilde\psi_l)\varphi_j,
$$
that is still an $h$-admissible operator. Thus, the proposition follows.\hfill$\bullet$
\vskip 0.3cm
In the sequel, if $A$ is a ${\mathcal  U}$-twisted $h$-admissible operator, then an expression of $A$ as in (\ref{adapt}) will be said ``adapted'' to ${\mathcal  U}$.
\vskip 0.3cm
One also has at disposal a notion of (full) symbol for such operators. In the sequels, we denote by $S(\Omega_j\times \R^n;{\mathcal  L}({\mathcal  H}))$ the space of ($h$-dependent) operator-valued symbols $a_j\in C^\infty (\Omega_j\times \R^n;{\mathcal  L}({\mathcal  H}))$ such that, for any $\alpha\in\Z_+^{2n}$, the quantity $\Vert \partial^\alpha a_j(x,\xi )\Vert_{{\mathcal  L}({\mathcal  H})}$ is bounded uniformly for $h$ small enough and for $(x,\xi)$ in any set of the form $\Omega_j'\times\R^n$, with $\Omega_j'\subset\Omega_j$, ${\rm dist} \left( \Omega_j', \R^n\backslash \Omega_j\right) >0$. We also set,
\begin{eqnarray*}
&& \Omega :=(\Omega_j)_{j=0, \dots, r};\\
&& {\bf S}(\Omega ; {\mathcal  L}({\mathcal  H})):= S(\Omega_0\times \R^n;{\mathcal  L}({\mathcal  H}))\times \dots\times S(\Omega_r\times \R^n;{\mathcal  L}({\mathcal  H})),
\end{eqnarray*}
and we write $a={\mathcal  O}(h^\infty )$ in ${\bf S}(\Omega ; {\mathcal  L}({\mathcal  H}))$ when $\Vert \partial^\alpha a_j(x,\xi)\Vert_{{\mathcal  L}({\mathcal  H})} ={\mathcal  O}(h^\infty )$ uniformly in any set $\Omega_j'\times\R^n$ as before.
\begin{lemma}\sl 
Let $A$ be a ${\mathcal  U}$-twisted $h$-admissible operator, where\\ ${\mathcal  U}=(U_j,\Omega_j)_{0\leq j\leq r}$ is some regular unitary covering. Then, for all $j=0,\dots,r$, there exists an operator-valued symbol $a_j\in S(\Omega_j\times \R^n;{\mathcal  L}({\mathcal  H}))$, unique up to ${\mathcal  O}(h^\infty )$, such that, for any $\i_j=\i_j(x)\in C_d^\infty (\Omega_j)$, the  symbol of the $h$-admissible operator $U_j {\i}_j A U_j^{-1}{\i}_j$ is ${\i}_j\sharp a_j\sharp\i_j$ (where $\sharp$ stands for the standard symbolic composition: see Appendix A). 
\end{lemma}
{\em Proof -- } Indeed, given two  functions $\i_j, \varphi_j\in C_d^\infty (\Omega_j)$ with $\varphi_j\i_j =\i_j$, one has
$$
U_j {\i}_j A U_j^{-1}{\i}_j={\i}_j\left( U_j\varphi_j A U_j^{-1}\varphi_j\right){\i}_j,
$$
and thus, denoting by $a_j^{\chi}$ the symbol of $U_j {\i} A U_j^{-1}{\i}$, one obtains
$$
a_j^{\chi_j}= \i_j\sharp a_j^{\varphi_j}\sharp \i_j.
$$
In particular, using the explicit expression of $\sharp$ (see Appendix A, Proposition \ref{sharp}), we see that $a_j^{\varphi_j}=a_j^{\chi_j}  + {\mathcal  O}(h^\infty )$ in the interior of $\{ \chi_j (x)=1\}$.
Then, the result follows by taking a non-decreasing sequence $(\varphi_{j,k})_{k\geq 1}$ in $C_d^\infty (\Omega_j)$, such that $\bigcup_{k\geq 0}\{x\in\Omega_j\, ;\,  \varphi_{j,k}(x)=1\} =\Omega_j$, and, for any $(x,\xi )\in \Omega_j\times\R^n$,  by defining $a_j(x,\xi)$ as the common value of the $a_j^{\varphi_{j,k}}(x,\xi)$'s for $k$ large enough.\hfill$\bullet$
\vskip 0.3cm
\begin{definition}[Symbol]\sl
Let $A$ be a ${\mathcal  U}$-twisted $h$-admissible operator, where\\ ${\mathcal  U}=(U_j,\Omega_j)_{0\leq j\leq r}$ is some regular unitary covering. Then, the family of operator-valued functions $\sigma(A):=(a_j)_{0\leq j\leq r}\in  {\bf S}(\Omega ; {\mathcal  L}({\mathcal  H}))$, defined in the previous lemma, is called the (full) symbol of $A$. Moreover,
$A$ is said to be elliptic if, for any $j=0,\cdots,r$ and $(x,\xi )\in\Omega_j\times\R^n$, the operator $a_j(x,\xi)$ is invertible on ${\mathcal  H}$, and verifies,
\be
\label{ell}
\Vert a_j(x,\xi )^{-1}\Vert_{{\mathcal  L}({\mathcal  H})} ={\mathcal  O}(1),
\ee
uniformly for $h$ small enough and for $(x,\xi)$ in any set of the form $\Omega_j'\times\R^n$, with $\Omega_j'\subset\Omega_j$, ${\rm dist} \left( \Omega_j', \R^n\backslash \Omega_j\right) >0$.  \end{definition}
In particular, it follows from the proof of Proposition \ref{eqdef} that, if such an operator $A$ is elliptic, then it can be written in the form (\ref{adapt}), with $A_j^N$ elliptic on $\{ {\i}_j\not= 0\}$ for all $j,N$. Moreover, we have the two following result on composition and parametrices:
\begin{proposition}[Composition]\sl 
 Let ${\mathcal  U}$ be a regular covering of  $L^{2}({\R^n};{\mathcal  H})$, and let
$A$ and $B$ be two ${\mathcal  U}$-twisted $h$-admissible operators. Then, the composition $AB$ is  a ${\mathcal  U}$-twisted $h$-admissible operator, too. Moreover, its  symbol is given by,
$$
\sigma(AB) =\sigma(A)\sharp \sigma(B),
$$
where the operation $\sharp$ is defined component by component, that is,
$$
(a_j)_{0\leq j\leq r}\sharp (b_j)_{0\leq j\leq r}:= (a_j\sharp b_j)_{0\leq j\leq r}.
$$
\end{proposition}
{\em Proof -- } First of all, since
 $$
{\rm ad}_{\i} (AB) = {\rm ad}_{\i} (A)B + A\hskip 1pt {\rm ad}_{\i} (B),
$$
one easily sees, by induction on $N$, that the first condition in Proposition \ref{eqdef} is satisfied. Moreover, if ${\i}_j \in  C_d^\infty (\Omega_j)$,  let $\varphi_j  \in  C_d^\infty (\Omega_j)$ such that  $\varphi_j{\i}_j ={\i}_j$. Then, if, for any operator $C$, we set $C_j:=U_j\varphi_j CU_j^{-1}\varphi_j$, we have,
\begin{eqnarray}
U_j{\i}_jABU_j^{-1}{\i}_j &=& {\i}_j A_j B_j {\i}_j + U_j{\i}_j \mbox{ad}_{(\varphi_j^2)}(A) BU_j^{-1}{\i}_j\nonumber\\
&=& {\i}_j A_j B_j {\i}_j + {\i}_j [\mbox{ad}_{(\varphi_j^2)}(A)]_j B_j{\i}_j + U_j{\i}_j \mbox{ad}^2_{(\varphi_j^2)}(A) BU_j^{-1}{\i}_j\nonumber\\
&=& \cdots \nonumber\\
\label{comp2}
&=& \sum_{k=0}^{N-1}{\i}_j [\mbox{ad}^k_{(\varphi_j^2)}(A)]_j B_j{\i}_j + U_j{\i}_j \mbox{ad}^N_{(\varphi_j^2)}(A) BU_j^{-1}{\i}_j
\end{eqnarray}
for all $N\geq 1$. Therefore, since $U_j{\i}_j \mbox{ad}^N_{(\varphi_j^2)}(A) BU_j^{-1}{\i}_j
= {\mathcal  O}(h^N)$, and the operator $[\mbox{ad}^k_{(\varphi_j^2)}(A)]_j = \mbox{ad}^k_{(\varphi_j^2)}(A_j)$ is a bounded $h$-admissible operator, we deduce from (\ref{comp2}) that $AB$ is a ${\mathcal  U}$-twisted $h$-admissible operator. Moreover, since $\varphi_j =1$ on the support of $\i_j$, we see that the symbol of ${\i}_j \mbox{ad}^k_{(\varphi_j^2)}(A_j)$ vanishes identically for $k\geq 1$, and thus, we also deduce from (\ref{comp2}) that the symbol $(c_j)_{0\leq j\leq r}$ of $AB$  verifies,
$$
\i_j\sharp c_j\sharp \i_j=\i_j\sharp a_j\sharp b_j\sharp \i_j, 
$$
for any $\i_j \in  C_d^\infty (\Omega_j)$, and the result follows.\hfill$\bullet$
\vskip 0.3cm
\begin{proposition}[Parametrix]\sl 
\label{param}
Let $A$ be a ${\mathcal  U}$-twisted $h$-admissible operator, and assume that $A$ is elliptic. Then, $A$ is invertible on $L^2(\R^n ;{\mathcal  H})$, and its inverse $A^{-1}$ is a ${\mathcal  U}$-twisted $h$-admissible operator. Moreover, its  symbol $\sigma (A^{-1})$ is related to the one $\sigma(A)=(a_j)_{0\leq j\leq r}$ of $A$ by the following formula:
$$
\sigma (A^{-1})=(\sigma(A))^{-1} + hb,
$$
where $(\sigma(A))^{-1}:=(a_j^{-1})_{0\leq j\leq r}$ and $b\in  {\bf S}(\Omega ; {\mathcal  L}({\mathcal  H}))$.
\end{proposition}
{\em Proof -- } We first prove that $A$ is invertible by following an idea of \cite{KMSW} (proof of Theorem 1.2). 
\vskip 0.3cm
For ${j=0,\cdots,r}$, let ${\i}_{j}, \varphi_j\in  C_d^{\infty}(\Omega_j)$ such that  $\varphi_j\i_j =\i_j$, and  $\sum_{j=0}^r{\i}_j =1$. Then, by assumption, the symbol of $U_j \varphi_j A U_j^{-1}\varphi_j$ can be written on the form $\varphi_j(x)\sharp a_j(x,\xi )\sharp \varphi_j(x)$ with $a_j(x,\xi)$ invertible, and the operator,
$$
B := \sum_{j=0}^r U_j^{-1}\varphi_j^3 {\rm Op}_h(\varphi_j a_j^{-1})U_j{\i}_j
$$
is well defined and bounded on $L^2(\R^n;{\mathcal  H})$. Moreover, using the standard symbolic calculus, we compute,
\begin{eqnarray}
AB &=&\sum_{j=0}^r AU_j^{-1}\varphi_j^3 {\rm Op}_h( \varphi_ja_j^{-1})U_j{\i}_j\nonumber\\
&=& \sum_{j=0}^r U_j^{-1}\varphi_jU_j\varphi_jAU_j^{-1}\varphi_j {\rm Op}_h( \varphi_ja_j^{-1})U_j{\i}_j \nonumber\\
&&\hskip 4cm+ [A, \varphi_j^2]U_j^{-1}\varphi_j {\rm Op}_h( \varphi_ja_j^{-1})U_j{\i}_j\nonumber\\
&=& \sum_{j=0}^rU_j^{-1}\varphi_j{\rm Op}_h(\varphi_j^2a_j){\rm Op}_h(\varphi_j a_j^{-1})U_j{\i}_j +{\mathcal  O}(h)\nonumber\\
\label{invapprox}
&=& \sum_{j=0}^rU_j^{-1}\varphi_j^4U_j{\i}_j+{\mathcal  O}(h) = \sum_{j=0}^r{\i}_j+{\mathcal  O}(h)=1+{\mathcal  O}(h).
\end{eqnarray}
In the same way, defining,
$$
B':=\sum_{j=0}^r U_j^{-1}{\i}_j {\rm Op}_h(\varphi_ja_j^{-1})U_j\varphi_j^3,
$$
we obtain $B'A = 1+{\mathcal  O}(h)$, and this proves the invertibility of $A$ for $h$ small enough. It remains to verify that $A^{-1}$ is  a ${\mathcal  U}$-twisted  $h$-admissible operator. We first prove,
\begin{lemma}\sl 
\label{suppdisj}
Let $A$ be a ${\mathcal  U}$-twisted $h$-admissible operator, and let ${\i} , \psi\in  C_b^\infty (\R^n)$ such that ${\rm dist} \left( \supp {\i}, \supp\psi \right) >0$. Then, $\Vert {\i} A\psi\Vert = {\mathcal  O}(h^\infty )$.
\end{lemma}
{\em Proof -- } Given $N\geq 1$, let $\varphi_1, \cdots, \varphi_N\in  C_b^\infty (\R^n)$, such that $\varphi_1{\i} ={\i}$, $\varphi_{k+1}\varphi_k = \varphi_k$ ($k=1,\cdots, N-1$), and $\varphi_N\psi =0$. Then, one has,
\begin{eqnarray*}
{\i} A \psi &=& \varphi_1 {\rm ad}_{\i} (A)\psi = \varphi_2{\rm ad}_{\varphi_1} \circ {\rm ad}_{\i} (A)\psi\\
&=& \cdots = {\rm ad}_{\varphi_N} \circ\cdots \circ{\rm ad}_{\varphi_1}\circ{\rm ad}_{\i} (A)\psi ={\mathcal  O}(h^{N+1}).
\end{eqnarray*}
{} \hfill$\bullet$
\vskip 0.3cm
Now, since,
$$
{\rm ad}_{\i} (A^{-1}) = -A^{-1}{\rm ad}_{\i} (A)A^{-1},
$$
it is easy to see, by induction on $N$, that $A^{-1}$ satisfies to the first property of Proposition \ref{eqdef}. Moreover, for $v\in  L^2(\R^n;{\mathcal  H})$ and for ${\i}_{j}\in  C_d^\infty (\Omega_j)$, let us set,
$$
u= A^{-1}U_j^{-1}{\i}_j v,
$$
and choose $\varphi_{j}, \in  C_d^{\infty}(\Omega_j;\R)$, $\psi_j\in C_b^\infty (\R^n;\R)$, such that  $\psi_j{\i}_j =0$, $\varphi_j^4+\psi_j^2\geq 1$, and ${\rm dist} \left( \supp (\varphi_{j}-1), \supp {\i}_j\right) >0$. Then, since the symbol of $A_j:=U_j \varphi_j A U_j^{-1}\varphi_j$ is of the form $\varphi_j\sharp a_j\sharp\varphi_j$ with $a_j(x,\xi)$ invertible for $x$ in $\supp\varphi_j$, we see that the bounded $h$-admissible operator $B_j:=A_j^*A_j + \psi_j^2$ is globally elliptic, and one has,
\begin{eqnarray}
B_jU_j{\i}_ju&=&A_j^*A_jU_j{\i}_ju=A_j^*U_j\varphi_j A{\i}_ju =A_j^*U_j{\i}_j Au + A_j^*U_j\varphi_j [A,{\i}_j]u\nonumber\\
&=&  A_j^*{\i}_j^2v+A_j^*U_j\varphi_j [A,{\i}_j]\varphi_j^2u + A_j^*U_j{\i}_j A (\varphi_j^2 -1)u\nonumber\\
\label{pseudinv}
&=& A_j^*{\i}_j^2v+A_j^* [A_j,{\i}_j]U_j\varphi_ju + {\mathcal  O}(h^\infty\Vert v\Vert),
\end{eqnarray}
where the last estimate comes from Lemma \ref{suppdisj}. In particular, since $B_j^{-1}$ is an $h$-admissible operator,  we obtain that $U_j{\i}_ju$ can be written on the form,
$$
U_j{\i}_ju = C_j v + hC'_jU_j\varphi_ju + {\mathcal  O}(h^\infty\Vert v\Vert)
$$
where $C_j, C_j'$ are bounded $h$-admissible operators. Repeating the same argument with  $U_j\varphi_ju$ instead of $U_j{\i}_ju$, and iterating the procedure, it easily follows that $U_j{\i}_jA^{-1}U_j^{-1}{\i}_j$ is an $h$-admissible operator. Moreover,  we see on (\ref{pseudinv}) that the symbol of $U_j\i_jA^{-1}U_j^{-1}\i_j$ coincides, up to ${\mathcal  O}(h)$, with that of $B_j^{-1}A_j^*\i_j^2$,  that is,
$$
(\varphi_j^4(x) a_j^*(x,\xi)a_j(x,\xi) + \psi_j^2(x))^{-1}a_j^*(x,\xi)\i_j(x)^2 = a_j(x,\xi)^{-1}\i_j(x)^2,
$$
since $\varphi_j =1$ and $\psi_j =0$ on the support of $\i_j$.
Thus, the proposition follows.\hfill$\bullet$
\begin{proposition}[Functional Calculus]\sl 
\label{fc}
Let $A$ be a selfadjoint ${\mathcal  U}$-twisted $h$-admissible operator, and let $f\in C_0^\infty (\R)$. Then, the operator $f(A)$ is a ${\mathcal  U}$-twisted $h$-admissible operator, and its symbol is related to that of $A$ by the formula,
$$
\sigma (f(A)) = f(\re \sigma (A)) + hb,
$$
where $f(\re (a_j)_{j=0,\dots ,r}):= (f(\re a_j))_{j=0,\dots ,r}$, $\re a_j:=\frac12(a_j +a_j^*)$,  and $b\in {\bf S}(\Omega ;{\mathcal  L}(H))$.
\end{proposition}
{\em Proof -- } We use a formula of representation of $f(A)$ due to B. Helffer and J. Sj\"ostrand. Denote by $\tilde f\in 
C^{\infty}_0({\C})$ an almost analytic extension of $f$, that is, a function verifying $\tilde f\left\vert_{\R} = f\right.$ and 
$\vert{\overline{\partial}}\tilde f(z)\vert ={\mathcal  O}(\vert {\im}z\vert^\infty)$ uniformly on $\C$. 
Then, we have (see, e.g., \cite{DiSj1, Ma2}),
\be
\label{calfonc}
f(A)=\frac1{\pi}\int_\C{\overline{\partial}}\tilde f(z)
(A-z)^{-1} d\re z\; d\im z.
\ee
Now, by Proposition \ref{param}, we see that, for $z\in  \C\backslash\R$, the operator $(A-z)^{-1}$ is a ${\mathcal  U}$-twisted $h$-admissible operator.  Moreover,  by standard rules on the operations ${\rm ad}_{{\i}}$, if $A$ and $B$ are two bounded operators, then, for any $N\geq 1$ and any ${\i}_1,\cdots,{\i}_N\in  C_b^\infty (\R^n)$, one has,
$$
{\rm ad}_{{\i}_1}\circ\cdots\circ{\rm ad}_{{\i}_N}(AB)=\sum_{\genfrac{}{}{0pt}{}{I\cup J =\{1,\dots,N\}}{I\cap J=\emptyset}}\left( \prod_{i\in I}{\rm ad}_{{\i}_i}\right)(A)\left( \prod_{j\in J}{\rm ad}_{{\i}_j}\right)(B).
$$
In particular, replacing $A$ and $B$ by $A-z$ and $(A-z)^{-1}$ respectively, one obtains,
\begin{eqnarray*}
&&{\rm ad}_{{\i}_1}\circ\cdots\circ{\rm ad}_{{\i}_N}((A -z)^{-1})\\
&& =-(A-z)^{-1}\sum_{\genfrac{}{}{0pt}{}{I\cup J =\{1,\dots,N\}}{I\cap J=\emptyset,\, I\not=\emptyset}}  \left( \prod_{i\in I}{\rm ad}_{{\i}_i}\right)(A-z)\left( \prod_{j\in J}{\rm ad}_{{\i}_j}\right)((A -z)^{-1}),
\end{eqnarray*}
and thus, an easy induction gives,
$$
{\rm ad}_{{\i}_1}\circ\cdots\circ{\rm ad}_{{\i}_N}((A -z)^{-1})={\mathcal  O}(h^N\vert\im z\vert^{-(N+1)}),
$$
uniformly with respect to $h$ and $z$. Therefore, it is immediate from (\ref{calfonc}) that $f(A)$ verifies the first condition in Proposition \ref{eqdef}.
\vskip 0.3cm
Moreover,  setting $(a_j)_{0\leq j\leq r}:=\sigma (A)$, for $\i_j\in C_d^\infty (\omega_j)$,  we denote
by $B_j(z)$ the $h$-admissible operator with symbol $ (\re a_j-\overline z)(\re a_j-z)\varphi_j^4 + \psi_j^2$, where $\varphi_j$ and $\psi_j$ are as at the end of the proof of Proposition \ref{param}. Then, using that $a_j=\re a_j + {\mathcal  O}(h)$, we see that
$$
B_j(z) = A_j(z)^*A_j(z)+ \psi_j^2 +hB'_j(z),
$$
with $A_j (z)= U_j\varphi_j (A-z)U_j^{-1}\varphi_j$, and $B'_j(z)$ is a uniformly bounded $h$-admissible operator. As a consequence, if $v\in L^2(\R^n; {\mathcal H})$, and for $\im z\not=0$, 
a computation similar to that of (\ref{pseudinv}) shows that,
\begin{equation}
\label{bprimej}
B_j(z) U_j\i_j u_j(z)= C_j(z)v + hC'_j(z)U_j\varphi_ju_j(z)+{\mathcal O}(h^\infty )\Vert v\Vert,
\end{equation}
where $u_j(z):= (A-z)^{-1}U_j^{-1}\i_j v$, and $C_j(z), C'_j(z)$  are uniformly bounded $h$-admissible operators. Then, denoting by $\tilde B_j(z)$ the $h$-admissible operator with symbol $ [(\re a_j-\overline z)(\re a_j-z)\varphi_j^4 + \psi_j^2]^{-1}$, the standard pseudodifferential calculus with operator-valued symbols shows that,
$$
 \Vert \tilde B_j(z)\Vert ={\mathcal O}(|\im z|^{-N_0})
$$
for some $N_0\geq 1$, and,
$$
\tilde B_j(z) B_j (z)=1+hR_j(z),
$$
where $R_j(z)$ is a $h$-admissible operator with symbol $r_j(z)$ verifying $\partial_{x,\xi}^\alpha r_j(z)={\mathcal  O}(|\im z|^{-N_{\alpha,j}})$, for all $\alpha\in\Z_+^{2n}$, and for some $N_{\alpha,j}\geq 1$. Thus, applying $B'_j(z)$ to (\ref{bprimej}), we obtain,
$$
U_j\i_j u_j(z)= C^{(1)}_j(z)v + hC^{(2)}_j(z)U_j\varphi_ju_j(z)+{\mathcal O}(h^\infty |\im z|^{-N_1})\Vert v\Vert,
$$
where $C^{(1)}_j(z), C^{(2)}_j(z)$  are  two $h$-admissible operators, uniformly bounded by some negative power of $ |\im z|$, and $N_1$ is some positive number. Again, iterating the procedure as in  the proof of Proposition \ref{param}, one can deduce that $f(A)$ also verifies the second condition in Proposition \ref{eqdef}, and therefore is a ${\mathcal U}$-twisted $h$-admissible operator.
\vskip 0.3cm
Finally, a computation similar to that of (\ref{invapprox}) shows that,
$$
(A-z)^{-1} =  \sum_{j=0}^r U_j^{-1}\varphi_j^3 {\rm Op}_h(\varphi_j (\re a_j -z)^{-1})U_j{\i}_j + hR
$$
where $\varphi_j$ and $\i_j$  are as in (\ref{invapprox}), and $R$ verifies,
$$
U_j\tilde\i_jRU_j^{-1}\tilde\i_j = {\rm Op}_h(\sum_{k=0}^Nh^kr_{k,j}(z)) +{\mathcal  O}(h^N|\im z|^{-N_1(N)}),
$$
for any $\tilde\i_j\in C_d^\infty (\Omega_j)$ such that $\tilde\i_j\varphi_j =\tilde\i_j\i_j =\tilde\i_j$, any $N\geq 1$, and for some $N_1(N)\geq 1$ and $r_{k,j}(z)\in C^\infty (T^*\Omega_j)$, $\partial^\alpha r_{k,j}(z)={\mathcal  O}(|\im z|^{-N_{\alpha,k,j}})$ uniformly. Then, one easily concludes that the symbol $b_j$ of $U_j\tilde\i_jf(A)U_j\tilde\i_j$ verifies,
$$
b_j= \tilde\i_j f(\re a_j )\tilde\i_j+{\mathcal  O}(h),
$$
and since the previous construction can be made for $\tilde\i_j\in C_d^\infty (\Omega_j)$ arbitrary, the result on the symbol of $f(A)$ follows. \hfill$\bullet$
\vskip 0.3cm
In order to complete the theory of bounded ${\mathcal  U}$-twisted $h$-admissible operators, it remains to generalize the notion of quantization. To this purpose, we first observe that, if $a =(a_j)_{j=0,\dots ,r}\in {\bf S}(\Omega ;{\mathcal  L}(H))$, then, the two operators $\varphi_j{\rm Op}_h(a_j)\varphi_j$ and $U_j^{-1}\varphi_j{\rm Op}_h(a_j)U_j\varphi_j$ are well defined  for any $\varphi_j  \in  C_d^\infty (\Omega_j)$. Moreover, if $a=\sigma (A)$ 
is the symbol of a ${\mathcal  U}$-twisted
$h$-admissible operator $A$, then, by construction, it necessarily verifies the following condition of compatibility: 
\be
\label{compatib}
U_j^{-1}\varphi {\rm Op}_h(a_j)U_j \varphi = U_k^{-1}\varphi {\rm Op}_h(a_k)U_k\varphi ,
\ee
for any $\varphi \in C_d^\infty (\Omega_j)\cap C_d^\infty (\Omega_k)$.
Then, we have,
\begin{theorem}[Quantization]\sl 
Let $a =(a_j)_{j=0,\dots ,r}\in {\bf S}(\Omega ;{\mathcal  L}(H))$ satisfying to the compatibility condition (\ref{compatib}). Then, there exists a ${\mathcal  U}$-twisted $h$-admissible operator $A$, unique up to ${\mathcal  O}(h^\infty)$, such that $a=\sigma (A)$. Moreover, $A$ is given by the formula,
\be
\label{quantif}
A=\sum_{j=0}^rU_j^{-1}\chi_j{\rm Op}_h(a_j)U_j\varphi_j,
\ee
where $\i_j, \varphi_j \in C_d^\infty (\Omega_j)$($j=0, \dots ,r$) is any family of functions such that $\sum_{j=0}^r\i_j =1$ and ${\rm dist} \left( \supp (\varphi_{j}-1), \supp {\i}_j\right) >0$.
\end{theorem}
{\em Proof -- } The unicity up to ${\mathcal  O}(h^\infty)$ is a direct consequence of the formulas (\ref{reconst1})-(\ref{reconst2}), where $A$ is expressed in terms of $U_j\varphi_j AU_j^{-1}\varphi_j$ and is  clearly ${\mathcal  O}(h^\infty)$ if these operators have identically vanishing symbols. For the existence, we define $A$ as in (\ref{quantif}) and we observe that, thanks to (\ref{compatib}), for any $k\in\{0,\dots, r\}$ and $\psi_k\in C_d^\infty (\Omega_k)$, one has,
\begin{eqnarray*}
U_k\psi_kAU_k^{-1}\psi_k &=&\sum_{j=0}^r\chi_j\psi_k{\rm Op}_h(a_k)\varphi_j\psi_k=\sum_{j=0}^r\chi_j\psi_k{\rm Op}_h(a_k)\psi_k +{\mathcal  O}(h^\infty)\\
&=&  \psi_k{\rm Op}_h(a_k)\psi_k+{\mathcal  O}(h^\infty).
\end{eqnarray*}
Thus, $A$ admits $(a_k)_{k=0,\dots ,r}$ as  its symbol, and the result follows.\hfill$\bullet$
\vskip 0.3cm
To end this chapter, let us go back to our operator $\tilde P$ defined at the end of Chapter \ref{modop}. We have,
\begin{proposition}\sl  
\label{putruc}
Assume (H1)-(H4). Then, the operator $\tilde P$ defined in (\ref{ptilde}) is such that 
$\tilde P (\boldsymbol{\omega} + Q_0)^{-1}$ is a ${\mathcal  U}$-twisted $h$-admissible operator on  $L^{2}({\R^n};{\mathcal  H})$,  where  ${\mathcal  U}=(U_j,\Omega_j)_{j=0,1,\cdots,r}$ is the regular covering defined in Chapter \ref{sect2}. Moreover, its  symbol $\tilde p=(\tilde p_j)_{j=0,1,\cdots,r}$ verifies,
$$
\tilde p_j(x,\xi ) = (\omega (x,\xi ) + \tilde Q_j(x) + \zeta (x)W(x))(\omega (x,\xi ) + Q_{0,j}(x))^{-1} +hb_j,
$$
where $(\tilde Q_j(x))_{j=0,1,\cdots,r}$ (resp. $(Q_{0,j}(x))_{j=0,1,\cdots,r}$) is the  symbol of $\tilde Q(x)$ (resp. $Q_0(x)$), and $(b_j)_{j=0,\dots ,r}\in {\bf S}(\Omega ;{\mathcal  L}(H))$.
\end{proposition}
{\em Proof -- } We must verify the two conditions of Proposition \ref{eqdef}. We have,
\begin{eqnarray*}
&&{\rm ad}_{{\i}}(\tilde P(\boldsymbol{\omega} + Q_0)^{-1})\\
&&\hskip 0.5cm={\rm ad}_{\i} (\tilde P)(\boldsymbol{\omega} + Q_0)^{-1} +\tilde P{\rm ad}_{{\i}}((\boldsymbol{\omega} + Q_0)^{-1})\\
&&\hskip 0.5cm= {\rm ad}_{\i} (\boldsymbol{\omega})(\boldsymbol{\omega} + Q_0)^{-1} -\tilde P(\boldsymbol{\omega} + Q_0)^{-1}{\rm ad}_{{\i}}(\boldsymbol{\omega})(\boldsymbol{\omega} + Q_0)^{-1}\\
&&\hskip 0.5cm= {\mathcal  O}(h),
\end{eqnarray*}
and an easy iteration shows that the first condition of Proposition \ref{eqdef} is satisfied. Moreover, if ${\i}_j ,\tilde{\i}_j\in  C_b^\infty (\R^n)$ are supported in $\Omega_j$ ($j=1\cdots,r$) and verify $\supp {\i}_j\cap \supp (1-\tilde{\i}_j) =\emptyset$, and if we set $P_j:= U_j{\i}_j\tilde PU_j^{-1}\tilde{\i}_j$, we have,
\begin{eqnarray*}
&&U_j{\i}_j\tilde P(\boldsymbol{\omega} + Q_0)^{-1}U_j^{-1}{\i}_j\\
&&\hskip 1cm =
U_j{\i}_j\tilde P\tilde{\i}_j^2(\boldsymbol{\omega} + Q_0)^{-1}U_j^{-1}{\i}_j+ U_j{\i}_j\boldsymbol{\omega}(1-\tilde{\i}_j^2)(\boldsymbol{\omega} + Q_0)^{-1}U_j^{-1}{\i}_j\\
&&\hskip 1cm =P_jU_j\tilde{\i}_j(\boldsymbol{\omega} + Q_0)^{-1}U_j^{-1}{\i}_j +{\mathcal  O}(h^\infty),
\end{eqnarray*}
and a slight generalization of the last argument in the proof of Proposition \ref{param} (this time with $B_j=U_j\varphi_j(\boldsymbol{\omega} +Q_0)U_j^{-1}\varphi_j +\psi_j(\boldsymbol{\omega} +Q_0)\psi_j$), shows that $P_jU_j\tilde{\i}_j(\boldsymbol{\omega} + Q_0)^{-1}U_j^{-1}{\i}_j$ is a bounded $h$-admissible operator on $L^2(\R^n;{\mathcal  H})$. Therefore, the second condition of Proposition \ref{eqdef} is satisfied, too, and the result follows.\hfill$\bullet$
\begin{corollary}\sl 
\label{pqutruc}
The two operators $(\tilde P+i)^{-1}$ and $(\boldsymbol{\omega} +Q_0)^{-1}$ are ${\mathcal  U}$-twisted $h$-admissible operators on  $L^{2}({\R^n};{\mathcal  H})$.
\end{corollary}
{\em Proof -- } First observe that the previous proof is still valid if $\tilde P$ is changed into $\tilde P +1$. This proves that $(\boldsymbol{\omega} +Q_0)^{-1} = (\tilde P +1)(\boldsymbol{\omega} +Q_0)^{-1} - \tilde P(\boldsymbol{\omega} +Q_0)^{-1}$ is a ${\mathcal  U}$-twisted $h$-admissible operator. Moreover, since $(\tilde P +i)(\boldsymbol{\omega} +Q_0)^{-1}$ is elliptic, by Proposition \ref{param} its inverse $(\boldsymbol{\omega} +Q_0)(\tilde P +i)^{-1}$ is a ${\mathcal  U}$-twisted $h$-admissible operator, too. Therefore, so is $(\tilde P +i)^{-1} = (\boldsymbol{\omega} +Q_0)^{-1}\left[ (\boldsymbol{\omega} +Q_0)(\tilde P +i)^{-1}\right]$.\hfill$\bullet$
\hskip 0.3cm
\begin{proposition}\sl 
\label{fdeP}
For any $f\in C_0^\infty (\R)$, the operator $f(\tilde P)$ is a ${\mathcal  U}$-twisted $h$-admissible.
\end{proposition}
{\em Proof -- }
By Proposition \ref{putruc} and Corollary \ref{pqutruc}, we  see that the operator $(\tilde P -z)(\boldsymbol{\omega} +Q_0)^{-1}$ is a ${\mathcal  U}$-twisted $h$-admissible operator, and it is elliptic 
for $z\in  \C\backslash\R$. Therefore, by Proposition \ref{param}, its inverse $(\boldsymbol{\omega} +Q_0)(\tilde P -z)^{-1}$ is a ${\mathcal  U}$-twisted $h$-admissible operator, too.  Moreover, for any $N\geq 1$ and any ${\i}_1,\cdots,{\i}_N\in  C_b^\infty (\R^n)$, one has,
$$
{\rm ad}_{{\i}_1}\circ\cdots\circ{\rm ad}_{{\i}_N}((\boldsymbol{\omega} +Q_0)(\tilde P -z)^{-1})={\mathcal  O}(h^N\vert\im z\vert^{-(N+1)})
$$
uniformly with respect to $h$ and $z$. Therefore, we deduce again from (\ref{calfonc}) that $(\boldsymbol{\omega} +Q_0)f(\tilde P)$, too, is a ${\mathcal  U}$-twisted $h$-admissible operator. As a consequence, so is $f(\tilde P)$. \hfill$\bullet$
\chapter{Twisted Partial Differential Operators}
\label{twpdo}
\setcounter{equation}{0}%
\setcounter{theorem}{0}%
For $\mu \geq 0$, we set,
$$
H^\mu_d(\Omega_j):=\{ u\in L^2(\Omega_j;{\mathcal  H})\, ;\, \forall \i_j\in C_d^\infty (\Omega_j), \,\i_j u\in H^\mu (\R^n;{\mathcal  H}) \},
$$
where $H^\mu (\R^n;{\mathcal  H}) $ stands for the usual Sobolev space of order $\mu$ on $\R^n$ with values in $\mathcal  H$. Moreover, if ${\mathcal  U}:=\left( U_j,\Omega_j\right)_{j=0,\cdots,r}$ is a regular unitary covering (in the previous sense) of $L^2(\R^n;{\mathcal  H})$, we introduce the vector-space,
$$
{\mathcal  H}^\mu_d({\mathcal  U}):= \{ u\in L^2(\R^n;{\mathcal  H})\, ;\, \forall\, j=0,\dots,r,\, U_ju\left\vert_{\Omega_j}\right. \in H^\mu_d(\Omega_j)\},
$$
endowed with the family of semi-norms,
$$
\Vert u\Vert_{\mu ,\chi}:=\Vert u\Vert_{L^2}+\sum_{j=0}^r\Vert U_j\i_j u\Vert_{H^\mu},
$$
where $\i :=(\i_j)_{j=0,\dots,r}$ is such that $\i_j\in C_d^\infty (\Omega_j)$ for all $j$. In particular, we have a notion of continuity for operators $A:{\mathcal  H}^\mu_d({\mathcal  U})\to {\mathcal  H}^\nu_d({\mathcal  U})$.
\vskip 0.3cm
Let us also remark that, for $\mu =0$, we recover ${\mathcal  H}^0_d({\mathcal  U})=L^2 (\R^n;{\mathcal  H})$, and, if $\mu \geq\nu$, then ${\mathcal  H}^\mu_d({\mathcal  U} )\subset {\mathcal  H}^\nu_d({\mathcal  U})$ with a continuous injection.
\begin{definition}\sl
\label{pdoutruc}
Let  ${\mathcal  U}:=\left( U_j,\Omega_j\right)_{j=0,\cdots,r}$ be a regular unitary covering (in the previous sense) of $L^2(\R^n;{\mathcal  H})$, and let $\mu \in\Z_+$. 
We say that an  operator $A:{\mathcal  H}^\mu_d({\mathcal  U})\to L^{2}({\R^n};{\mathcal  H})$ is a  (semiclassical) ${\mathcal  U}$-twisted partial differential operator  up to regularizing unitary conjugation (in short: ${\mathcal  U}$-twisted PDO) of degree $\mu$, if $A$ is local with respect to the variable $x$ (that is, $ \supp (Au)\subset \supp u$ for all $u$, where $\supp$ stands for the support with respect to $x$), and, for all $j=0,\dots ,r$, the operator  $U_jAU_j^{-1}$ (well defined on $H_d^\mu (\Omega_j)$) is of the form,
$$
U_j AU_j^{-1} = \sum_{|\alpha|\leq \mu} a_{\alpha ,j} (x;h)(hD_x)^\alpha
$$
with $a_{\alpha ,j} \in S(\Omega_j ; {\mathcal  L}({\mathcal  H}))$.
\end{definition}
In particular, for any partition of unity $({\i}_{j})_{j=0,\dots ,r}$ on $\R^n$ with ${\i}_{j}\in  C_d^{\infty}(\Omega_j)$, $A$ can be
written as,
\be
\label{opd}
A=\sum_{j=0}^rU_j^{-1} A_jU_j{\i}_j ,
\ee
with $A_j:=U_j AU_j^{-1}$. As a consequence, one also has ${\rm ad}_{{\i}_1}\circ\cdots\circ {\rm ad}_{{\i}_{\mu +1}}(A) = 0$ for any  functions ${\i}_1,\cdots,{\i}_{\mu +1}\in  C_b^\infty (\R^n)$.
\vskip 0.3cm
Of course, we also  have an obvious notion of (full) symbol for such operators, namely, the family,
$$
\sigma (A):= (a_j)_{0\leq j\leq r},\quad a_j(x,\xi ;h):=\sum_{|\alpha|\leq \mu} a_{\alpha ,j} (x;h)\xi^\alpha.
$$
Moreover, if $A$ and $B$ are two ${\mathcal  U}$-twisted PDO's on  $L^{2}({\R^n};{\mathcal  H})$, of respective degrees $\mu$ and $\mu'$, by writing $U_jABU_j^{-1}= (U_jAU_j^{-1})(U_jBU_j^{-1})$ and by using a partition of unity as before, we immediately see that $AB$ is well defined on ${\mathcal  H}^{\mu+\mu'}_d({\mathcal  U})$, and is a ${\mathcal  U}$-twisted PDO, too, with symbol,
$$
\sigma(AB) =\sigma(A)\sharp \sigma(B).
$$
Now, we turn back again to the operator $\tilde P$ defined at the end of Chapter \ref{modop}, and the regular covering defined in Chapter \ref{sect2}.
\begin{proposition}\sl 
\label{regul}
 Let $A$  be a ${\mathcal  U}$-twisted PDO on  $L^{2}({\R^n};{\mathcal  H})$ of degree $\mu$,  where ${\mathcal  U}$ is the regular covering defined in Chapter \ref{sect2}. Then, for any integers $k,\ell$ such that $k+\ell \geq\mu /m$, the operator $(\tilde P +i)^{-k}A(\tilde P +i)^{-\ell}$ is a ${\mathcal  U}$-twisted $h$-admissible operator.
\end{proposition}
{\em Proof -- } We first consider the case $k =0$. For $\varphi_j,\psi_j\in C_d ^\infty (\Omega_j)$, such that $\dist (\supp (\psi_j-1), \supp\varphi_j)>0$, we have,
\be
\label{apml}
U_j\varphi_j A(\tilde P +i)^{-\ell}U_j^{-1}\varphi_j =U_j\varphi_j AU_j^{-1}\psi_jU_j\psi_j(\tilde P +i)^{-\ell}U_j^{-1}\varphi_j ,
\ee
and, as in the proof of Proposition \ref{param}, we see that the inverse of $(\tilde P +i)^\ell$ can be written as,
\be
\label{pml1}
(\tilde P +i)^{-\ell} = B (1+hR)
\ee
where $R$ is uniformly bounded, and  $B$ is of the form,
\be
\label{pml2}
B=\sum_{\nu =0}^r U_\nu^{-1}\tilde {\i}_\nu{\rm Op}_h((p_\nu+i)^{-\ell})U_\nu \i_\nu,
\ee
where $(\i_\nu)_{\nu=0,\dots ,r}$ is an arbitrary partition of unity with $\i_\nu\in C_d^\infty (\Omega_\nu)$, $\tilde{\i}_\nu\in C_d^\infty (\Omega_\nu)$ is such that $\tilde {\i}_\nu \i_\nu =\i_\nu$, and $p_\nu(x,\xi ;h) = \omega (x,\xi ;h) + \tilde Q_\nu(x) + \zeta (x) W(x)$.
\begin{lemma}\sl  Let $j\in \{ 0,\dots ,r\}$ and $\psi_j\in C_d^\infty (\Omega_j)$ be fixed.
Then, there exists a partition of unity $(\i_\nu)_{\nu=0,\dots ,r}$ of $\R^n$  with $\i_\nu\in C_d^\infty (\Omega_\nu)$, and there exists $\tilde{\i}_\nu\in C_d^\infty (\Omega_\nu)$ with $\tilde {\i}_\nu \i_\nu =\i_\nu$ ($\nu=0,\dots ,r$), such that $\i_j \psi_j= \psi_j$ and $\tilde{\i}_\nu \psi_j=0$ if $\nu\not= j$.
\end{lemma}
{\em Proof -- } \hskip 0.3cm It is enough to construct a partition of unity in such a way that  $\dist (\supp\psi_j , \supp (\i_j -1))>0$ (and thus, automatically, one will also have $\dist (\supp\psi_j , \supp \i_\nu )>0$ for $\nu\not= j$). Let $(\i'_\nu)_{\nu=0,\dots,r}$ be a partition of unity as in Definition \ref{regcov}, and let $\i''_j\in C_d^\infty (\Omega_j ; [0,1])$ such that $\i''_j =1$ in a neighborhood of $\supp\psi_j \cup \supp \i_j$. Then,  the result is obtained by taking $\i_\nu := (1-\i''_j)\i'_\nu$ if $\nu\not= j$, and $\i_j:= \i''_j$.\hfill$\bullet$
\vskip 0.3cm
Taking the $\i_\nu$'s and $\tilde{\i}_\nu$'s as in the previous lemma, we obtain from (\ref{pml1})-(\ref{pml2}),
$$
U_j\psi_j(\tilde P +i)^{-\ell} = \psi_j{\rm Op}_h((p_j+i)^{-\ell})U_j\i_j(1+hR),
$$
and thus, since $U_j\varphi_j AU_j^{-1}\psi_j$ is a differential operator of degree $\mu$ with operator-valued symbol, we easily deduce from (\ref{apml}) that if $m\ell \geq \mu$, then $A(\tilde P +i)^{-\ell}$ is bounded on $L^2(\R^n;{\mathcal  H})$, uniformly with respect to $h>0$. Moreover, writing,
$$
U_j\varphi_j A(\tilde P +i)^{-\ell}U_j^{-1}\varphi_j = [ U_j\varphi_j AU_j^{-1}\psi_j\la hD_x\ra^{-m\ell}] [ \la hD_x\ra^{m\ell}U_j\psi_j(\tilde P +i)^{-\ell}U_j^{-1}\varphi_j],
$$
and using the standard pseudodifferential calculus with operator-valued symbol for the first factor, and a slight refinement of (\ref{pseudinv}) for the second one, we see that $U_j\varphi_j A(\tilde P +i)^{-\ell}U_j^{-1}\varphi_j$ is an $h$-admissible operator on $L^2(\R^n);{\mathcal  H})$. Then, it only remains to verify the first property of Proposition \ref{eqdef}. We first prove,
\begin{lemma}\sl  For any $\alpha_1,\dots,\alpha_N\in C_b^\infty (\R^n)$, one has,
\be
\label{compchi}
{\rm ad}_{\alpha_1}\circ \dots \circ {\rm ad}_{\alpha_N}((\tilde P +i)^{-\ell}) = h^N(\tilde P +i)^{-\ell}R_N,
\ee
with $R_N={\mathcal  O}(1)$ on $L^2(\R^n;{\mathcal  H})$. 
\end{lemma}
{\em Proof -- } Since ${\rm ad}_{\alpha_N}((\tilde P +i)^{-\ell})=-(\tilde P +i)^{-\ell}{\rm ad}_{\alpha_N}((\tilde P +i)^{\ell})(\tilde P +i)^{-\ell}$, by an easy iteration we see that it is enough to prove that $h^{-N}{\rm ad}_{\alpha_1}\circ \dots \circ {\rm ad}_{\alpha_N}((\tilde P +i)^{\ell})(\tilde P +i)^{-\ell}$ is uniformly bounded. Moreover, since $ {\rm ad}_{\alpha_N}((\tilde P +i)^{\ell})(\tilde P +i)^{-\ell}$ is a sum of terms of the type $(\tilde P +i)^k{\rm ad}_{\alpha_N}(\boldsymbol{\omega})(\tilde P +i)^{-k-1}$ ($0\leq k\leq \ell -1$), another easy iteration shows that it is enough to prove that $h^{-N}(\tilde P +i)^\ell {\rm ad}_{\alpha_1}\circ \dots \circ {\rm ad}_{\alpha_N}(\boldsymbol{\omega})(\tilde P +i)^{-\ell-1}$ is uniformly bounded. Now, by (H4), we see that, for any partition of unity $(\i_j)$ as before, $(\tilde P +i)^\ell $ can be written as,
$$
(\tilde P +i)^\ell = \sum_{j=0}^r  U_j^{-1}P_{j,\ell}U_j \i_j,
$$
where  $P_{j,\ell}$ is of the form,
$$
P_{j,\ell} =\sum_{|\alpha|\leq m\ell}\rho_{j,\ell,\alpha}(x;h)(hD_x)^\alpha,
$$
with $\rho_{j,\ell,\alpha}Q_0^{\frac{|\alpha|}{m}-\ell}\in C^\infty (\Omega_j; {\mathcal  H})$. Moreover, by (\ref{conjomega}), the operator $U_j{\rm ad}_{\alpha_1}\circ \dots \circ {\rm ad}_{\alpha_N}(\boldsymbol{\omega})U_j^{-1}={\rm ad}_{\alpha_1}\circ \dots \circ {\rm ad}_{\alpha_N}(U_j\boldsymbol{\omega}U_j^{-1})$ is of the form,
$$
U_j{\rm ad}_{\alpha_1}\circ \dots \circ {\rm ad}_{\alpha_N}(\boldsymbol{\omega})U_j^{-1}=h^N\sum_{|\alpha|\leq (m-N)_+}\tau_{j,\alpha}(x;h)(hD_x)^\alpha,
$$
with $\tau_{j,\alpha}Q_0^{\frac{|\alpha|}{m}-1}\in C^\infty (\Omega_j; {\mathcal  H})$. In particular, we obtain,
$$
(\tilde P +i)^\ell {\rm ad}_{\alpha_1}\circ \dots \circ {\rm ad}_{\alpha_N}(\boldsymbol{\omega})=h^N\sum_{j=0}^r  \sum_{|\alpha|\leq m(\ell +1)}U_j^{-1}\lambda_{j,\ell,\alpha}(x;h)(hD_x)^\alpha U_j \varphi_j,
$$
with $\varphi_j\in C_d^\infty (\Omega_j)$ and $\lambda_{j,\ell,\alpha}Q_0^{\frac{|\alpha|}{m}-\ell -1}\in C^\infty (\Omega_j; {\mathcal  H})$, and the result follows as before by using (\ref{pml1})-(\ref{pml2}), and by observing that, for $|\alpha|\leq m(\ell + 1)$, the operator  $Q_0^{1+\ell -\frac{|\alpha|}{m}}(hD_x)^{\alpha}
(\la hD_x\ra^m + Q_0)^{-\ell -1}$ is uniformly bounded, and thus so is the operator
$Q_0^{1+\ell -\frac{|\alpha|}{m}}(hD_x)^{\alpha}\varphi_j {\rm Op}_h((p_j+i)^{-\ell -1})U_j\i_j$. \hfill$\bullet$
\vskip 0.3cm
On the other hand, we see on (\ref{opd}) that ${\rm ad}_{{\i}_1}\circ \dots \circ {\rm ad}_{{\i}_N}(A)$ is a ${\mathcal  U}$-twisted PDO of degree $(\mu - N)_+$, and the  first property of Proposition \ref{eqdef} for  $A(\tilde P +i)^{-\ell}$ follows easily. 
\vskip 0.3cm
For the case $k>0$, by taking a partition of unity, we first observe that,
$$
 (\tilde P +i)^{-k}A(\tilde P +i)^{-\ell} = \sum_{j=0}^r (\tilde P +i)^{-k}U_j^{-1}A_jU_j\i_j(\tilde P +i)^{-\ell}
 $$
 where $A_j = U_jAU_j^{-1}$ can be written as,
 $$
 A_j=\sum_{\genfrac{}{}{0pt}{}{|\alpha| \leq mk}{|\beta|\leq m\ell}} (hD_x)^\alpha a_{\alpha ,\beta ,j}(x;h)(hD_x)^\beta.
 $$
 Then, by using (in addition to (\ref{pml1})-(\ref{pml2})) that,
$$
(\tilde P +i)^{-k} =  (1+hR')B'
$$
where $R'$ is uniformly bounded, and  $B'$ is of the form,
$$
B'=\sum_{\nu =0}^r U_\nu^{-1} {\i}_\nu{\rm Op}_h((p_\nu+i)^{-\ell})U_\nu \tilde\i_\nu,
$$
the same previous arguments show that the operator $(\tilde P +i)^{-k}A(\tilde P +i)^{-\ell}$ is bounded on $L^2(\R^n;{\mathcal  H})$, uniformly with respect to $h>0$.
\vskip 0.3cm
Then, let $N\geq 1$ and $\alpha_1\dots,\alpha_N\in C_d^\infty (\Omega_j)$, such that $\alpha_1\varphi_j =\varphi_j$, $\alpha_2\alpha_1 =\alpha_1$, ... , $\alpha_N\alpha_{N-1} =\alpha_{N-1}$, and $\alpha_N(\psi_j -1)=0$.
We have,
\begin{eqnarray*}
&&U_j\varphi_j (\tilde P +i)^{-k}A(\tilde P +i)^{-\ell}U_j^{-1}\varphi_j\\
&&\hskip 1cm =U_j\varphi_j (\tilde P +i)^{-k}A\psi_j(\tilde P +i)^{-\ell}U_j^{-1}\varphi_j\\
&&\hskip 2cm +U_j\varphi_j (\tilde P +i)^{-k}A(\psi_j -1){\rm ad}_{\alpha_1}\circ \dots \circ {\rm ad}_{\alpha_N}((\tilde P +i)^{-\ell})U_j^{-1}\varphi_j
\end{eqnarray*}
and thus, by (\ref{compchi}),
\begin{eqnarray*}
&&U_j\varphi_j (\tilde P +i)^{-k}A(\tilde P +i)^{-\ell}U_j^{-1}\varphi_j\\
&&\hskip 3cm =U_j\varphi_j (\tilde P +i)^{-k}A\psi_j(\tilde P +i)^{-\ell}U_j^{-1}\varphi_j+{\mathcal  O}(h^N).
\end{eqnarray*}
Then, writing $A\psi_j = U_j^{-1}\tilde\psi_j A_j U_j\psi_j$, with $A_j = U_j AU_j^{-1}$ and $\tilde\psi_j\in C_d^\infty (\Omega_j)$ such that $\tilde\psi_j \psi_j =\psi_j$, the result is obtained along the same lines as before.\hfill$\bullet$
\begin{proposition}\sl 
\label{omegaPDO}
 The two operators $\boldsymbol{\omega}Q_0^{-1}$ and  $Q_0^{-1}\boldsymbol{\omega}$ are ${\mathcal  U}$-twisted PDO's  of degree m. Moreover, if $A$ is a ${\mathcal  U}$-twisted PDO such that $Q_0A$ and $AQ_0$ are  ${\mathcal  U}$-twisted PDO's, too, of degree $\mu$, then the operator $h^{-1}[\boldsymbol{\omega}, A]$ is a  ${\mathcal  U}$-twisted PDO of degree at most $\mu + m-1$.
\end{proposition}
{\em Proof -- } Thank to (H4), the fact that $\boldsymbol{\omega}Q_0^{-1}$ and  $Q_0^{-1}\boldsymbol{\omega}$ are ${\mathcal  U}$-twisted PDO's   of degree $m$ is obvious. Moreover, the fact that  $Q_0A$ and $AQ_0$ are both ${\mathcal  U}$-twisted PDO's  implies that $U_jAU_j^{-1}$ can be written as,
$$
U_j AU_j^{-1} = \sum_{|\alpha|\leq \mu} a_{\alpha ,j} (x;h)(hD_x)^\alpha
$$ 
with $Q_0a_{\alpha ,j}$ and  $a_{\alpha ,j}Q_0$ in $S(\Omega_j ; {\mathcal  L}({\mathcal  H}))$. Then, using (H4), we have,
\begin{eqnarray*}
U_j\boldsymbol{\omega}AU_j^{-1} &=& \sum_{\genfrac{}{}{0pt}{}{|\alpha| \leq m}{|\beta|\leq \mu}} c_\alpha (x;h)(hD_x)^\alpha a_{\beta ,j} (x;h)(hD_x)^\beta \\
&&\hskip 1cm +h\sum_{\genfrac{}{}{0pt}{}{|\alpha| \leq m-1}{|\beta|\leq \mu}} \omega_{\alpha ,j}(x;h)(hD_x)^\alpha a_{\beta ,j} (x;h)(hD_x)^\beta
\end{eqnarray*}
and
\begin{eqnarray*}
U_jA\boldsymbol{\omega}U_j^{-1} &=& \sum_{\genfrac{}{}{0pt}{}{|\alpha| \leq m}{|\beta|\leq \mu}} a_{\beta ,j} (x;h)(hD_x)^\beta c_\alpha (x;h)(hD_x)^\alpha \\
&&\hskip 1cm +h\sum_{\genfrac{}{}{0pt}{}{|\alpha| \leq m-1}{|\beta|\leq \mu}} a_{\beta ,j} (x;h)(hD_x)^\beta\omega_{\alpha ,j}(x;h)(hD_x)^\alpha.
\end{eqnarray*}
Moreover, by (H4) (and the fact that  $U_j\boldsymbol{\omega}U_j^{-1}$ is symmetric), we know that $c_\alpha $ is scalar-valued,  and $Q_0^{-1}\omega_{\alpha ,j}$,  $\omega_{\alpha ,j}Q_0^{-1}$ are bounded operators on $ {\mathcal  H}$ (together with all their derivatives). Thus,  it is clear that $h^{-1}U_j[\boldsymbol{\omega},A]U_j^{-1} $ is a PDO of degree $\leq$ $\mu + m -1$, and the result follows.\hfill$\bullet$
\chapter{Construction of a Quasi-Invariant Subspace}
\label{qis}
\setcounter{equation}{0}%
\setcounter{theorem}{0}%
\begin{theorem}\sl 
\label{th:constrPi}
Assume (H1)-(H4), and denote by  ${\mathcal  U}:=(U_j, \Omega_j)_{j=0,\cdots,r}$  the regular unitary covering of $L^2(\R^n ;{\mathcal  H})$ constructed from the operators $U_j$ and the open sets $\Omega_j$ defined in Chapter \ref{sect2}. Then,  for any $g\in  C_0^\infty (\R)$,
there exists a ${\mathcal  U}$-twisted $h$-admissible operator 
$\Pi_g$ on $L^2(\R^n ;{\mathcal  H})$, such that $\Pi_g$ is an orthogonal projection that verifies, 
\begin{equation}
\label{eq:pi-pi0}
\Pi_g =\tilde\Pi_0 +{\mathcal  O}(h)
\end{equation}
and, for any $f\in C_0^\infty (\R)$ with $\supp f \subset \{g=1\}$, and any $\ell\geq 0$,
\begin{equation}
\label{eq:comm}
\tilde P^\ell [f(\tilde P),\Pi_g]  ={\mathcal  O}(h^\infty).
\end{equation}
Moreover, $\Pi_g$ is uniformly bounded as an operator : $L^2(\R^n ;{\mathcal  H})\to L^2(\R^n ;{\mathcal  D}_Q)$ and, for any $\ell\geq 0$, any $N\geq 1$, and any  functions ${\i}_1,\cdots,{\i}_N\in  C_b^\infty (\R^n)$, one has,
\be
\label{adpi}
\tilde P^\ell{\rm ad}_{{\i}_1}\circ\cdots\circ {\rm ad}_{{\i}_N}(\Pi_g)  ={\mathcal  O}(h^N).
\ee
\end{theorem}
{\em Proof -- }: 
We first perform a formal construction, by essentially following a procedure taken from  \cite{Ne1} (see also  \cite{BrNo} in the case $L=1$).  In the sequel, all the twisted PDO's  that are involved are  associated with the regular covering ${\mathcal  U}$ constructed in Chapter \ref{sect2}, and we will omit to specify it all the time. We say that a twisted PDO is symmetric when it is formally selfadjoint with respect to the scalar product in $L^2(\R^n; {\mathcal  H})$.
\vskip 0.3cm
Since $\boldsymbol{Q} =\tilde Q(x) + \zeta (x)W(x)$ commutes with $\tilde\Pi_0$, we have,
\be
\label{comppi0}
[\tilde P ,\tilde\Pi_0] = [\boldsymbol{\omega} ,\tilde\Pi_0].
\ee
Moreover, denoting by $\gamma (x)$ a complex oriented single loop surrounding the set $\{\tilde\lambda_{L'+1}(x), \dots , \tilde\lambda_{L'+L}(x)\}$ and leaving the rest of the spectrum of  $\tilde Q(x)$ in its exterior, we have,
\be
\label{reppi}
\tilde\Pi_0(x) =\frac1{2i\pi}\int_{\gamma(x)}(z-\tilde Q(x))^{-1}dz,
\ee
and thus, it  results from Proposition \ref{qtilde} and assumption (H4) that  $Q_0\tilde\Pi_0(x)$ is a ${\mathcal  U}$-twisted PDO of degree 0. Therefore, applying Proposition  \ref{omegaPDO}, we immediately obtain, 
\be
\label{step1}
[\tilde P ,\tilde\Pi_0] =-ihS_0,
\ee
where $S_0$  is a symmetric twisted PDO (of degree $m-1$). Moreover, 
setting $\tilde\Pi_0^\perp := 1-\tilde\Pi_0$, we observe that, 
\be
\label{decS}
 S_0=\tilde\Pi_0 S_0\tilde\Pi_0^\perp + \tilde\Pi_0^\perp S_0\tilde\Pi_0.
\ee
Then, we set,
\be
\label{pi1}
\tilde\Pi_1 := -\frac1{2\pi}\oint_{\gamma (x)}(z-\tilde Q(x))^{-1}\left[ \tilde\Pi_0^\perp(x) S_0\tilde\Pi_0(x) - \tilde\Pi_0(x) S_0\tilde\Pi_0^\perp(x)\right](z-\tilde Q(x))^{-1}dz.
\ee
Thus, $\tilde\Pi_1$ is a symmetric ${\mathcal  U}$-twisted PDO  (of  degree  $m-1$), and is such that  $Q_0\tilde\Pi_1$ is a twisted PDO, too.
Therefore, using Proposition \ref{omegaPDO} again, we have,
$$
[\tilde P, \tilde\Pi_1]  = [\boldsymbol{Q} , \tilde\Pi_1] + hB,
$$
where $B$ is a twisted PDO (of degree $2(m-1)$). Then, using that $\tilde Q(x)(z-\tilde Q(x))^{-1} = (z-\tilde Q(x))^{-1}\tilde Q(x) = z(z-\tilde Q(x))^{-1} -1$, one computes,
\begin{eqnarray*}
[\tilde Q(x), \tilde\Pi_1] &=& \frac1{2\pi}\oint_{\gamma (x)}\left[ \tilde\Pi_0^\perp(x) S_0\tilde\Pi_0(x) - \tilde\Pi_0(x) S_0\tilde\Pi_0^\perp(x)\right](z-\tilde Q(x))^{-1}dz\\
&&  - \frac1{2\pi}\oint_{\gamma (x)}(z-\tilde Q(x))^{-1}\left[ \tilde\Pi_0^\perp(x) S_0\tilde\Pi_0(x) - \tilde\Pi_0(x) S_0\tilde\Pi_0^\perp(x)\right]dz\\
&=& i\left[ \tilde\Pi_0^\perp(x) S_0\tilde\Pi_0(x) - \tilde\Pi_0(x) S_0\tilde\Pi_0^\perp(x)\right]\tilde\Pi_0(x)\\
&&
\hskip 2cm -i\tilde\Pi_0(x)\left[ \tilde\Pi_0^\perp(x) S_0\tilde\Pi_0(x) - \tilde\Pi_0(x) S_0\tilde\Pi_0^\perp(x)\right]\\
&=& i(\tilde\Pi_0^\perp S_0\tilde\Pi_0 + \tilde\Pi_0 S_0\tilde\Pi_0^\perp),
\end{eqnarray*}
that gives,
\be
\label{comqpi1}
[\boldsymbol{Q} , \tilde\Pi_1]=i(\tilde\Pi_0^\perp S_0\tilde\Pi_0 + \tilde\Pi_0 S_0\tilde\Pi_0^\perp)+ [\zeta W, \tilde \Pi_1],
\ee
and thus, using (\ref{decS}), one obtains,
\be
\label{corr1}
[\tilde P,\tilde\Pi_1] = iS_0 -ihS_1,
\ee
where $S_1$  is a symmetric twisted PDO (of degree $2(m-1)$).
Hence, setting,  
$$
\Pi^{(1)}:=\tilde\Pi_0+h\tilde\Pi_1,
$$
we deduce from (\ref{step1}) and (\ref{corr1}),
\be
[\tilde P,\Pi^{(1)}] =- ih^2S_1.
\ee
Moreover, 
$$
(\Pi^{(1)})^2-\Pi^{(1)}=h(\tilde\Pi_0\tilde\Pi_1+\tilde\Pi_1\tilde\Pi_0-\tilde\Pi_1)+h^2\tilde\Pi_1^2
=h^2\tilde\Pi_1^2=:h^2T_1,
$$
where $T_1$ is a symmetric twisted PDO (of degree $2(m-1)$), such that $Q_0T_1$ is a twisted PDO, too. \\
\vskip 0.3cm
Now, by induction on $M$, suppose that we have constructed a symmetric twisted PDO $\Pi^{(M)}$ as, 
$$
\Pi^{(M)}=\sum_{k=0}^M h^k \tilde\Pi_k,
$$
where the $Q_0\tilde\Pi_k$'s are twisted PDO's, such that,
\begin{eqnarray}
\label{projform}
 && (\Pi^{(M)})^2-\Pi^{(M)}=h^{M+1} T_M;\\
 \label{comform}
&&[\tilde P,\Pi^{(M)}]=- ih^{M+1} S_M,
\end{eqnarray}
with
$S_M$ and $Q_0T_M$  twisted PDO's.
\vskip 0.3cm
We set,
$$
\Pi^{(M+1)}=\Pi^{(M)} +h^{M+1}\tilde\Pi_{M+1},
$$
with,
\begin{eqnarray}
\tilde\Pi_{M+1}&:=& -\frac1{2\pi}\oint_{\gamma (x)}(z-\tilde Q(x))^{-1}\left[ \tilde\Pi_0^\perp S_M\tilde\Pi_0 - \tilde\Pi_0 S_M\tilde\Pi_0^\perp\right](z-\tilde Q(x))^{-1}dz\nonumber\\
\label{pim+1}
&&\hskip 4cm + \tilde\Pi_0^\perp T_M\tilde\Pi_0^\perp -\tilde\Pi_0T_M\tilde\Pi_0.
\end{eqnarray}
Then, $\Pi^{(M+1)}$ is again a symmetric twisted PDO, and, using the induction assumption, we immediately see that $\tilde Q(x)\tilde\Pi_{M+1}$ (and thus also $Q_0\tilde\Pi_{M+1}$) is a twisted PDO. Moreover,
since $T_M$ and $\Pi^{(M)}$ commute, we have,
$$
\Pi^{(M)} T_M(1-\Pi^{(M)})=(1-\Pi^{(M)})T_M\Pi^{(M)}   =-h^{M+1}T_M^2,
$$
and thus, since $\Pi^{(M)}=\tilde \Pi_0 + hR_M$ with $Q_0R_M$ twisted PDO, we first obtain,
\be
\label{projT}
\tilde \Pi_0^\perp T_M\tilde \Pi_0 + \tilde \Pi_0T_M\tilde \Pi_0^\perp = hR'_M,
\ee
with $Q_0R'_M$  twisted PDO.
On the other hand, one can check  that,
$$
\tilde\Pi_{M+1}-(\tilde\Pi_0\tilde\Pi_{M+1}+\tilde\Pi_{M+1}\tilde\Pi_0)=\tilde\Pi_0T_{M}\tilde\Pi_0 +\tilde\Pi_0^\perp T_{M}\tilde\Pi_0^\perp,
$$ 
and thus, with (\ref{projT}),
$$
\tilde\Pi_{M+1}-(\tilde\Pi_0\tilde\Pi_{M+1}+\tilde\Pi_{M+1}\tilde\Pi_0)=T_{M} -hR'_M.
$$ 
As a consequence, we obtain,
\be
(\Pi^{(M+1)})^2-\Pi^{(M+1)}=h^{M+2} T_{M+1},
\ee
where $Q_0T_{M+1}$ is a twisted PDO. Applying Proposition \ref{omegaPDO}, we also have,
$$
[\boldsymbol{\omega}, \tilde\Pi_{M+1}] =hR''_M,
$$
with $R''_M$ twisted PDO, and thus,
\begin{eqnarray}
[\tilde P, \tilde\Pi_{M+1}]&=& [\boldsymbol{Q},\tilde \Pi_{M+1}]+hR''_M\nonumber\\
&=& i(\tilde\Pi_0S_{M}\tilde\Pi_0^\perp + \tilde\Pi_0^\perp S_{M}\tilde\Pi_0)\nonumber\\
\label{comP}
&& \hskip 0.6cm+ \tilde\Pi_0^\perp [\boldsymbol{Q},T_{M}]\tilde\Pi_0^\perp - \tilde\Pi_0 [\boldsymbol{Q},T_{M}]\tilde\Pi_0 + hR_M^{(3)}
\end{eqnarray}
with $R_M^{(3)}$ twisted PDO,
and, using the hypothesis of induction (and, again, the twisted symbolic calculus),
\begin{eqnarray}
&& \tilde\Pi_0^\perp [\boldsymbol{Q},T_{M}]\tilde\Pi_0^\perp \nonumber \\
&&\hskip 0.5cm = \tilde\Pi_0^\perp [\tilde P,T_{M}] \tilde\Pi_0^\perp +hR_M^{(4)}\nonumber\\
&&\hskip 0.5cm = h^{-(M+1)}\tilde\Pi_0^\perp [\tilde P,(\Pi^{(M)})^2-\Pi^{(M)}]\tilde\Pi_0^\perp +hR_M^{(4)}\nonumber\\
&&\hskip 0.5cm = h^{-(M+1)}\tilde\Pi_0^\perp ([\tilde P,\Pi^{(M)}]\Pi^{(M)} + \Pi^{(M)}[\tilde P,\Pi^{(M)}]- [\tilde P,\Pi^{(M)}]) \tilde\Pi_0^\perp +hR_M^{(4)}\nonumber\\
&&\hskip 0.5cm =-i\tilde\Pi_0^\perp (S_M\Pi^{(M)} + \Pi^{(M)}S_M- S_M) \tilde\Pi_0^\perp +hR_M^{(4)}\nonumber\\
\label{comQ1}
&&\hskip 0.5cm =i\tilde\Pi_0^\perp S_M\tilde\Pi_0^\perp  + hR_M^{(5)},
\end{eqnarray}
and, in the same way,
\be
\label{comQ2}
\tilde\Pi_0 [\boldsymbol{Q},T_{M}]\tilde\Pi_0=-i\tilde\Pi_0 S_M\tilde\Pi_0  + hR_M^{(6)},
\ee
where the operators $R_M^{(k)}$'s are all twisted PDO's. Inserting (\ref{comQ1})-(\ref{comQ2}) into (\ref{comP}), we finally obtain,
$$
[\tilde P, \tilde\Pi_{M+1}]=iS_M +hR_M^{(7)},
$$
that implies,
$$
[\tilde P, \Pi^{(M+1)}]=-ih^{M+2} S_{M+1},
$$
where $S_{M+1}$ is a twisted PDO. Therefore, the induction is established.
\vskip 0.3cm
From this point, we follow an idea of \cite{So}. Let $g\in C_0^\infty (\R)$. Using Propositions \ref{regul} and \ref{fdeP}, and writing $g(\tilde P)\tilde\Pi_k =g(\tilde P)(\tilde P+i)^N(\tilde P+i)^{-N}\tilde\Pi_k$, we see that the operators $g(\tilde P)\tilde\Pi_k$ ($k\geq 0$) are all twisted $h$-admissible operators. In particular, they are all bounded, uniformly with respect to $h$. Moreover,  for any $ \ell, \ell' \geq 0$, any $N\geq 1$, and any  functions ${\i}_1,\cdots,{\i}_N\in  C_b^\infty (\R^n)$, by construction, $h^{-N}{\rm ad}_{{\i}_1}\circ\cdots\circ {\rm ad}_{{\i}_N}(\tilde\Pi_k)$
is a twisted PDO, and thus, by Propositions \ref{regul} and \ref{fdeP}, $h^{-N}\tilde P^\ell g(\tilde P){\rm ad}_{{\i}_1}\circ\cdots\circ {\rm ad}_{{\i}_N}(\tilde\Pi_k)\tilde P^{\ell'}$ is uniformly bounded. It is also easy to show (e.g., by using (\ref{funcreg}) hereafter) that,
 \be
 \label{adgP}
 \tilde P^\ell{\rm ad}_{{\i}_1}\circ\cdots\circ {\rm ad}_{{\i}_N}(g(\tilde P))\tilde P^{\ell'} ={\mathcal  O}(h^N),
 \ee
  and therefore, we obtain,
$$
h^{-N}\tilde P^\ell{\rm ad}_{{\i}_1}\circ\cdots\circ {\rm ad}_{{\i}_N}(g(\tilde P)\tilde\Pi_k)\tilde P^{\ell'}={\mathcal  O}(1),
$$
uniformly with respect to $h$. As a consequence, 
we can resum in a standard way the formal series of operators $\sum_{k= 0}^\infty h^kg(\tilde P)\tilde\Pi_k$ (see, e.g., \cite{Ma2} Lemma 2.3.3), in such a way that, if we denote by $\Pi (g)$ such a resummation,  we have,
\be
\label{commg2}
\Vert \tilde P^\ell{\rm ad}_{{\i}_1}\circ\cdots\circ {\rm ad}_{{\i}_N}(\Pi (g)-\sum_{k=0}^{M-1}h^kg(\tilde P)\tilde\Pi_k)\tilde P^{\ell'}\Vert_{{\mathcal  L}(L^2(\R^n;{\mathcal  H}))}  ={\mathcal  O}(h^{M+N}),
\ee
for any $\ell, \ell'\geq 0$, $M, N\geq 0$ and any ${\i}_1,\cdots,{\i}_N\in  C_b^\infty (\R^n)$ (with the conventions ${\rm ad}_{{\i}_1}\circ\cdots\circ {\rm ad}_{{\i}_N}(\Pi (g))=\Pi (g)$ if $N=0$, and $\sum_{k=0}^{M-1}=0$ if $M=0$).
\vskip 0.2cm
Then, we prove,
\begin{lemma}\sl 
\label{commg} For any $\ell\geq 0$, one has,
\be
\label{commg1}
\Vert \tilde P^\ell(\Pi (g)-\Pi (g)^*)\Vert_{{\mathcal  L}(L^2(\R^n;{\mathcal  H}))} ={\mathcal  O}(h^\infty ).
\ee
\end{lemma}
{\em Proof -- } In view of (\ref{commg2}),  it is enough to show that, for any $M\geq 1$, one has,
\be
(\tilde P+i)^\ell [g(\tilde P),\Pi^{(M)}] ={\mathcal  O}(h^{M+1}).
\ee
For $N\geq 1$ large enough, we set $g_N(s):=g(s)(s+i)^N \in C_0^\infty (\R)$, and we observe that,
\be
\label{funcreg}
g(\tilde P) = g_N(\tilde P)(\tilde P +i)^{-N} =\frac1{\pi}\int{\overline{\partial}}\tilde g_N(z)
(\tilde P-z)^{-1} (\tilde P +i)^{-N} dz\; d\bar z,
\ee
where $\tilde g_N$ is an almost analytic extension of $g_N$.
Therefore, we obtain,
\vskip 0.3cm\noindent
$(\tilde P+i)^\ell [g(\tilde P),\Pi^{(M)}]$
\be
\label{regul2}
 = \frac1{\pi}\int{\overline{\partial}}\tilde g_N(z)
(\tilde P-z)^{-1} (\tilde P +i)^{\ell -N}[\Pi^{(M)} , (\tilde P-z)(\tilde P+i)^N]  (\tilde P-z)^{-1} (\tilde P +i)^{-N} dz\; d\bar z,
\ee
and it follows from (\ref{comform}) and the twisted PDO calculus, that, 
\be
\label{regul3}
[\Pi^{(M)} , (\tilde P-z)(\tilde P+i)^N]=h^{M+1}R_{M,N}
\ee
where $R_{M,N}$ is a twisted PDO of degree $\mu_M + mN$,  with $\mu_M$ the degree of $S_M$. Therefore, if we choose $N$ such that $2mN-m\ell \geq \mu_M + mN$, that is, $N\geq \ell+\mu_M/m$, then (\ref{regul2})-(\ref{regul3}) and  Proposition  \ref{regul} tell us that $h^{-(M+1)}[g(\tilde P),\Pi^{(M)}]$ is a twisted $h$-admissible operator, and the result follows.
\hfill$\bullet$
\vskip 0.3cm
We set,
\be
\label{defpigtilde}
\tilde \Pi_g := \Pi (g) +\Pi (g)^*-\frac12(g(\tilde P))\Pi (g)^* + \Pi (g)g(\tilde P))+ (1-g(\tilde P))\tilde\Pi_0 (1-g(\tilde P)).
\ee
Then, $\tilde \Pi_g$ is a selfadjoint twisted $h$-admissible operator, and since $\Pi (g) = g(\tilde P) \tilde\Pi_0 + {\mathcal  O}(h)$, we have,
\be
\label{oh}
\Vert\tilde \Pi_g - \tilde\Pi_0\Vert_{{\mathcal  L}(L^2(\R^n;{\mathcal  H}))} + \Vert\tilde \Pi_g^2 - \tilde \Pi_g\Vert_{{\mathcal  L}(L^2(\R^n;{\mathcal  H}))} = {\mathcal  O}(h).
\ee
 By construction, we also have $\tilde P^\ell (g(\tilde P)\Pi (g)^*-\Pi (g)g(\tilde P))= {\mathcal  O}(h^\infty)$ for all $\ell\geq 0$, and thus, by Lemma \ref{commg},
\be
\label{simp1}
\tilde P^\ell\tilde \Pi_g =\tilde P^\ell\left[ \Pi (g) + (1-g(\tilde P))\left(\Pi(g) + \tilde\Pi_0(1-g(\tilde P))\right)\right] + {\mathcal  O}(h^\infty).
\ee
Moreover, if $f\in C_0^\infty (\R)$ is such that $\supp f\subset \{ g=1\}$, and if we denote by $\Pi(f )$ a resummation of the formal series $\sum_{k\geq 0}h^kf (\tilde P)\tilde\Pi_k$ as before, since $f(\tilde P)(1-g(\tilde P))=0$, $ f(\tilde P)\Pi (g) - \Pi (f)= {\mathcal  O}(h^\infty)$, and $\tilde P^\ell(1-g(\tilde P)\Pi(g)f(\tilde P) =\tilde P^\ell (1-g(\tilde P)\Pi(g)^*f(\tilde P)+ {\mathcal  O}(h^\infty) =\tilde P^\ell(1-g(\tilde P)\Pi(f)+ {\mathcal  O}(h^\infty) ={\mathcal  O}(h^\infty)$, we deduce from (\ref{simp1}) and Lemma \ref{commg},
$$
\tilde P^\ell [f(\tilde P), \tilde \Pi_g] =  \tilde P^\ell\left(\Pi (f) - \Pi(g)^*f(\tilde P)\right)+ {\mathcal  O}(h^\infty) 
=  \tilde P^\ell \left(\Pi (f) - \Pi(f)^*\right)+ {\mathcal  O}(h^\infty),
$$
and thus,
\be
\label{commt}
 \Vert  \tilde P^\ell [f(\tilde P),\tilde \Pi_g ]\Vert_{{\mathcal  L}(L^2(\R^n;{\mathcal  H}))} ={\mathcal  O}(h^\infty ).
\ee
On the other hand, we deduce from Lemma \ref{commg} and (\ref{projform}),
\begin{eqnarray}
\tilde P^\ell (\Pi (g)^2-\Pi (g^2)) &=& \tilde P^\ell(\Pi(g)\Pi(g)^*-\Pi (g^2))+ {\mathcal  O}(h^\infty) \nonumber\\
&=&\tilde P^\ell( \Pi (g )g(\tilde P)-\Pi (g^2))+ {\mathcal  O}(h^\infty)\nonumber\\
\label{simp2}
&=& {\mathcal  O}(h^\infty),
\end{eqnarray}
and thus, using (\ref{simp1})-(\ref{simp2}),
\be
\label{progsf}
\tilde P^\ell( \tilde \Pi_g^2 -\tilde \Pi_g )f(\tilde P) ={\mathcal  O}(h^\infty).
\ee
Then, following the arguments of \cite{Ne1,Ne2,NeSo,So}, for $h$ small enough we can define the following orthogonal projection:
\be
\label{pig}
\Pi_g :=\frac{1}{2i\pi}\int_{\vert z -1\vert=\frac{1}{2}}(\tilde \Pi_g -z)^{-1}\;dz,
\ee
and it verifies (see \cite{So}, Formula (3.9), and \cite{Ne1}, Proposition 3),
\be
\label{nen}
\Pi_g -\tilde \Pi_g = \frac1{2i\pi }(\tilde \Pi_g^2 -\tilde \Pi_g)\int_{\vert z -1\vert=\frac{1}{2}}(\tilde \Pi_g -z)^{-1}
(2\tilde \Pi_g - 1)(1-\tilde \Pi_g -z)^{-1}(1-z)^{-1}\;dz.
\ee
In particular, we obtain from (\ref{progsf}) and (\ref{nen}),
\be
\tilde P^\ell (\Pi_g -\tilde \Pi_g)f(\tilde P)={\mathcal  O}(h^\infty),
\ee
and thus, we deduce from (\ref{oh}) and (\ref{commt})
that
(\ref{eq:pi-pi0}) and (\ref{eq:comm}) hold.
\vskip 0.3cm
In order to prove (\ref{adpi}), we first observe that, by using (\ref{adgP}), (\ref{commg2}) and the fact that ${\rm ad}_{{\i}_k}(\tilde\Pi_0)=0$, we obtain,
\be
\label{ptildead}
 \tilde P^\ell{\rm ad}_{{\i}_1}\circ\cdots\circ {\rm ad}_{{\i}_N}(\tilde\Pi_g) ={\mathcal  O}(h^N),
\ee
for any $N\geq 1$. On the other hand, we have,
\begin{lemma}\sl  
\label{commPpitilde}For any $\ell \geq 0$ and $z\in\C$ such that $| z-1| = 1/2$,
the operator $\tilde P^\ell (\tilde\Pi_g -z)^{-1}(\tilde P+i)^{-\ell}$ is uniformly bounded on $L^2(\R^n ; {\mathcal  H})$.
\end{lemma}
{\em Proof -- } Writing, for $\ell >0$,
\begin{eqnarray*}
H_\ell :&=& (\tilde P+i)^\ell (\tilde\Pi_g -z)^{-1}(\tilde P+i)^{-\ell}\\
&=&H_{\ell -1}+ (\tilde P+i)^{\ell -1}[\tilde P, (\tilde\Pi_g -z)^{-1}](\tilde P+i)^{-\ell}\\
&=& H_{\ell -1} + H_{\ell -1}(\tilde P+i)^{\ell -1}[\tilde \Pi_g, \tilde P](\tilde P+i)^{-\ell }H_\ell,
\end{eqnarray*}
and performing an easy induction, we see that it is enough to prove that $(\tilde P+i)^{\ell -1}[\tilde \Pi_g, \tilde P](\tilde P+i)^{-\ell }$ is ${\mathcal  O}(h)$. Due to (\ref{simp1}), it is enough to study the two terms $(\tilde P+i)^{\ell -1}[\tilde \Pi(g), \tilde P](\tilde P+i)^{-\ell }$ and $(\tilde P+i)^{\ell -1}[\tilde \Pi_0, \tilde P](\tilde P+i)^{-\ell }$. By (\ref{comform}), the first one is ${\mathcal  O}(h^\infty )$, while the second one is equal to $(\tilde P+i)^{\ell -1}[\tilde \Pi_0, \boldsymbol{\omega}](\tilde P+i)^{-\ell }$ and thus, by Propositions \ref{omegaPDO} and \ref{regul}, is ${\mathcal  O}(h)$.
\hfill$\bullet$
\vskip 0.2cm
Combining (\ref{ptildead}), (\ref{pig}) and Lemma \ref{commPpitilde}, we easily obtain (\ref{adpi}), and this completes the proof of Theorem \ref{th:constrPi}.\hfill$\bullet$
\vskip 0.3cm
\begin{remark}\sl  Observe that the previous proof also provides a way of computing the full symbol of $\tilde\Pi_g$ (and thus of $\Pi_g$, too) up to ${\mathcal  O}(h^M)$, for any $M\geq 1$. Indeed, formulas (\ref{projform}), (\ref{comform}), and (\ref{pim+1}) permit to do it inductively. 
\end{remark}
\begin{remark}\sl For this proof, we did not succeed in adapting the elegant argument of \cite{Sj2} (as this was done for smooth interactions in \cite{So}), because of a technical problem. Namely, this argument involves a translation in the spectral variable $z$, of the type $z\mapsto z + \omega(x,\xi)$, inside the symbol of the resolvent of $\tilde P$. In our case, this would have led to consider a symbol $\tilde a=(\tilde a_j)_{0\leq j\leq r}$ of the type $\tilde a_j =a_j(x, \xi, z +\omega_j(x,\xi))$, where $\omega_j$ is the symbol of $U_j \boldsymbol{\omega}U_j^{-1}$ and $a(x,\xi ,z)=(a_j(x,\xi ,z))_{0\leq j\leq r}$ is the symbol of $(z-\tilde P)^{-1}$. But then, it is not clear to us (and probably may be wrong) that the compatibility conditions (\ref{compatib}) are verified by $\tilde a$, and this prevents us from quantizing it in order to continue the argument.
\end{remark}
\chapter{Decomposition of the Evolution for the Modified Operator}
\label{decmodop}
\setcounter{equation}{0}%
\setcounter{theorem}{0}%
In this chapter we prove a general result on the quantum evolution of $\tilde P$.
\begin{theorem}\sl 
\label{th2}
Under the same assumtions as for Theorem \ref{th:constrPi}, let $g\in C_0^\infty (\R)$. Then, one has the following results:\\
1) Let $\varphi_0\in L^2(\R^n ;{\mathcal  H})$  verifying,
\be
\label{locen1}
\varphi_0=f(\tilde P)\varphi_0,
\ee
for some $f\in C_0^\infty (\R)$ such that $\supp f \subset \{ g=1\}$.  Then, with the projection $\Pi_g$ constructed in Theorem \ref{th:constrPi}, one has,
\be
\label{reduc1}
e^{-it\tilde P/h} \varphi_0  =  e^{-it\tilde P^{(1)}/h}\Pi_g\varphi_0 + e^{-it\tilde P^{(2)}/h}(1-\Pi_g)\varphi_0 +
{\mathcal  O}(\vert t\vert h^\infty  \Vert\varphi_0\Vert)
\ee
uniformly with respect to $h$ small enough, $t\in \R$ and $\varphi_0$ verifying (\ref{locen1}), with,
$$
\tilde P^{(1)}: =\Pi_g \tilde P\Pi_g \quad ;\quad \tilde P^{(2)}: =(1-\Pi_g) \tilde P(1-\Pi_g).
$$
2) Let $\varphi_0\in L^2(\R^n ;{\mathcal  H})$ (possibly $h$-dependent) verifying $\Vert \varphi_0\Vert =1$, and,
\be
\label{locen2}
\varphi_0=f(\tilde P)\varphi_0 + {\mathcal  O}(h^\infty ),
\ee
for some $f\in C_0^\infty (\R)$ such that $\supp f \subset \{ g=1\}$.  Then,  one has,
\be
\label{reduc2}
e^{-it\tilde P/h} \varphi_0  =  e^{-it\tilde P^{(1)}/h}\Pi_g\varphi_0 + e^{-it\tilde P^{(2)}/h}(1-\Pi_g)\varphi_0 +
{\mathcal  O}(\la t\ra h^\infty  )
\ee
uniformly with respect to $h$ small enough and $t\in \R$.\\
3)  There exists a bounded operator ${\mathcal  W} : L^2(\R^n ;{\mathcal  H})\rightarrow  L^2(\R^n)^{\oplus L}$ with the following properties:
\begin{itemize}
\item For any $j\in \{ 0,1,\dots,r\}$, and any $\varphi_j\in C_d^\infty (\Omega_j)$, the operator ${\mathcal  W}_j : = {\mathcal  W}U_j^{-1}\varphi_j$ is an $h$-admissible operator from $L^2(\R^n ;{\mathcal  H})$ to $L^2(\R^n)^{\oplus L}$;
\item ${\mathcal  W}{\mathcal  W}^* = 1$ and ${\mathcal  W}^*{\mathcal  W}=\Pi_g$;
\item The operator $A:={\mathcal  W}\tilde P{\mathcal  W}^*={\mathcal  W}\tilde P^{(1)}{\mathcal  W}^*$ is an $h$-admissible operator on $L^2(\R^n)^{\oplus L}$ with domain $H^m(\R^n)^{\oplus L}$, and its  symbol $a(x,\xi ;h )$ verifies,
$$
a(x,\xi ;h) = \omega (x,\xi ;h){\bf I}_L+ {\mathcal  M}(x) + \zeta (x)W(x){\bf I}_L+hr(x,\xi ;h )$$
where ${\mathcal  M}(x)$ is a $L\times L$ matrix depending smoothly on $x$, with spectrum $\{\tilde\lambda_{L'+1}(x),\dots,\tilde\lambda_{L'+L}(x)\}$, and $r(x,\xi :h)$ verifies,
$$
\partial^\alpha r (x,\xi ;h ) ={\mathcal  O}(\langle \xi\rangle^{m-1})
$$
for any multi-index $\alpha$ and uniformly with respect to $(x,\xi )\in T^*\R^{n}$ and $h>0$ small enough. 
\end{itemize}
In particular, ${\mathcal  W}\left\vert_{{\rm Ran}\hskip 1pt\Pi_g}\right. \; :\; {\rm Ran}\hskip 1pt\Pi_g\rightarrow  L^2(\R^n)^{\oplus L}$ is unitary, and $e^{-it\tilde P^{(1)}/h}\Pi_g = {\mathcal  W}^*e^{-itA/h} {\mathcal  W}\Pi_g={\mathcal  W}^*e^{-itA/h} {\mathcal  W}$ for all $t\in \R$.
\end{theorem}
\begin{remark}\sl  In Chapter \ref{expexpl}, we give a way of computing easily the expansion of $A$ up to any power of $h$. As an example, we compute explicitly its first three terms (that is, up to ${\mathcal  O}(h^4)$).
\end{remark}
{\em Proof -- }
1) Setting $\varphi := e^{-it\tilde P/h}\varphi_0$, we have $f(\tilde P)\varphi =\varphi$, and thus
\be
\label{evol1}
ih\partial_t \Pi_g \varphi =\Pi_g\tilde Pf(\tilde P)\varphi = \Pi_g^2\tilde Pf(\tilde P)\varphi.
\ee
Moreover, writing $[\Pi_g, \tilde P]f(\tilde P) = [\Pi_g, \tilde Pf(\tilde P)] + \tilde P[f(\tilde P),\Pi_g]$, Theorem \ref{th:constrPi} tells us that $\Vert [\Pi_g, \tilde P]f(\tilde P)\Vert ={\mathcal  O}(h^\infty)$. Therefore, we obtain from (\ref{evol1}),
$$
\label{evol2}
ih\partial_t \Pi_g \varphi =\Pi_g\tilde P\Pi_gf(\tilde P)\varphi +{\mathcal  O}(h^\infty\Vert\varphi\Vert ) = \tilde P^{(1)}\Pi_g\varphi +{\mathcal  O}(h^\infty\Vert\varphi_0\Vert ),
$$
uniformly with respect to $h$ and $t$. This equation can be re-written as,
$$
ih\partial_t (e^{it\tilde P^{(1)}/h}\Pi_g \varphi) = {\mathcal  O}(h^\infty\Vert\varphi_0\Vert ),
$$
and thus, integrating from $0$ to $t$, we obtain,
$$
\Pi_g \varphi = e^{-it\tilde P^{(1)}/h}\Pi_g \varphi_0 +{\mathcal  O}(|t|h^\infty\Vert\varphi_0\Vert ),
$$
uniformly with respect to $h$, $t$ and $\varphi_0$. 
\vskip 0.3cm
Reasoning in the same way with $1-\Pi_g$ instead of $\Pi_g$, we also obtain,
$$
(1-\Pi_g )\varphi = e^{-it\tilde P^{(2)}/h}(1-\Pi_g) \varphi_0 +{\mathcal  O}(|t|h^\infty\Vert\varphi_0\Vert ),
$$
and (\ref{reduc1}) follows.
\vskip 0.2cm
2) Formula (\ref{reduc2}) follows exactly in the same way.
\vskip 0.2cm
3) Since $\Pi_g -\tilde \Pi_0 ={\mathcal  O}(h)$, for $h$ small enough we can consider the operator ${\mathcal  V}$ defined by the Nagy formula,
\be
\label{defV}
{\mathcal  V}= \left( \tilde \Pi_0\Pi_g +(1-\tilde \Pi_0 )(1-\Pi_g )\right)\left( 1-(\Pi_g - \tilde \Pi_0)^2\right)^{-1/2}.
\ee
Then, ${\mathcal  V}$ is a twisted $h$-admissible operator, it differs from the identity by ${\mathcal  O}(h)$, and  standard computations (using that $(\Pi_g - \tilde \Pi_0)^2$ commutes with both $\tilde \Pi_0\Pi_g$ and $(1-\tilde \Pi_0 )(1-\Pi_g )$: see, e.g.,   \cite{Ka} Chap.I.4) show that,
$$
{\mathcal  V}^*{\mathcal  V} = {\mathcal  V}{\mathcal  V}^* = 1\quad {\rm and\quad}  \tilde \Pi_0{\mathcal  V} = {\mathcal  V}\Pi_g .
$$
Now, with $\tilde u_k$ as in Lemma \ref{sectionscinf}, we define $Z_L : L^2(\R^n ;{\mathcal  H})\rightarrow  L^2(\R^n)^{\oplus L}$ by,
$$
Z_L\psi (x)=\bigoplus_{k=L'+1}^{L'+L}\langle \psi (x),
 \tilde u_k(x)\rangle_{{\mathcal  H}},
$$
and we set,
\be
\label{defW}
{\mathcal  W}:= Z_L \circ {\mathcal  V} =Z_L +{\mathcal  O}(h).
\ee
Thanks to the properties of ${\mathcal  V}$, we see that  ${\mathcal  W}\Pi_g ={\mathcal  W}$, and, since $Z_L^*Z_L =\tilde \Pi_0$ and $Z_LZ_L^* =1$, we also obtain:
$$
{\mathcal  W}^*{\mathcal  W}={\mathcal  V}^*\tilde \Pi_0{\mathcal  V} = \Pi_g\;\; ;\;\; {\mathcal  W}{\mathcal  W}^* = 1.
$$
Moreover, for any $\varphi_j, \i_j\in C_d^\infty (\Omega_j)$ such that $\i_j =1$ near $\supp \varphi_j$,
and for any $\psi\in L^2(\R^n ;{\mathcal  H})$, we have,
$$
{\mathcal  W}U_j^{-1}\varphi_j \psi (x)= \bigoplus_{k=L'+1}^{L'+L}\langle {\mathcal  V}_j\psi (x),
 \tilde u_{k,j}(x)\rangle_{{\mathcal  H}},
 $$
 with ${\mathcal  V}_j:=U_j \i_j{\mathcal  V}U_j^{-1}\varphi_j$ and $\tilde u_{k,j}(x):=U_j(x)\tilde u_k(x) \in C^\infty (\Omega_j, {\mathcal  H})$. Therefore, ${\mathcal  W}U_j^{-1}\varphi_j$ is an $h$-admissible operator from $L^2(\R^n ;{\mathcal  H})$ to $L^2(\R^n)^{\oplus L}$, and the first two properties stated on ${\mathcal  W}$ are proved. (Actually, one can easily see that ${\mathcal  W}$ also verifies a property analog to the first one in Proposition \ref{eqdef}, and thus, with an obvious extension of the notion of twisted operator, that  ${\mathcal  W}$ is, indeed,  a twisted $h$-admissible operator from $L^2(\R^n ;{\mathcal  H})$ to $L^2(\R^n)^{\oplus L}$.)
 \vskip 0.3cm
Then, defining
\be
\label{defA}
A:=  {\mathcal  W}\tilde P   {\mathcal  W}^*={\mathcal  W}\tilde P^{(1)}   {\mathcal  W}^*,
\ee
we want  to prove that $A$ is an $h$-admissible operator and study its symbol. We first need the following result:
\begin{lemma}\sl 
For any $\ell \geq 0$, any $N\geq 1$ and any ${\i}_1,\cdots,{\i}_N\in  C_b^\infty (\R^n)$, one has,
\be
\label{adW}
\Vert \tilde P^\ell{\rm ad}_{{\i}_1}\circ\cdots\circ {\rm ad}_{{\i}_N}({\mathcal  W}^*)\Vert_{{\mathcal  L}(L^2(\R^n); L^2(\R^n;{\mathcal  H})}  ={\mathcal  O}(h^N).
\ee
\end{lemma}
{\em Proof -- } Since ${\mathcal  W}^*={\mathcal  V}^*Z_L^*$ and $Z_L^*$ commutes with the multiplication by any function of $x$, it is enough to prove,
$$
\tilde P^\ell{\rm ad}_{{\i}_1}\circ\cdots\circ {\rm ad}_{{\i}_N}({\mathcal  V}^*) ={\mathcal  O}(h^N),
$$
on $L^2(\R^n;{\mathcal  H})$. Moreover, using (\ref{adpi}) and and the fact that $\tilde \Pi_0$ commutes with the multiplication by any function of $x$, too, we see on  (\ref{defV}) that it is enough to show  that,
\begin{eqnarray}
\label{commrac}
&&( \tilde P+i)^\ell (1-(\Pi_g - \tilde \Pi_0)^2)^{-1/2}(\tilde P+i)^{-\ell} ={\mathcal  O}(1);\\
\label{adrac}
&& \tilde P^\ell{\rm ad}_{{\i}_1}\circ\cdots\circ {\rm ad}_{{\i}_N}\left( (1-(\Pi_g - \tilde \Pi_0)^2)^{-1/2}\right) ={\mathcal  O}(h^N).
\end{eqnarray}
By construction, we have $\tilde P^\ell (\Pi(g) - g(\tilde P)\tilde\Pi_0) ={\mathcal  O}(h)$, and thus, we immediately see on (\ref{simp1}) that $\tilde P^\ell (\tilde\Pi_g - \tilde\Pi_0) ={\mathcal  O}(h)$. Then, writing
$$
\Pi_g - \tilde\Pi_0 =\frac{1}{2i\pi}\int_{\vert z -1\vert=\frac{1}{2}}(\tilde \Pi_g -z)^{-1}(\tilde\Pi_0 - \tilde\Pi_g)(\tilde \Pi_0 -z)^{-1}\;dz,
$$
and using Lemma \ref{commPpitilde}, we also obtain,
\be
\label{plpi}
\tilde P^\ell (\Pi_g - \tilde\Pi_0) ={\mathcal  O}(h),
\ee
for all $\ell \geq 0$.
In particular, $( \tilde P+i)^\ell (\Pi_g - \tilde \Pi_0)( \tilde P+i)^{-\ell}={\mathcal  O}(h)$, and therefore, for $h$ sufficiently small, we can write,
$$
( \tilde P+i)^\ell (1-(\Pi_g - \tilde \Pi_0)^2)^{-1/2}(\tilde P+i)^{-\ell} =\left(1-[(\tilde P+i)^{\ell}(\Pi_g - \tilde \Pi_0)(\tilde P+i)^{-\ell}]^2\right)^{-1/2},
$$
and (\ref{commrac}) follows.
\vskip 0.3cm
To prove (\ref{adrac}), we write $(1-(\Pi_g - \tilde \Pi_0)^2)^{-1/2}$ as,
$$
(1-(\Pi_g - \tilde \Pi_0)^2)^{-1/2}=1 +\sum_{k=1}^{\infty} \alpha_k(\Pi_g -\tilde\Pi_0)^{k},
$$
where the radius of convergence of the power series $\sum_{k=1}^{\infty} \alpha_kz^k$ is 1. Thus,
$$
\tilde P^{\ell}{\rm ad}_{{\i}_1}\circ\cdots\circ {\rm ad}_{{\i}_N}\left( (1-(\Pi_g - \tilde \Pi_0)^2)^{-1/2}\right) =\sum_{k=1}^{\infty} \alpha_k{\mathcal  A}_{N,k}
$$ 
where ${\mathcal  A}_{N,k}:=\tilde P^{\ell}{\rm ad}_{{\i}_1}\circ\cdots\circ {\rm ad}_{{\i}_N}((\Pi_g -\tilde\Pi_0)^{k})$ is the sum of $k^N$ terms of the form,
$$
\tilde P^{\ell}[{\rm ad}_{{\i}_{i_{1,1}}}\cdots {\rm ad}_{{\i}_{i_{1,n_1}}}(\Pi_g -\tilde\Pi_0)]\dots [{\rm ad}_{{\i}_{i_{k,1}}}\cdots {\rm ad}_{{\i}_{i_{k, n_k}}}(\Pi_g -\tilde\Pi_0)],
$$
with $n_1, \dots ,n_k\geq 0$, $n_1 +\dots +n_k=N$. Then, using (\ref{commg2}) together with (\ref{plpi}), we see that all these terms have a norm bounded by $(C_N)^kh^{k+N}$, for some constant $C_N >0$ independent of $k$. Therefore, $\Vert {\mathcal  A}_{N,k}\Vert \leq k^N(C_N)^kh^{k+N}$, and (\ref{adrac}) follows.\hfill$\bullet$
\vskip 0.3cm
Then, proceeding as in the proof of Lemma \ref{suppdisj}, we deduce from Lemma \ref{adW} that, if ${\i} , \psi\in  C_b^\infty (\R^n)$ are such that ${\rm dist} \left( \supp {\i}, \supp\psi \right) >0$, then,  $\Vert \tilde P^\ell {\i} {\mathcal  W}^*\psi\Vert = {\mathcal  O}(h^\infty )$.
As a consequence, taking a partition of unity $(\i_j)_{j=0,\dots,r}$ on $\R^n$ with $\i_j\in C_d^\infty (\Omega_j)$, and choosing $\varphi_j\in C_d^\infty (\Omega_j)$ such that $\dist (\supp (\varphi_j -1), \supp\i_j) >0$ ($j=0,\dots,r$),
we have (using also that $\tilde P$ is local in the variable $x$),
$$
A= \sum_{j=0}^r {\mathcal  W}\i_j\tilde P   {\mathcal  W}^*=\sum_{j=0}^r \varphi_j {\mathcal  W}\i_j\tilde P \varphi_j^2  {\mathcal  W}^*\varphi_j+R(h),
$$ 
with $\Vert R(h)\Vert_{{\mathcal  L}(L^2(\R^n))} ={\mathcal  O}(h^\infty )$. Thus,
$$
A= \sum_{j=0}^r \varphi_j {\mathcal  W}U_j^{-1}\i_j\tilde P_j U_j\varphi_j {\mathcal  W}^*\varphi_j+R(h),
$$
where $\tilde P_j = U_j \tilde PU_j^{-1}\varphi_j$ is an $h$-admissible (differential) operator from $H^m(\R^n; {\mathcal  D}_Q)$ to $L^2(\R^n ;{\mathcal  H})$, while $U_j\varphi_j {\mathcal  W}^*\varphi_j$ is an $h$-admissible operator from
$H^m(\R^n)^{\oplus L}$ to $H^m(\R^n; {\mathcal  D}_Q)$, and
$\varphi_j {\mathcal  W}U_j^{-1}\i_j$ is an $h$-admissible operator from $L^2(\R^n ;{\mathcal  H})$ to $L^2(\R^n)^{\oplus L}$. 
\vskip 0.3cm
Therefore, $A$ is an $h$-admissible operator from $H^m(\R^n)^{\oplus L}$ to $L^2(\R^n)^{\oplus L}$, and, if we set,
$$
\tilde p_j(x,\xi ;h) = \omega (x,\xi ;h) + \tilde Q_j(x)+ \zeta (x)W(x)+h\sum_{|\beta|\leq m-1}\omega_{\beta ,j}(x;h)\xi^\beta,
$$
and if we denote by $v_j(x,\xi )$ (resp. $v_j^*(x,\xi )$) the symbol of $U_j{\mathcal  V}U_j^{-1})$ (resp. $U_j{\mathcal  V}U_j^{-1}$),
then, the (matrix) symbol $a=(a_{k,\ell})_{1\leq k,\ell\leq L}$ of $A$, is given by,
$$
a_{k,\ell}(x,\xi ,h) = \sum_{j=0}^r \la \i_j(x)v_j(x,\xi )\sharp \tilde p_j(x,\xi)\sharp v_j^*(x,\xi )\sharp \tilde u_{L'+k,j}(x), \tilde u_{L'+\ell,j}(x)\ra_{\mathcal  H}.
$$
In particular, since $\partial^\alpha (v_j -1)$ and $\partial^\alpha (v_j^* -1)$ are ${\mathcal  O}(h)$, we obtain,
\begin{eqnarray*}
 a_{k,\ell}(x,\xi ,h)=  \sum_{j=0}^r \la \i_j(x)(\omega(x,\xi ) + \tilde Q_j(x)+ \zeta (x)W(x)) \tilde u_{L'+k,j}(x), \tilde u_{L'+\ell,j}(x)\ra_{\mathcal  H}\\
+ r_{k,\ell}(h)
\end{eqnarray*}
with $\partial^\alpha r_{k,\ell}(h) ={\mathcal  O}(h\la\xi\ra^{m-1})$, and thus, using the fact that  
$$
\la \tilde Q_j(x)\tilde u_{L'+k,j}(x), \tilde u_{L'+\ell,j}(x)\ra = \varphi_j(x)\la \tilde Q(x)\tilde u_{L'+k}(x), \tilde u_{L'+\ell} (x)\ra,
$$
this finally gives,
\begin{eqnarray*}
a_{k,\ell}(x,\xi ,h) 
&=&  \sum_{j=0}^r\i_j(x)(\omega(x,\xi )\delta_{k,\ell}+m_{k,\ell}(x)+ \zeta (x)W(x)\delta_{k,\ell}) + r_{k,\ell}(h)\\
&=& (\omega(x,\xi )+ \zeta (x)W(x))\delta_{k,\ell} +m_{k,\ell}(x) + r_{k,\ell}(h),
\end{eqnarray*}
with $m_{k,\ell}(x):=\la \tilde Q(x)\tilde u_{L'+k}(x), \tilde u_{L'+\ell }(x)\ra$.
This completes the proof of
Theorem \ref{th2}.\hfill$\bullet$
\chapter{Proof of Theorem \ref{MAINTH}}
\label{mainproof}
\setcounter{equation}{0}%
\setcounter{theorem}{0}%
In view of Theorem \ref{th2}, it is enough to prove,
\begin{theorem}\sl  
\label{th3} Let  $\varphi_0\in L^2(\R^n;{\mathcal  H})$
such that $\Vert \varphi_0\Vert =1$, and,
\be
\label{dinit}
\Vert\varphi_0\Vert_{L^2(K_0^c;{\mathcal  H})}+\Vert (1-\Pi_g )\varphi_0\Vert +\Vert (1-f(P))\varphi_0\Vert ={\mathcal  O}(h^\infty ),
\ee
for some $K_0\subset\subset\Omega'\subset\subset\Omega$, $f,g \in  C_0^{\infty}({\R})$, $gf=f$, and let  $\tilde P$ be  the operator constructed in Chapter \ref{sect2} with $K=\overline{\Omega'}$, and 
$\Pi_g$ be the projection constructed in Theorem \ref{th:constrPi}.
Then, with the notations of Theorem \ref{th2}, we have,
\begin{equation}
\label{sol2}
e^{-itP/h}\varphi_0= {\mathcal  W}^*e^{-itA/h}{\mathcal  W}\varphi_0  +{\mathcal  O}\left( \la t\ra h^\infty  \right),
\end{equation}
uniformly with respect to $h>0$ small enough and  $t\in  [0, T_{\Omega'}(\varphi_0))$.
\end{theorem}
{\bf Proof :} Denote by $\i\in C_0^\infty (\Omega'_{K})$ (where $\Omega'_{K}$ is the same as in  Proposition \ref{qtilde}) a cutoff function such that $\i =1$ on $K$. We first prove,
\begin{lemma}\sl 
\label{compP}
$$
\Vert (f(P) -f(\tilde P))\i \Vert_{{\mathcal  L}(L^2(\R^n;{\mathcal  H})} ={\mathcal  O}(h^\infty ).
$$
\end{lemma}
{\em Proof -- } Using (\ref{calfonc}), we obtain,
$$
(f(P) -f(\tilde P))\i =\frac1{\pi}\int{\overline{\partial}}\tilde f(z)
(P-z)^{-1}(\tilde P - P)(\tilde P-z)^{-1}\i dz\; d\bar z.
$$
Moreover, if $\psi \in C_0^\infty (\Omega'_{K})$ is such that $\psi =1$ on a neighborhood of $\supp \i$, Corollary \ref{pqutruc} and Lemma \ref{suppdisj} tell us,
$$
(\psi -1)(\tilde P-z)^{-1}\i ={\mathcal  O}(h^N\vert\im z\vert^{-(N+1)}),
$$
for any $N\geq 1$. As a consequence,
$$
(f(P) -f(\tilde P))\i =\frac1{\pi}\int{\overline{\partial}}\tilde f(z)
(P-z)^{-1}(\tilde P - P)\psi (\tilde P-z)^{-1}\i dz\; d\bar z + {\mathcal  O}(h^\infty),
$$
and since $(\tilde P - P)\psi =(\tilde Q -Q)\psi =0$, the result follows.\hfill$\bullet$
\vskip 0.2cm
Now, by (\ref{dinit}), we have,
$$
\varphi_0 = f(P)\varphi_0+{\mathcal  O}(h^\infty) = f(P)\i\varphi_0+{\mathcal  O}(h^\infty),
$$
and thus, by Lemma \ref{compP},
$$
\varphi_0  = f(\tilde P)\i\varphi_0+{\mathcal  O}(h^\infty)=  f(\tilde P)\varphi_0+{\mathcal  O}(h^\infty).
$$
This means that (\ref{locen2}) is satisfied, and thus, by Theorem \ref{th2}, the decomposition (\ref{reduc2}) is true. Using (\ref{dinit}) again, this gives,
\be
\label{evpti}
e^{-it\tilde P/h} \varphi_0  =  e^{-it\tilde P^{(1)}/h}\Pi_g\varphi_0  +
{\mathcal  O}(\vert t\vert h^\infty  ) = {\mathcal  W}^*e^{-itA/h} {\mathcal  W}\varphi_0  +
{\mathcal  O}(\la t\ra h^\infty  ),
\ee
uniformly with respect to $h$ and $t$.
\vskip 0.3cm
On the other hand, if we set $\varphi (t) := e^{-itP/h}\varphi_0$, then, by assumption, $\varphi (t) = f(P) \varphi (t) + {\mathcal  O}(h^\infty )$ and $\varphi (t) = \i \varphi (t) + {\mathcal  O}(h^\infty )$ uniformly for $t\in [0,T_{\Omega'}(\varphi_0)]$. Therefore, applying Lemma  \ref{compP} again, we obtain as before,
$\varphi (t) = f(\tilde P) \varphi (t) + {\mathcal  O}(h^\infty )$, and thus also,
\be
\label{locfi}
\varphi (t) = f(\tilde P) \i \varphi (t) + {\mathcal  O}(h^\infty ),
\ee 
uniformly with respect to $h$ and $t\in [0,T_{\Omega'}(\varphi_0)]$.
Moreover, since $P$ and $\tilde P$ coincide on the support of $\i$, we can write,
$$
ih\partial_t f(\tilde P)\i\varphi(t) = f(\tilde P)\i P\varphi (t) = f(\tilde P)\tilde P\i\varphi (t) + f(\tilde P)[\i ,\tilde P]\varphi (t),
$$
and thus, since $f(\tilde P)[\i ,\tilde P]=f(\tilde P)[\i ,\boldsymbol{\omega}]$ is bounded, and $[\i ,\boldsymbol{\omega}]$ is a differential operator with coefficients supported in $\supp\nabla\i$ (where $\varphi$ is ${\mathcal  O}(h^\infty)$), we obtain,
$$
ih\partial_t f(\tilde P)\i\varphi(t) = f(\tilde P)\i P\varphi (t) = \tilde Pf(\tilde P)\i\varphi (t) + {\mathcal  O}(h^\infty).
$$
As  a consequence,
$$
f(\tilde P)\i\varphi(t)=e^{-it\tilde P/h}f(\tilde P)\i\varphi_0+ {\mathcal  O}(| t|h^\infty),
$$
and therefore, by (\ref{locfi}),
\be
\label{compevo}
\varphi(t)=e^{-it\tilde P/h}\varphi_0+ {\mathcal  O}(\la t\ra h^\infty),
\ee
uniformly with respect to $h$ and $t\in [0,T_{\Omega'}(\varphi_0))$.
Then, Theorem \ref{th3} follows from (\ref{evpti}) and (\ref{compevo}).
\hfill$\bullet$
\chapter{Proof of Corollary \ref{CORT}}
\label{proofcor}
\setcounter{equation}{0}%
\setcounter{theorem}{0}%
First of all, let us recall the (standard) notion of frequency set $FS(v)$ of some  (possibly $h$-dependent) $v\in L_{\rm loc}^2(\Omega)$ (see, e.g., \cite{Ma2} and references therein). It is said that a point $(x_0,\xi_0)\in T^*\Omega$ is not in $FS(v)$ if there exist $\i_1 \in C_0^\infty (\omega)$ and $\i_2\in C_0^\infty (\R^n)$ such that $\i_1(x_0)=\i_2(\xi_0) =1$ and $\Vert\i_2(hD_x)\i_1v\Vert_{L^2(\R^n)}={\mathcal  O}(h^\infty)$. This is also equivalent to say that there exists an open neighborhood ${\mathcal  N}$ of $(x_0,\xi_0)$ in $T^*\R^n$, such that, for {\it any} $\i\in C_0^\infty ({\mathcal  N})$ and any $\i_1\in C_0^\infty (\Omega)$, one has $\Vert{\rm Op}_h(\i)\i_1v\Vert_{L^2(\R^n)}={\mathcal  O}(h^\infty)$.
\vskip 0.3cm
As one can see, this notion can be extended in an  obvious way to functions in $ L^2_{\rm loc}(\Omega;{\mathcal  H})$, and  it is easy to see (e.g., as in \cite{Ma2} Section 2.9) that the latter property still holds with operator-valued functions $\i\in C_0^\infty ({\mathcal  N} ;{\mathcal  L}({\mathcal  H}))$, or even more generally, $\i\in C_0^\infty ({\mathcal  N} ;{\mathcal  L}({\mathcal  H};{\mathcal  H}'))$ where ${\mathcal  H}'$ is an arbitrary Hilbert-space. 
\vskip 0.3cm
We first prove,
\begin{lemma}\sl  
\label{FSW}
Let ${\mathcal  W} : L^2(\R^n ;{\mathcal  H})\rightarrow  L^2(\R^n)$ be the operator given in Theorem \ref{th2}. Then, for any  $j\in \{0,1,\dots,r\}$, any $\varphi \in L^2(\R^n ;{\mathcal  H})$ and $v\in L^2(\R^n)$, such that $\Vert\varphi\Vert =\Vert v\Vert =1$, one has,
\begin{eqnarray*}
&&FS({\mathcal  W}\varphi)\cap T^*\Omega_j=FS(U_j\Pi_g\varphi)\cap T^*\Omega_j;\\
&&FS(U_j{\mathcal  W}^*v)\cap T^*\Omega_j=FS(v)\cap T^*\Omega_j.
\end{eqnarray*}
\end{lemma}
{\em Proof -- } Since ${\mathcal  W}{\mathcal  W}^*=1$ and ${\mathcal  W}^*{\mathcal  W}=\Pi_g$, it is enough to prove the two inclusions $FS({\mathcal  W}\varphi)\cap T^*\Omega_j\subset FS(U_j\Pi_g\varphi)\cap T^*\Omega_j$ and $FS(U_j{\mathcal  W}^*v)\cap T^*\Omega_j\subset FS(v)\cap T^*\Omega_j$.
\vskip 0.3cm
Therefore, let $(x_0,\xi_0)\in T^*\Omega_j$, and
assume first that  $(x_0,\xi_0)\notin FS(U_j\Pi_g\varphi)$. In particular, this implies that, if ${\mathcal  N}\subset\subset T^*\Omega_j$ is a small enough neighborhood of $(x_0,\xi_0)$, then $\Vert{\rm Op}_h(\i_1)U_j\Pi_g\varphi\Vert ={\mathcal  O}(h^\infty)$ for all $\i_1\in C_0^\infty ({\mathcal  N}; {\mathcal  L}({\mathcal  H};\C))$. Then, taking $\i\in C_0^\infty ({\mathcal  N})$ and $\psi_j\in C_0^\infty (\Omega_j)$ such that $\psi_j(x)=1$ near $\pi_x(\supp\i )$ and $\i (x_0,\xi_0)=1$, we write,
\begin{eqnarray*}
{\rm Op}_h(\i){\mathcal  W}\varphi &=& {\rm Op}_h(\i){\mathcal  W}\Pi_g\varphi = {\rm Op}_h(\i){\mathcal  W}\psi_j^2\Pi_g\varphi +{\mathcal  O}(h^\infty )\\
&=&{\rm Op}_h(\i){\mathcal  W}U_j^{-1}\psi_jU_j\psi_j\Pi_g\varphi +{\mathcal  O}(h^\infty ),
\end{eqnarray*}
and since ${\rm Op}_h(\i){\mathcal  W}U_j^{-1}\psi_j$ is an $h$-admissible operator  from $L^2(\R^n;{\mathcal  H})$ to $L^2(\R^n)$, with symbol supported  in ${\mathcal  N}$ (that is, modulo ${\mathcal  O}(h^\infty)$ in $C_b^\infty (\R^n; {\mathcal  L}({\mathcal  H};\C))$), we obtain $\Vert {\rm Op}_h(\i){\mathcal  W}\varphi\Vert ={\mathcal  O}(h^\infty)$, and thus $(x_0,\xi_0)\notin FS({\mathcal  W}\varphi)$.
\vskip 0.3cm
Now, assume that $(x_0,\xi_0)\notin FS(v)$. Since $U_j\psi_j{\mathcal  W}^*$ is an $h$-admissible operator, we obtain in the same way that $\Vert {\rm Op}_h(\i)U_j\psi_j{\mathcal  W}^*v\Vert ={\mathcal  O}(h^\infty)$, and thus  $(x_0,\xi_0)\notin FS(U_j{\mathcal  W}^*v)$.\hfill$\bullet$
\vskip 0.3cm
Without loss of generality, we can assume $ T_{\Omega '}(\varphi_0)<+\infty$. By Theorem \ref{th3}, we have,
$$
e^{-itP/h}\varphi_0= {\mathcal  W}^*e^{-itA/h}{\mathcal  W}\varphi_0  +{\mathcal  O}\left(  h^\infty  \right),
$$
uniformly for $t\in [0, T_{\Omega'}(\varphi_0)]$, where ${\mathcal  W}$ and $A$  are given in Theorem \ref{th2}. Thus, by Lemma \ref{FSW}, we immediately obtain,
$$
FS(U_je^{-itP/h}\varphi_0)\cap T^*\Omega_j=FS(e^{-itA/h}{\mathcal  W}\varphi_0)\cap T^*\Omega_j.
$$
On the other hand, since $A$ is an $h$-admissible operator on $L^2(\R^n)$, a well-known result of propagation (see, e.g., \cite{Ma2} Section 4.6, Exercise 12) tells us,
$$
FS(e^{-itA/h}{\mathcal  W}\varphi_0)=\exp tH_{a_0}(FS({\mathcal  W}\varphi_0)).
$$
Therefore, applying Lemma \ref{FSW} again, we obtain,
\be
\label{estfs}
FS(U_je^{-itP/h}\varphi_0)\cap T^*\Omega_j=T^*\Omega_j\cap \exp tH_{a_0}\left(\cup_{k=0}^rFS(U_k\Pi_g\varphi_0)\cap T^*\Omega_k\right).
\ee
By assumption, we also have,
\be
\label{estfsfi1}
\cup_{k=0}^rFS(U_k\Pi_g\varphi_0)=\cup_{k=1}^rFS(U_k\varphi_0)\subset K_0\times\R^n.
\ee
In order to conclude, we need the following result:
\begin{lemma}\sl 
\label{locenergy}
For any  $f\in C_0^\infty (\R)$, $\psi \in C_0^\infty (\R^n)$, $\i_j\in C_0^\infty (\Omega_j)$,  $\varepsilon >0$, and $\rho \in C_b^\infty (\R)$ with $\supp\rho\subset [C_f -\gamma +\varepsilon,+\infty)$ (where $C_f$ is as in Corollary \ref{CORT}),  one has,
$$
\Vert \rho (\i_j\boldsymbol{\omega}\i_j )\psi f(U_j\i_j\tilde PU_j^{-1}\i_j)\Vert ={\mathcal  O}(h^\infty ).
$$
\end{lemma}
{\em Proof -- } We set $\boldsymbol{\omega}_j := \i_j\boldsymbol{\omega}\i_j$ and $\tilde P_j := U_j\i_j\tilde PU_j^{-1}\i_j$. Using Assumptions (H1), (H2), (H4) and Proposition \ref{qtilde}, we see that $\tilde P_j \geq (1-Ch)\boldsymbol{\omega}_j +\gamma -Ch$ for some constant $C>0$ independent of $h$. As a consequence, we have,
$$
 \rho (\boldsymbol{\omega}_j )\tilde P_j \rho (\boldsymbol{\omega}_j )\geq  \rho (\boldsymbol{\omega}_j) ((1-Ch)\boldsymbol{\omega}_j +\gamma -Ch) \rho (\boldsymbol{\omega}_j )\geq (C_f+\varepsilon -C'h)\rho (\boldsymbol{\omega}_j )^2,
 $$
 with $C'= C+CC_f$.
Therefore, we can write,
 $$
 \Vert \rho (\boldsymbol{\omega}_j )\psi f(\tilde P_j)u\Vert^2\leq \frac1{C_f+\varepsilon -C'h}\la \tilde P_j\rho (\boldsymbol{\omega}_j )\psi f(\tilde P_j)u, \rho (\boldsymbol{\omega}_j )\psi f(\tilde P_j)u\ra,
 $$
 for any $u\in L^2(\R^n ;{\mathcal  H})$, and thus,
\begin{eqnarray}
 \Vert \rho (\boldsymbol{\omega}_j )\psi f(\tilde P_j)\Vert &\leq& \frac1{C_f+\varepsilon -C'h}\Vert \tilde P_j\rho (\boldsymbol{\omega}_j )\psi f(\tilde P_j)\Vert\nonumber\\ 
 &\leq&  \frac1{C_f+\varepsilon -C'h}\left(\Vert \rho (\boldsymbol{\omega}_j )\psi \tilde P_jf(\tilde P_j)\Vert + \Vert [\tilde P_j,\rho (\boldsymbol{\omega}_j )\psi ]f(\tilde P_j)\Vert\right).\nonumber\\
  \label{lemfs1}
 && {}
\end{eqnarray}
 Now, on the one hand, since $\supp f$ is included in $[-C_f,C_f]$, we have,
\begin{eqnarray}
 \frac1{C_f+\varepsilon -C'h}\Vert \rho (\boldsymbol{\omega}_j )\psi \tilde P_jf(\tilde P_j)\Vert &=&\frac1{C_f+\varepsilon -C'h}\Vert \tilde P_jf(\tilde P_j)\psi \rho (\boldsymbol{\omega}_j )\Vert\nonumber\\
 \label{lemfs2}
 &\leq& \frac{C_f}{C_f+\varepsilon -C'h}\Vert f(\tilde P_j)\psi \rho (\boldsymbol{\omega}_j )\Vert.
\end{eqnarray}
 On the other hand, since $\tilde P_j$ and $\boldsymbol{\omega}_j $ are both differential operators with respect to $x$ with  smooth (operator-valued) coefficients, and $\rho (\boldsymbol{\omega}_j) \psi$ is  a scalar operator, by standard symbolic calculus, we have,
\be
\label{lemfs3}
 [\tilde P_j,\rho (\boldsymbol{\omega}_j) \psi ]f(\tilde P_j)={\mathcal  O}(h)\rho_1 (\boldsymbol{\omega}_j )\psi_1 f(\tilde P_j)+ {\mathcal  O}(h^\infty ),
 \ee
where $\rho_1\in C_b^\infty (\R)$ and $\psi_1\in C_0^\infty (\R^n)$ are arbitrary functions verifying $\rho_1\rho =\rho$ and $\psi_1\psi =\psi$. Inserting (\ref{lemfs2})-(\ref{lemfs3}) into (\ref{lemfs1}), we obtain,
$$
 \Vert \rho (\boldsymbol{\omega}_j )\psi f(\tilde P_j)\Vert ={\mathcal  O}(h\Vert \rho_1 (\boldsymbol{\omega}_j )\psi_1 f(\tilde P_j)\Vert) + {\mathcal  O}(h^\infty ).
$$
Iterating the procedure, we clearly obtain the lemma.\hfill$\bullet$
 \vskip 0.2cm
 Now, using, e.g., (\ref{locfi}), we know that $e^{-itP/h}\varphi_0 = f(\tilde P)e^{-itP/h}\varphi_0+{\mathcal  O}(h^\infty)$. Moreover, if $\i_j, \psi_j\in C_0^\infty (\Omega_j)$ are such that $\i_j =1$ near $\supp \psi_j$, by Lemma \ref{suppdisj}, we have,
 $$
 U_j\psi_j f(\tilde P)=U_j\psi_j f(\tilde P)\i_j^2 +{\mathcal  O}(h^\infty ) = U_j\psi_j f(\tilde P)U_j^{-1}\i_jU_j\i_j +{\mathcal  O}(h^\infty ),
 $$
 and therefore,
 $$
 U_j\psi_je^{-itP/h}\varphi_0 = U_j\psi_j f(\tilde P)U_j^{-1}\i_jU_j\i_je^{-itP/h}\varphi_0 +{\mathcal  O}(h^\infty ).
 $$
 Then, using lemma \ref{tec2}, we obtain,
 $$
 U_j\psi_je^{-itP/h}\varphi_0=\psi_j f(\tilde P_j)U_j\i_je^{-itP/h}\varphi_0 +{\mathcal  O}(h^\infty ),
 $$
 with $\tilde P_j = U_j\i_j\tilde PU_j^{-1}\i_j$.
Therefore,  using   Lemma \ref{locenergy},  this gives,
$$
\Vert \rho (\i_j\boldsymbol{\omega}\i_j)U_j\psi_j e^{-itP/h}\varphi_0\Vert ={\mathcal  O}(h^\infty ),
$$
and thus, by Lemma \ref{tec1},
\be
\label{estfsfi2}
\Vert \rho (\boldsymbol{\omega})U_j\psi_j e^{-itP/h}\varphi_0\Vert ={\mathcal  O}(h^\infty ).
\ee
Since the principal symbol of $\rho (\boldsymbol{\omega})$ is  $\rho (\omega )$, 
we deduce from (\ref{estfsfi1}), (\ref{estfsfi2}), and standard results on $FS$, that, 
$$
\cup_{k=0}^rFS(U_k\Pi_g\varphi_0)\subset K(f):=\{ (x,\xi )\, ;\, x\in K_0\,, \, \omega (x,\xi) \leq C_f-\gamma\},
$$
and thus, by (\ref{estfs}),
\be
\label{estfsuj}
FS(U_je^{-itP/h}\varphi_0)\cap T^*\Omega_j\subset  \exp tH_{a_0}\left(K(f)\right)\cap T^*\Omega_j,
\ee
for all $t\geq 0$.
\vskip 0.3cm
Then,   for any $j\in\{0,1,\dots,r\}$,  $\psi_j, \tilde\psi_j\in C_0^\infty (\Omega_j)$ with $\tilde\psi_j \psi_j =\psi_j$, and any $\alpha \in C_0^\infty (\R^n)$, we write,
$$
U_j\psi_je^{-itP/h}\varphi_0 = \alpha (hD_x)\tilde\psi_j (x)U_j\psi_je^{-itP/h}\varphi_0 + (1-\alpha (hD_x))U_j\psi_je^{-itP/h}\varphi_0,
$$
and therefore, if $\alpha (\xi )=1$ in a sufficiently large compact set,
$$
U_j\psi_je^{-itP/h}\varphi_0 = \alpha (hD_x)\tilde\psi_j (x)U_j\psi_je^{-itP/h}\varphi_0 + {\mathcal  O}(h^\infty ).
$$
Finally,  if $\supp\tilde\psi_j \cap \pi_x\left( \exp tH_{a_0}\left(K(f)\right)\right) =\emptyset$ (or, more generally, $\supp\tilde\psi_j \cap \pi_x\left(\cup_{k=0}^r \exp tH_{a_0}(FS(U_k\Pi_g\varphi_0))\right) =\emptyset$), then, (\ref{estfs}) and (\ref{estfsuj}) tell us,
$$
\Vert \alpha (hD_x)\tilde\psi_j (x)U_j\psi_je^{-itP/h}\varphi_0\Vert = {\mathcal  O}(h^\infty ),
$$
and thus, by the unitarity of $U_j$,
$$
\Vert \psi_je^{-itP/h}\varphi_0\Vert=\Vert U_j\psi_je^{-itP/h}\varphi_0\Vert = {\mathcal  O}(h^\infty ),
$$
uniformly for $t\in [0,T_{\Omega '}(\varphi_0)]$.
Since we also know that $\Vert e^{-itP/h}\varphi_0\Vert_{K^c}= {\mathcal  O}(h^\infty )$ for some compact set $K\subset \R^n$ (by definition of $T_{\Omega '}(\varphi_0)$), 
this proves that we can actually take for $K$ any compact neighborhood  of $\pi_x\left( \exp tH_{a_0}\left(K(f)\right)\right)$. 
Thus, if $T_{\Omega '}(\varphi_0) < \sup \{ T>0\, ;\, \pi_x(\cup_{t\in [0,T]}\hskip 1pt \exp tH_{a_0}(K(f)))\subset \Omega'\}$, clearly (e.g., by using Theorem \ref{th:appA}),  one can find $T>T_{\Omega '}(\varphi_0)$ and $K_T\subset\subset \Omega'$, such that $\sup_{t\in [0,T]}\Vert e^{-itP/h}\varphi_0\Vert_{K_T^c}={\mathcal  O}(h^{\infty})$. This is in contradiction with the definition of $T_{\Omega '}(\varphi_0)$, and therefore, necessarily, 
$$
T_{\Omega '}(\varphi_0) \geq \sup \{ T>0\, ;\, \pi_x(\cup_{t\in [0,T]}\hskip 1pt \exp tH_{a_0}(K(f)))\subset \Omega'\}.$$
This proves Corollary \ref{CORT}, and also Remark \ref{esttemps2} since, in the last argument, one can replace $K(f)$   by $\cup_{k=0}^r \exp tH_{a_0}(FS(U_k\Pi_g\varphi_0))$ everywhere. \hfill$\bullet$
 \chapter{Computing the  Effective Hamiltonian}
\label{expexpl}
\setcounter{equation}{0}%
\setcounter{theorem}{0}%
Now that we know the existence of an effective Hamiltonian describing the evolution of those states $\varphi_0$ that verify (\ref{condloc}), the problem remains of computing its symbol up to any arbitrary power of $h$ (in Theorem \ref{MAINTH}, only the principal symbol of $A$ is given). Because of the conditions of localization (\ref{condloc}), it is clear that such an effective Hamiltonian is not unique (for instance, the three operators $A$, $Af(A)$ or ${\mathcal  W}f(\tilde P){\mathcal  W}^*A{\mathcal  W}f(\tilde P){\mathcal  W}^*$ could  indifferently be taken). However, its symbol is certainly uniquely determined in the relevant region of the phase space where $\tilde\varphi (t):= {\mathcal  W}e^{-itP/h}\varphi_0$ lives (that is,  on  $FS (\tilde\varphi (t))$ in the sense of the previous chapter, and for $t\in [0, T_{\Omega'}(\varphi_0))$). Therefore, as long as we deal with $h$-admissible operators (that is, with operators that do not move the Frequency Set), or even with {\it twisted} $h$-admissible operators (that become standard $h$-admissible operators once conjugated with ${\mathcal  W}$ or $Z_L$) it is enough, for computing the  symbol $A$ in this region, to start by performing formal computations on the operators themselves (instead of immediately using the twisted symbolic calculus, that appears to be a little bit too heavy at the beginning).
\vskip 0.3cm
In this chapter, we describe a rather easy way to perform these computations, and we give a simple expression of  the effective Hamiltonian up to ${\mathcal  O}(h^4)$. Moreover, as an example, we also compute its symbol, up to ${\mathcal  O}(h^3)$, in the case $L=1$. Let us inform the reader that the results of this  chapter are not used in the rest of the paper (except for Theorem \ref{remsymbol}), and thus can be skipped without problem at a first reading.
\vskip 0.2cm
We start from the definition of $A$ given in Chapter \ref{decmodop} (in particular  (\ref{defA})):
$$
A= {\mathcal  W}\tilde P{\mathcal  W}^* = Z_L{\mathcal  V}\tilde P{\mathcal  V}^* Z_L^*.
$$
Since $Z_L$ is rather explicit, the problem mainly consists in determining the expansion of ${\mathcal  V}$. Setting,
$$
\Delta := h^{-1}(\Pi_g -\tilde\Pi_0),
$$
and using that $\Pi_g^2 -\Pi_g = \tilde\Pi_0^2 -\tilde\Pi_0=0$, we immediately obtain,
\be
\label{pid+dpi}
\Pi_g\Delta + \Delta\Pi_g = \Delta + h\Delta^2.
\ee
Thus, we deduce from (\ref{defV}),
\begin{eqnarray*}
{\mathcal  V} &=& ((\Pi_g -h\Delta)\Pi_g + (1-\Pi_g +h\Delta)(1-\Pi_g))(1-h^2\Delta^2)^{-\frac12} \\
&=& (1 + h[\Pi_g,\Delta] -h^2\Delta^2)(1-h^2\Delta^2)^{-\frac12}.
\end{eqnarray*}
Then, using the (convergent) series expansion,
$$
(1-h^2\Delta^2)^{-\frac12} = 1 +\sum_{k=1}^\infty \nu_k h^{2k}\Delta^{2k},
$$
with,
$$
\nu_k =\frac12 (\frac12 +1)(\frac12 +2)\dots (\frac12 +k-1) \frac1{k!}=\frac{(2k-1)!}{2^{2k-1}k!(k-1)!},
$$
we obtain,
$$
{\mathcal  V} =1 - ih{\mathcal  V}_1+h^2{\mathcal  V}_2,
$$
where the two selfadjoint operators ${\mathcal  V}_1$ and ${\mathcal  V}_2$ are given by,
\begin{eqnarray*}
{\mathcal  V}_1 &=& i[\Pi_g,\Delta] (1+ \sum_{k=1}^\infty \nu_k h^{2k}\Delta^{2k});\\
{\mathcal  V}_2 &=& -\frac12\Delta^2 + \sum_{k=1}^\infty (\nu_{k+1}-\nu_k)h^{2k}\Delta^{2(k+1)},
\end{eqnarray*}
that is, observing that $\nu_k -\nu_{k+1} = \nu_k/(2k+2)$,
\begin{eqnarray*}
{\mathcal  V}_1 &=& i[\Pi_g,\Delta] F_1(\Delta^2);\\
{\mathcal  V}_2 &=& F_2(\Delta^2),
\end{eqnarray*}
with, (setting also $\nu_0:=1$),
\begin{eqnarray*}
F_1(s) &=&  \sum_{k=0}^\infty \nu_k h^{2k}s^k;\\
F_2 (s) &=& - \sum_{k=0}^\infty \frac{\nu_k}{2(k+1)}h^{2k}s^{k+1}.
\end{eqnarray*}
As a consequence,
$$
{\mathcal  V}^* =1 +ih{\mathcal  V}_1+h^2{\mathcal  V}_2,
$$
and  therefore,
\begin{eqnarray*}
{\mathcal  V}\tilde P{\mathcal  V}^* = \tilde P + ih [\tilde P, {\mathcal  V}_1] + h^2({\mathcal  V}_1\tilde P{\mathcal  V}_1+{\mathcal  V}_2\tilde P+\tilde P{\mathcal  V}_2) + ih^3({\mathcal  V}_2\tilde P{\mathcal  V}_1-{\mathcal  V}_1\tilde P{\mathcal  V}_2)\\
 + h^4{\mathcal  V}_2\tilde P{\mathcal  V}_2,
\end{eqnarray*}
that is,
\begin{eqnarray*}
A=Z_L(\tilde P + ih [\tilde P, {\mathcal  V}_1] + h^2({\mathcal  V}_1\tilde P{\mathcal  V}_1+{\mathcal  V}_2\tilde P+\tilde P{\mathcal  V}_2) + ih^3({\mathcal  V}_2\tilde P{\mathcal  V}_1-{\mathcal  V}_1\tilde P{\mathcal  V}_2)\\
 + h^4{\mathcal  V}_2\tilde P{\mathcal  V}_2)Z_L^*.
\end{eqnarray*}
From now on, we work modulo ${\mathcal  O}(h^5)$ error-terms, and, as we observed at the beginning of this chapter, if we restrict our attention to the relevant region of the phase space,  then formal computations are sufficient and $\Pi_g$ can be replaced by the formal series $\tilde \Pi := \sum_{k\geq 0}h^k\tilde\Pi_k$ constructed in Chapter \ref{qis}. In particular, $\tilde P$ formally commutes with $\tilde \Pi$ and thus, since $[\tilde P, \tilde \Pi_0]=-ihS_0$ (see Chapter \ref{qis}),
\be
\label{com1}
[\tilde P, [\tilde \Pi, \Delta]]= -h^{-1}[\tilde P, [\tilde \Pi, \tilde\Pi_0]] = -h^{-1}[\tilde \Pi , [ \tilde P, \tilde\Pi_0]]= i[\tilde \Pi, S_0],
\ee
where, from now on, $\Delta$ stands for $h^{-1}(\tilde \Pi -\tilde\Pi_0) =\sum_{k\geq 1}h^k\tilde\Pi_k$.
\vskip 0.3cm
 Moreover, from the identities $[\tilde P, \tilde \Pi]=0$, $\tilde \Pi = \tilde \Pi_0 +h\Delta$, we deduce,
 $$
 [ \tilde P, \Delta] = -h^{-1}[ \tilde P, \tilde\Pi_0]= iS_0,
 $$
and therefore,
\begin{eqnarray*}
[\tilde P, {\mathcal  V}_1] &=& [S_0, \tilde \Pi ]F_1(\Delta^2) + i[\tilde\Pi,\Delta][\tilde P, F_1(\Delta^2)];\\
{} [ \tilde P, F_1(\Delta^2)] &=& i\sum_{k=1}^\infty \nu_k h^{2k}\sum_{j=0}^{2k-1}\Delta^jS_0\Delta^{2k-1-j}.
\end{eqnarray*}
Since $\nu_0 =1$ and $\nu_1 = 1/2$, this gives,
\be
[\tilde P, {\mathcal  V}_1] = [S_0, \tilde \Pi ](1+\frac{h^2}2\Delta^2)  -\frac{h^2}2[\tilde\Pi,\Delta](S_0\Delta +\Delta S_0)+ {\mathcal  O}(h^4)
\ee
Moreover,  (\ref{pid+dpi}) implies $\tilde\Pi \Delta\tilde\Pi = h \Delta^2\tilde\Pi= h \tilde\Pi\Delta^2$, and thus, in particular, $\Delta^2$ commutes with $\tilde\Pi$. As  a consequence, we can write,
\begin{eqnarray*}
{\mathcal  V}_1 \tilde P{\mathcal  V}_1&=& F_1( \Delta^2 ) [\tilde\Pi, \Delta]\tilde P[ \Delta, \tilde\Pi] F_1( \Delta^2 )\\
&=& [\tilde\Pi, \Delta]\tilde P[ \Delta, \tilde\Pi]+h^2\re \Delta^2 [\tilde\Pi, \Delta]\tilde P[ \Delta, \tilde\Pi]+{\mathcal  O}(h^4),
\end{eqnarray*}
and, still using (\ref{pid+dpi}), we have,
\begin{eqnarray*}
[\tilde\Pi, \Delta]\tilde P[ \Delta, \tilde\Pi] &=&\tilde\Pi \Delta\tilde P \Delta \tilde\Pi + \Delta\tilde\Pi\tilde P \tilde\Pi\Delta -\tilde\Pi \Delta\tilde P\tilde\Pi \Delta -  \Delta \tilde\Pi\tilde P \Delta\tilde\Pi\\
&=& (\tilde\Pi \Delta + \Delta \tilde\Pi )\tilde P (\Delta \tilde\Pi +\tilde\Pi \Delta)-2\tilde\Pi \Delta\tilde P\tilde\Pi \Delta - 2 \Delta \tilde\Pi\tilde P \Delta\tilde\Pi\\
&=&  (\Delta + h\Delta^2)\tilde P (\Delta + h\Delta^2)-2h\tilde\Pi\Delta^2\tilde P\Delta-2h\Delta\tilde P\Delta^2\tilde\Pi\\
&=& \Delta \tilde P \Delta +h(1-2\tilde\Pi)\Delta^2\tilde P\Delta + h\Delta\tilde P\Delta^2(1-2\tilde\Pi)\\
&=& \frac12(\Delta^2P + P\Delta^2) + \frac{i}2[\Delta, S_0]+2h\re \Delta^2(1-2\tilde\Pi)\tilde P\Delta.
\end{eqnarray*}
Therefore,
\begin{eqnarray*}
{\mathcal  V}_1 \tilde P{\mathcal  V}_1=\re \Delta^2P + \frac{i}2[\Delta, S_0]+2h\re \Delta^2(1-2\tilde\Pi)\tilde P\Delta\\
+h^2\re \Delta^2 (\re \Delta^2P + \frac{i}2[\Delta, S_0])
+{\mathcal  O}(h^3).
\end{eqnarray*}
and, since ${\mathcal  V}_2 = -\frac12\Delta^2-\frac18h^2\Delta^4 + {\mathcal  O}(h^4)$, we obtain,
\begin{eqnarray*}
{\mathcal  V}_1 \tilde P{\mathcal  V}_1+{\mathcal  V}_2\tilde P + \tilde P{\mathcal  V}_2&=&\frac{i}2[\Delta, S_0]+2h\re \Delta^2(1-2\tilde\Pi)\tilde P\Delta\\
&&+h^2\left(\re \Delta^2 (\re \Delta^2P + \frac{i}2[\Delta, S_0])-\frac14\re \Delta^4\tilde P\right)\\
&&+{\mathcal  O}(h^3)\\
&=&\frac{i}2[\Delta, S_0]+2h\re \Delta^2(1-2\tilde\Pi)\tilde P\Delta\\
&&+\frac12h^2\left(\re (i\Delta^2 [\Delta, S_0])+\Delta^2\tilde P\Delta^2+\frac14\re \Delta^4\tilde P\right)\\
&&+{\mathcal  O}(h^3)
\end{eqnarray*}
Finally, since, obviously, $\Delta^2$ also commutes with $\Delta$, thus with $[\tilde\Pi,\Delta]$, too, we see that ${\mathcal  V}_1$ and ${\mathcal  V}_2$ commute together, and therefore,
\begin{eqnarray*}
{\mathcal  V}_2\tilde P {\mathcal  V}_1 - {\mathcal  V}_1\tilde P {\mathcal  V}_2 &=& [\tilde P, {\mathcal  V}_1] {\mathcal  V}_2- [\tilde P, {\mathcal  V}_2] {\mathcal  V}_1\\
&=&-\frac12 [S_0, \tilde\Pi]\Delta^2  +\frac{i}2[\tilde P,\Delta^2][\tilde\Pi, \Delta]+{\mathcal  O}(h^2)\\
&=& -\frac12 [S_0, \tilde\Pi]\Delta^2  -\frac12(S_0\Delta+\Delta S_0)[\tilde\Pi, \Delta]+{\mathcal  O}(h^2).
\end{eqnarray*}
Summing up, we have found,
$$
{\mathcal  V}\tilde P{\mathcal  V}^* = B_0 +hB_1+h^2B_2+h^3B_3+h^4B_4 +{\mathcal  O}(h^5),
$$
with,
\begin{eqnarray*}
B_0 &=& \tilde P\\
B_1 &=& i[S_0, \tilde\Pi]\\
B_2 &=& \frac{i}2[\Delta, S_0]\\
B_3 &=& - \re i[\tilde\Pi,\Delta](S_0\Delta+\Delta S_0)+2\re \Delta^2(1-2\tilde\Pi)\tilde P\Delta\\
B_4&=& \frac12\left(\re (i\Delta^2 [\Delta, S_0])+\Delta^2\tilde P\Delta^2+\frac14\re \Delta^4\tilde P\right)
\end{eqnarray*}
Then, writing $\tilde\Pi = \sum_{k=0}^3h^k\tilde\Pi_k +{\mathcal  O}(h^4)$ and $\Delta =\sum_{k=1}^3h^{k-1}\tilde\Pi_k +{\mathcal  O}(h^3)$, we obtain,
$$
{\mathcal  V}\tilde P{\mathcal  V}^* = C_0 +hC_1+h^2C_2+h^3C_3+h^4C_4 +{\mathcal  O}(h^5),
$$
with,
\begin{eqnarray*}
C_0 &=& \tilde P\\
C_1 &=& i[S_0, \tilde\Pi_0]\\
C_2 &=& \frac{i}2[S_0,\tilde \Pi_1]\\
C_3 &=&\frac{i}2[S_0, \tilde\Pi_2] - \re i[\tilde\Pi_0,\tilde\Pi_1](S_0\tilde\Pi_1+\tilde\Pi_1 S_0)+2\re \tilde\Pi_1^2(1-2\tilde\Pi_0)\tilde P\tilde\Pi_1\\
C_4&=& \frac{i}2[S_0, \tilde\Pi_3]  - \re i[\tilde\Pi_0,\tilde\Pi_2](S_0\tilde\Pi_1+\tilde\Pi_1 S_0)-\re  i[\tilde\Pi_0,\tilde\Pi_1](S_0\tilde\Pi_2+\tilde\Pi_2 S_0)\\ 
&&+2\re (\tilde\Pi_1\tilde\Pi_2+\tilde\Pi_2\tilde\Pi_1)(1-2\tilde\Pi_0)\tilde P\tilde\Pi_1-4\re \tilde\Pi_1^3\tilde P\tilde\Pi_1\\
&&+2\re 
\tilde\Pi_1^2(1-2\tilde\Pi_0)\tilde P\tilde\Pi_2+\frac12\left(\re (i\tilde\Pi_1^2 [\tilde\Pi_1, S_0])+\tilde\Pi_1^2\tilde P\tilde\Pi_1^2+\frac14\re \tilde\Pi_1^4\tilde P\right)
\end{eqnarray*}
Now, due to (\ref{decS})-(\ref{pi1}), we observe that,
 $$\tilde\Pi_0S_0\tilde\Pi_0 =\tilde\Pi_0^\perp S_0\tilde\Pi_0^\perp =\tilde\Pi_0\tilde\Pi_1\tilde\Pi_0=\tilde\Pi_0^\perp \tilde\Pi_1\tilde\Pi_0^\perp =0.
 $$
  As a consequence,
$$
\tilde \Pi_0C_1\tilde\Pi_0 = i\tilde\Pi_0[S_0, \tilde\Pi_0]\tilde\Pi_0=0,
$$
and,
$$
 \tilde \Pi_0[\tilde\Pi_0,\tilde\Pi_1](S_0\tilde\Pi_1+\tilde\Pi_1 S_0)\tilde \Pi_0=\tilde \Pi_0[\tilde\Pi_0,\tilde\Pi_1]\tilde \Pi_0(S_0\tilde\Pi_1+\tilde\Pi_1 S_0)\tilde \Pi_0=0;
 $$
\begin{eqnarray*}
 \tilde \Pi_0\tilde\Pi_1^2(1-2\tilde\Pi_0)\tilde P\tilde\Pi_1\tilde \Pi_0&=&
\tilde \Pi_0\tilde\Pi_1^2\tilde\Pi_0^\perp\tilde P\tilde\Pi_1\tilde \Pi_0+ \tilde \Pi_0\tilde\Pi_1^2(1-2\tilde\Pi_0)[\tilde P, \tilde\Pi_0^\perp]\tilde\Pi_1\tilde \Pi_0\\
&=& ih\tilde \Pi_0\tilde\Pi_1^2(1-2\tilde\Pi_0)S_0\tilde\Pi_1\tilde \Pi_0\\
&=& -ih\tilde \Pi_0\tilde\Pi_1^2S_0\tilde\Pi_1\tilde \Pi_0.
\end{eqnarray*}
(In the last two steps we have used that $\tilde\Pi_0\tilde\Pi_1^2\tilde\Pi_0^\perp =\tilde\Pi_0^\perp S_0\tilde\Pi_1\tilde\Pi_0=0$.)
Since we also have $Z_L = Z_L\tilde\Pi_0$ and $Z_L^* =\tilde\Pi_0Z_L^*$, we deduce,
\begin{eqnarray}
&& Z_L C_1Z_L^* =0;\nonumber\\
\label{c3h}
&& Z_LC_3Z_L^* = \frac{i}2Z_L[S_0, \tilde\Pi_2]Z_L^*+2h\im\tilde \Pi_0\tilde\Pi_1^2S_0\tilde\Pi_1\tilde \Pi_0.
\end{eqnarray}
In particular, since $A=Z_L{\mathcal  V}\tilde P{\mathcal  V}^*Z_L^*$, we have proved,
\begin{proposition} The effective Hamiltonian $A$ verifies, 
\be
\label{exphameff}
A=A_0+h^2A_2+h^3A_3+{\mathcal  O}(h^4),
\ee
with,
\begin{eqnarray*}
A_0 &=& Z_L\tilde PZ_L^*\\
A_2 &=& \frac{i}2 Z_L[S_0,\tilde \Pi_1] Z_L^*\\
A_3 &=&\frac{i}2Z_L[S_0, \tilde\Pi_2]Z_L^*.
\end{eqnarray*}
\end{proposition}
It is interesting to observe that, at this level, the absence of a term in $h$ (that is, an extra-term of the form $hA_1$)  is completely general and, in particular, is not related to any particular form of $\boldsymbol{\omega}$ (however, some term in $h$ may be hidden in $A_0$, as we shall see in the sequels).
\vskip 0.3cm
Here, we have stopped the computation of $A$ at the third power of $h$, but it is clear from the expression of $C_4$ and (\ref{c3h}) that the coefficient of $h^4$ can be written down, too (but has a more complicated form). Of course, pushing forward the series and spending more time in the calculation would permit to also obtain the next terms. 
\vskip 0.3cm
From that point, in order to have an even more explicit expression of $A$ (in particular to compute its symbol), one must use the expressions of $\tilde\Pi_1$ and $\tilde\Pi_2$ obtained in Chapter \ref{qis}. Let us do it in the case $L=1$. In that case, setting $\lambda (x):= \lambda_{L'+1}(x)$, one has $\tilde \Pi_0 (z-\tilde Q(x))^{-1} =(z-\lambda (x))^{-1} \tilde \Pi_0 $, and thus,
$$
\tilde\Pi_0^\perp\tilde\Pi_1\tilde\Pi_0=-\frac1{2\pi}\oint_{\gamma (x)}\frac{(z-\tilde Q(x))^{-1}\tilde\Pi_0^\perp(x) S_0\tilde\Pi_0(x)}{ z-\lambda  (x)}dz =-iR'(\lambda (x))S_0,
$$
where $R'(x,z):= \tilde\Pi_0^\perp(x)(z-\tilde Q(x))^{-1}\tilde\Pi_0^\perp(x)$ is the so-called {\it reduced resolvent} of $\tilde Q(x)$.
\vskip 0.3cm
As a consequence,
$$
\tilde\Pi_0[S_0,\tilde \Pi_1]\tilde\Pi_0=S_0\tilde\Pi_0^\perp\tilde\Pi_1\tilde\Pi_0-\tilde\Pi_0\tilde\Pi_1\tilde\Pi_0^\perp S_0 =-2i S_0R'(x, \lambda (x))S_0,
$$
that leads to,
$$
A_2 = Z_1S_0R'(x,\lambda (x))S_0Z_1^*.
$$
In the same way,
$$
\tilde\Pi_0^\perp\tilde\Pi_2\tilde\Pi_0=-iR'(x,\lambda (x))S_1\tilde\Pi_0,
$$
and therefore,
$$
A_3 = \re Z_1S_0R'(x,\lambda (x))S_1Z_1^*.
$$
Now, we can start to use the twisted symbolic calculus introduced in Chapter \ref{twao}. We denote by $s_0 =(s_0^j)_{0\leq j\leq r}$ and $\pi_0 =(\pi_0^j)_{0\leq j\leq r}$ the (twisted) symbols of $S_0$ and $\tilde\Pi_0$ respectively. We also set $\tilde \omega =(\tilde\omega_j)_{0\leq j\leq r}$, where,
$$
\tilde\omega_j(x,\xi ) := \omega (x,\xi) +h\sum_{|\beta|\leq m-1}\omega_{\beta ,j}(x)\xi^\beta, \quad ((x.\xi)\in T^*\Omega_j),
$$
is the symbol of the operator introduced in (\ref{conjomega}) (we remind that we work with the standard quantization of symbols, as described in Chapter \ref{stpdo}). From (\ref{comppi0})-(\ref{step1}) and the considerations of Chapter \ref{twao} (and since $\pi_0^j=\pi_0^j(x)$ does not depend on $\xi$), it s easy to see that,
$$
s_0^j =\partial_\xi\tilde\omega_j \partial_x\pi_0^j +i\sum_{|\beta|\leq m-1}[\omega_{\beta ,j}(x), \pi_0^j(x)]\xi^\beta -\frac{ih}2\sum_{|\alpha|=2}(\partial_\xi^\alpha \omega)(\partial_x^\alpha\pi_0^j )+{\mathcal  O}(h^2)
$$
If we also set, 
$$
\tilde Q_j(x):= U_j(x)\tilde Q(x)U_j(x)^{-1},
$$
then,  the symbol $\rho=(\rho_j)_{0\leq j\leq r}$ of
$R'(x,\lambda (x))$ is simply given by,
$$
\rho_j(x)= (1-\pi_0^j(x))(\lambda (x) - \tilde Q_j(x))^{-1}(1-\pi_0^j(x)),
$$
and thus, the symbol $\sigma_2 =(\sigma_2^j)_{0\leq j\leq r}$ of $S_0R'(x,\lambda (x))S_0$ verifies,
$$
\sigma_2^j(x,\xi )= s_0^j(x,\xi)\rho_j(x)s_0^j(x,\xi) + \frac{h}{i}\partial_\xi s_0^j(x,\xi)\partial_x(\rho_j(x)s_0^j(x,\xi))+{\mathcal  O}(h^2).
$$
From (\ref{pi1})-(\ref{corr1}), we also obtain,
\begin{eqnarray*}
&&\tilde \Pi_1= i[S_0, R'(x,\lambda (x))]\\
&& S_1 = \frac{i}{h}[\boldsymbol{\omega} +\zeta W,\tilde \Pi_1].
\end{eqnarray*}
Therefore, since $\omega$ and $\zeta W$ are scalar operators,  the respective symbols $\pi_1 =(\pi_1^j)_{0\leq j\leq r}$ and $s_1=(s_1^j)_{0\leq j\leq r}$ of $\tilde \Pi_1$ and $S_1$,   verify,
\begin{eqnarray*}
&& \pi_1^j(x,\xi)= i[s_0^j(x,\xi), \rho_j (x)]+{\mathcal  O}(h)=  i\partial_\xi\omega(x,\xi)[\partial_x\pi_0^j(x), \rho_j (x)]+{\mathcal  O}(h)\\
&& s_1^j = \{\omega+\zeta W, \pi_1^j\}+{\mathcal  O}(h)= \partial_\xi\omega\cdot \partial_x\pi_1^j-\partial_\xi\pi_1^j\cdot \partial_x(\omega + \zeta W)+{\mathcal  O}(h),
\end{eqnarray*}
and thus,
\begin{eqnarray}
s_1^j =i\sum_{k,\ell=1}^n \left( (\partial_{\xi_k}\omega)\partial_{x_k}(\partial_{\xi_\ell}\omega[\partial_{x_\ell}\pi_0^j, \rho_j ])-(\partial_{\xi_k}\partial_{\xi_\ell}\omega )[\partial_{x_\ell}\pi_0^j, \rho_j ]\partial_{x_k}(\omega +\zeta W)\right)\nonumber\\
\label{symbs1}
 +{\mathcal  O}(h).\hskip 1cm{}
\end{eqnarray}
This permits to compute the symbol $\sigma_3 =(\sigma_3^j)_{0\leq j\leq r}$ of $\re S_0R'(x,\lambda (x))S_1$, by using  the formula,
\be
\label{sigma3}
\sigma_3^j(x,\xi )= \frac12\partial_\xi\omega \cdot \left((\partial_x\pi_0^j)\rho_js_1^j + s_1^j\rho_j(\partial_x\pi_0^j)\right)+{\mathcal  O}(h).
\ee
Observe that one also has,
$$
\partial_x\pi_0^j(x) =\la \cdot, \nabla_xu_j(x)\ra_{\mathcal  H} u_j(x) + \la \cdot, u_j(x)\ra_{\mathcal  H} \nabla_x u_j(x),
$$
where $\la \cdot, u\ra_{\mathcal  H}$ stands for the operator $w \mapsto \la w, u\ra_{\mathcal  H}$, and $u_j=:U_j(x)u_{L'+1}(x)$ is the normalized eigenfunction of $\tilde Q_j(x)$ associated with $\lambda (x)$.
\vskip 0.3cm
Finally, we use the following elementary remark: let $B$ is a twisted $h$-admissible (or PDO) operator on $L^2(\R^n; {\mathcal  H})$, with symbol $b=(b_j)_{0\leq j\leq r}$, and let $u(x), v(x)\in{\mathcal  H}$  such that, for all $j=0, \dots,r$, $u_j(x):=U_j(x)u(x)$ and $v_j(x):=U_j(x)v(x)$ are in $C^\infty (\Omega_j ;{\mathcal  H})$. Denote by ${\mathcal  Z}_u$, ${\mathcal  Z}_v$ the operators  $L^2(\R^n; {\mathcal  H})\to L^2(\R^n)$ defined by ,
$$
{\mathcal  Z}_uw:= \la w, u\ra_{\mathcal  H}\,\,\, ;\,\,\, {\mathcal  Z}_vw:= \la w, v\ra_{\mathcal  H}.
$$
Then, the symbol $\check b$ of the (standard) $h$-admissible operator ${\mathcal  Z}_vB{\mathcal  Z}_u^*$ verifies,
$$
\forall\, (x,\xi)\in T^*\Omega_j,\,\,\, \check b (x,\xi ) =\la  b_j(x,\xi )\sharp u_j(x), v_j(x)\ra_{\mathcal  H},
$$
where the operation $\sharp$ is defined in an obvious way, by substituting the usual product with the action of an operator (here, the various derivatives of  $b_j(x,\xi )$) on a function (here, the various derivatives of  $u_j(x )$).
\vskip 0.3cm
We can clearly apply this remark to compute the symbol of $A_2$ and $A_3$, but also that of $A_0$, since we have,
$$
A_0 = Z_1\tilde PZ_1^* = {\mathcal  Z}_u\tilde P{\mathcal  Z}_u^*= {\mathcal  Z}_{Q_0u}Q_0^{-1}\tilde P{\mathcal  Z}_u^*,
$$
with $u:=\tilde u_{L'+1}$ (defined in Chapter \ref{modop}), and, by Proposition \ref{omegaPDO}, we know that $Q_0^{-1}\tilde P$ is a twisted PDO.
\vskip 0.3cm
Combining all the previous computations, using that $\tilde Q_j(x) u_j(x) = \lambda (x)u_j(x)$ for all $j=0,\dots ,r$ and $x\in\Omega_j$, and gathering (as far as possible) the terms with same homogeneity in $h$, we finally arrive to the following result (leaving some details to the reader):
\begin{proposition}
\label{symhameff}
In the case ${\rm Rank }\hskip 1pt \Pi_0(x)=1$, the effective Hamiltonian $A$ verifies (\ref{exphameff}) with,
\begin{eqnarray}
A_0 &=& Z_1\tilde PZ_1^*;\nonumber\\
\label{hameffl1}
A_2&=& \frac1{h^2}Z_1[\tilde P, \tilde\Pi_0]R'(x,\lambda (x))[ \tilde\Pi_0, \tilde P]Z_1^*;\\
A_3&=&  \frac1{h^3}\re Z_1[\tilde P, \tilde\Pi_0]R'(x,\lambda (x))[ [ [\tilde P, \tilde\Pi_0],R'(x,\lambda (x)) ] , {\boldsymbol{\omega}}+\zeta W] Z_1^*,\nonumber
\end{eqnarray}
where $\lambda (x)$ is the (only) eigenvalue of $\tilde Q(x)\tilde \Pi_0$, and $R'(x,\lambda (x))=\tilde\Pi_0^\perp(x)(\lambda (x)-\tilde Q(x))^{-1}\tilde\Pi_0^\perp(x)$ is the reduced resolvent of $\tilde Q(x)$.
\vskip 0.3cm
Moreover, the symbol $a(x,\xi ;h)$ of $A$ verifies,
$$
a(x,\xi ;h)= a_0(x,\xi ) +ha_1(x,\xi) +h^2a_2(x,\xi)  +{\mathcal  O}(h^3),
$$
with, for any $(x,\xi)\in T^*\Omega_j$ ($j=0, \dots ,r$ arbitrary),
\begin{eqnarray*}
 a_0(x,\xi ) &=& \omega  (x,\xi ;h) + \lambda (x) + \zeta (x)W(x) ;\\
 a_1(x,\xi ) &=&  \sum_{|\beta|\leq m-1}\la \omega_{\beta ,j}(x)u_j(x), u_j(x)\ra\xi^\beta -i\la\nabla_{\xi}\omega(x,\xi)\nabla_x u_j(x),u_j(x)\ra ;\\
 a_2(x,\xi ) &=& \sum_{k,\ell =1}^n(\partial_{\xi_k}\omega )(\partial_{\xi_\ell }\omega)\la \rho_j(x)\partial_{x_k}u_j, \partial_{x_\ell}u_j\ra -\frac12\sum_{|\alpha| =2}(\partial_\xi^\alpha\omega)\la\partial_x^\alpha u_j , u_j\ra\\
&& -i\sum_{|\beta|\leq m-1}\la \omega_{\beta ,j}(x)\nabla_xu_j(x), u_j(x)\ra \cdot \nabla_\xi(\xi^\beta)\\
&& -2\im \sum_{ |\beta|\leq m-1}\nabla_\xi\omega(x,\xi)\la \omega_{\beta ,j}(x)\rho_j(x)\nabla_xu_j(x), u_j(x)\ra \xi^\beta \\
&& +\sum_{|\beta|, |\gamma| \leq m-1}
\la \omega_{\beta ,j}(x)\rho_j(x)u_j(x), 
\omega_{\beta ,j}(x)u_j(x)\ra  \xi^{\beta +\gamma}.
\end{eqnarray*}
\end{proposition}
\begin{remark}\sl Although some of these terms may seem to depend on the choice of $j$ verifying $(x,\xi)\in T^*\Omega_j$, actually we know that this cannot be the case. In fact, the independency with respect to $j$ is due to the compatibility conditions (\ref{compatib}) satisfied by the symbols of twisted pseudodifferential operators.
\end{remark}
\begin{remark}\sl Actually, it results from the previous computations that (\ref{hameffl1}) is still valid in the (slightly) more general case where $L$ is arbitrary and $\lambda_{L'+1}(x) = \dots =\lambda_{L'+L}(x)$ for all $x\in \Omega$.
\end{remark}
\begin{remark}\sl Using (\ref{symbs1})-(\ref{sigma3}), one can find an expression for  the $h^3$-term of the symbol of $A$, too. We leave it as an exercise to the reader.
\end{remark}
\chapter{Propagation of Wave-Packets}
\label{wp}
\setcounter{equation}{0}%
\setcounter{theorem}{0}%
In this chapter, we assume $L=1$ and we make the following additional assumption on the coefficients $c_\alpha$ of $\boldsymbol{\omega}$:
\be
\label{coefclas}
c_\alpha (x;h) \sim \sum_{k=0}^\infty h^kc_{\alpha ,k}(x),
\ee
with $c_{\alpha ,k}$ independent of $h$. Then, in a similar spirit as in \cite{Ha6},
we investigate the evolution of an initial state of the form,
\be
\label{cohststart}
\varphi_0(x) = (\pi h)^{-n/4}f(P)\Pi_g(e^{ix\xi_0/h - (x-x_0)^2/2h} u_{L'+1}(x)),
\ee
where $(x_0,\xi_0)\in T^*\Omega$ is fixed,  $f, g\in C_0^\infty (\R)$ are such that $f=1$ near $a_0(x_0,\xi_0)$  (here, $a_0(x,\xi )$ is the same as in Corollary \ref{CORT}), $g=1$ near $\supp f$, and $\Pi_g$ is constructed as in Chapter \ref{qis}, starting from the operator $\tilde P$ constructed in Chapter \ref{modop} with $K\ni x_0$.  In particular, since $e^{-(x-x_0)^2/2h}$ is exponentially small for $x$ outside any neighborhood of $x_0$, by Lemma \ref{compP}, we have,
$$
\varphi_0(x) = (\pi h)^{-n/4}f(\tilde P)\Pi_g(e^{ix\xi_0/h - (x-x_0)^2/2h} \tilde u_{L'+1}(x)) +{\mathcal  O}(h^\infty),
$$
in $L^2(\R^n;{\mathcal  H})$.
Moreover,
due to the properties of $\Pi_g$, and the fact that  the coherent state $\phi_0:= (\pi h)^{-n/4}e^{ix\xi_0/h - (x-x_0)^2/2h}$ is normalized in ${L^2(\R^n)}$, we also obtain,
$$
\varphi_0(x) =(\pi h)^{-n/4}f(\tilde P)e^{ix\xi_0/h - (x-x_0)^2/2h}\tilde u_{L'+1}(x) +{\mathcal  O}(h),
$$ 
and thus, in particular, $\Vert \varphi_0\Vert =1+{\mathcal  O}(h)$. Actually, we even have the following better result:
\begin{proposition}\sl  The function $\varphi_0$ admits, in $L^2(\R^n;{\mathcal  H})$, an asymptotic expansion of the form,
\begin{eqnarray}
\label{asexpphi1}
\varphi_0(x) \sim (\pi h)^{-n/4}e^{ix\xi_0/h - (x-x_0)^2/2h}\sum_{k=0}^\infty h^kv_k(x) +{\mathcal  O}(h^\infty),
\end{eqnarray}
with 
$v_k \in  L^\infty (\R^n; {\mathcal  H})$ ($k\geq 0$), and $v_0(x) =\tilde u_{L'+1}(x)+{\mathcal  O}(|x-x_0|)$ in ${\mathcal  H}$, uniformly with respect to $x\in \R^n$. Moreover, for any $j\in\{0, 1,\dots,r\}$ and any $\i_j\in C_d^\infty (\Omega_j)$, 
the function $U_j \i_j\varphi_0$ admits, in $C_d^\infty(\Omega_j;{\mathcal  H})$, an asymptotic expansion of the form,
\be
\label{asexpphi2}
U_j (x)\i_j(x)\varphi_0(x) \sim (\pi h)^{-n/4}e^{ix\xi_0/h - (x-x_0)^2/2h} \sum_{k=0}^\infty h^k\i_j(x)v_{j,k}(x)  +{\mathcal  O}(h^\infty), 
\ee
with $v_{j,k}\in C^\infty(\Omega_j;{\mathcal  H})$, $v_{j,0} (x) =U_j(x) \tilde u_{L'+1}(x)+{\mathcal  O}(|x-x_0|)$.
\end{proposition}
{\em Proof -- } For $j=0,1,\dots,r$, let $\i_j\in C_d^\infty (\Omega_j)$, such that $\sum\i_j =1$, and let $\tilde \i_j\in C_d^\infty (\Omega_j)$, such that $\tilde \i_j=1$ near $\supp \i_j$. Then,  since $f(\tilde P)$ and $\Pi_g$ are twisted $h$-admissible operators,  have,
\begin{eqnarray*}
\varphi_0 &=& \sum_j \i_j\varphi_0\\
&=& \sum_j U_j^{-1}\tilde \i_j U_j\i_j f(\tilde P)\tilde\i_j^2\Pi_g\tilde\i_j^2(\phi_0(x)\tilde u_{L'+1}(x))+{\mathcal  O}(h^\infty)\\
&=& \sum_j U_j^{-1}\tilde \i_j U_j\i_j f(\tilde P)U_j^{-1}\tilde\i_jU_j\tilde\i_j\Pi_g\tilde\i_j^2(\phi_0(x)\tilde u_{L'+1}(x))+{\mathcal  O}(h^\infty),
\end{eqnarray*}
and thus, by Lemma \ref{tec2}, and setting $\tilde P_j := U_j\tilde \i_j\tilde P U_j^{-1}\tilde\i_j$, $\Pi_{g,j}:=U_j\tilde \i_j \Pi_g U_j^{-1}\tilde\i_j$, and $u_{L'+1,j}(x):=U_j(x)\tilde\i_j(x) \tilde u_{L'+1}(x)$ ($\in C_d^\infty (\Omega_j ;{\mathcal  H})$), we obtain,
\be
\label{expphi0}
\varphi_0 = \sum_{j=0}^rU_j^{-1}\i_jf(\tilde P_j)\Pi_{g,j} (\phi_0(x)u_{L'+1,j}(x))+{\mathcal  O}(h^\infty).
\ee
Now, using the results of Chapters 4 and 6, we see that $f(\tilde P_j)\Pi_{g,j}$ is an $h$-admissible operator on $L^2(\R^n;{\mathcal  H})$, with symbol $b_j$ verifying,
\begin{eqnarray*}
&& b_j(x,\xi;h ) \sim \sum_{k=0}^\infty h^kb_{j,k}(x,\xi );\\
&& b_{j,0}(x,\xi ) = f(\tilde\i_j(x)^2(\omega_0 (x,\xi )+\tilde Q_j(x) +W(x)))\tilde\i_j(x)^2\tilde\Pi_{0,j}(x),
\end{eqnarray*}
where $\omega_0(x,\xi ):= \sum_{|\alpha|\leq m}c_{\alpha ,0}(x)\xi^\alpha$, $\tilde Q_j(x)=U_j(x)\tilde Q(x)U_j(x)^{-1}$, and $\tilde\Pi_{0,j}(x)=U_j(x) \tilde\Pi_0(x)U_j(x)^{-1}$. Moreover, we have,
$$
{\rm Op}_h(b_j)(\phi_0u_{L'+1,j})(x;h)  =\frac1{(2\pi h)^n}\int e^{i(x-y)\xi /h+iy\xi_0 /h}\rho (x,y,\xi ;h)dyd\xi,
$$
with,
$$
\rho (x,y,\xi ;h)=(\pi h)^{-n/4}e^{ -(y-x_0)^2/2h}b_j(x,\xi ;h)u_{L'+1,j}(y),
$$
and it is easy to check that, for any $\alpha, \beta\in \Z_+^n$, one has,
$$
\Vert (hD_y)^\alpha(hD_\xi)^\beta \rho (x,y,\xi ;h)\Vert_{\mathcal  H} ={\mathcal  O}(h^{|\alpha|/2 +|\beta| }),
$$
uniformly for $(x,y,\xi )\in\R^{3n}$and $h>0$ small enough. As a consequence, we can perform a standard stationary phase expansion in the previous (oscillatory) integral (see, e.g., \cite{DiSj1, Ma2}), and since the unique critical point is given by $y=x$ and $\xi =\xi_0$, we obtain,
$$
{\rm Op}_h(b_j)(\phi_0v_j)(x;h)= e^{ix\xi_0 /h}w_j(x;h) +{\mathcal  O}(h^\infty ),
$$
with,
$$
w_j(x;h) \sim \sum_{k=0}^\infty \frac{h^k}{i^kk!}(\nabla_y\cdot\nabla_\xi)^k\rho (x,y,\xi ;h)\left\vert_{\genfrac{}{}{0pt}{}{y=x}{\xi =\xi_0}} \right..
$$
Therefore, since $e^{ (y-x_0)^2/2h}\nabla_y e^{ -(y-x_0)^2/2h} = \nabla_y - \frac{y-x_0}{h}$, and, for any $k\in\N$,  $|y-x_0|^ke^{-(y-x_0)^2/2h}={\mathcal  O}(h^{k/2})$, we also obtain, 
$$
{\rm Op}_h(b_j)(\phi_0u_{L'+1,j})(x;h) = (\pi h)^{-n/4}e^{ ix\xi_0/h-(x-x_0)^2/2h}\tilde w_j(x;h),
$$
with,
\be
\label{wtilde}
\tilde w_j(x,h)=
\sum_{k=0}^N \frac{h^k}{i^kk!}((\nabla_y- h^{-1}(y-x_0))\cdot\nabla_\xi)^k b_j(x,\xi ;h)u_{L'+1,j}(y)\left\vert_{\genfrac{}{}{0pt}{}{y=x}{\xi =\xi_0}} \right. +{\mathcal  O}(h^{N/2}),
\ee
for any $N\geq 0$. Then,  taking a ressummation of  the formal series in $(x-x_0)$ obtained for each degree of homogeneity in $h$ in  (\ref{wtilde}), we obtain an asymptotic expansion of $\tilde w_j$, of the form,
$$
\tilde w_j(x,h)\sim \sum_{k=0}^\infty h^k\tilde w_{j,k}(x).
$$
(Alternatively -- and equivalently -- one could have used instead the stationary phase theorem with complex-valued phase function \cite{MeSj1} Theorem 2.3, with the phase $(x-y)\xi +y\xi_0 +i(y-x_0)^2/2$.) In particular, the first coefficient $\tilde w_{j,0}(x)$ is obtained as a resummation of the formal series $\sum_{k\geq 0}\frac{i^k}{k!}((y-x_0))\cdot\nabla_\xi)^kb_j(x,\xi ;h)u_{L'+1,j}(y)\left\vert_{\genfrac{}{}{0pt}{}{y=x}{\xi =\xi_0}} \right.$, and thus,
\begin{eqnarray*}
\tilde w_{j,0}(x)(x) &=& b_j(x,\xi_0 ;h)u_{L'+1,j}(x) +{\mathcal  O}(|x-x_0|)\\
&=& f(\tilde\i_j(x)^2(\omega_0 (x,\xi_0 )+\tilde Q_j(x) +W(x)))\tilde\i_j(x)^2u_{L'+1,j}(x)\\
&&\hskip 6cm +{\mathcal  O}(|x-x_0|)\\
&=& f(\tilde\i_j(x)^2(\omega_0 (x,\xi_0 )+\lambda_{L'+1} (x) +W(x)))\tilde\i_j(x)^2u_{L'+1,j}(x)\\
&&\hskip 6cm +{\mathcal  O}(|x-x_0|)\\
&=& f(\tilde\i_j(x)^2(a_0 (x_0,\xi_0 ))\tilde\i_j(x)^2u_{L'+1,j}(x) +{\mathcal  O}(|x-x_0|).
\end{eqnarray*}
Going back to (\ref{expphi0}), this gives an asymptotic expansion for $\varphi_0$ of the form (\ref{asexpphi1}), with,
\begin{eqnarray*}
v_0(x)&=&\sum_{j=0}^rU_j(x)^{-1}\i_j(x)f(\tilde\i_j(x)^2a_0 (x_0,\xi_0 ))u_{L'+1,j}(x) +{\mathcal  O}(|x-x_0|)\\
&=&\sum_{j=0}^rU_j(x)^{-1}\i_j(x)f(a_0 (x_0,\xi_0 ))u_{L'+1,j}(x)+{\mathcal  O}(|x-x_0|)\\
&=& \sum_{j=0}^rU_j(x)^{-1}\i_j(x)u_{L'+1,j}(x)+{\mathcal  O}(|x-x_0|)\\
&=& \tilde u_{L'+1}(x) + {\mathcal  O}(|x-x_0|).
\end{eqnarray*}
The asymptotic expansion (\ref{asexpphi2}) is obtained exactly in  the same way.\hfill$\bullet$
\vskip 0.2cm
As a consequence, we also obtain,
\begin{proposition}\sl 
\label{FSCS}
 For any $j\in\{0,1,\dots,r\}$, one has,
$$
FS(U_j\varphi_0) = \{ (x_0,\xi_0 )\}\cap T^*\Omega_j.
$$
\end{proposition}
{\em Proof -- } For $\i_j\in C_d^\infty (\Omega_j)$ fixed, we denote by  $w_j(x;h)$ a resummation of the formal series $\sum_{k\geq 0}h^kU_j(x)\i_j(x) v_{j,k}(x)$ in $C_d^{\infty}(\Omega_j ;{\mathcal  H})$, where the $v_{j,k}$'s are those in (\ref{asexpphi2}).
Then, defining,
\begin{eqnarray*}
&&A = A(x,hD_x) := (hD_x -\xi_0)^2 + (x-x_0)^2\\
 &&\hskip 2.5cm = (hD_x -\xi_0+i(x-x_0))\cdot (hD_x -\xi_0-i(x-x_0)) +nh,
\end{eqnarray*}
a straightforward computation gives,
$$
A(U_j\varphi_0) = A(\phi_0(x)w_j(x;h)) +{\mathcal  O}(h^\infty)= h\phi_0(x)Bw_j(x;h)+{\mathcal  O}(h^\infty)
$$
with $Bw_j(x;h) :=2i(x-x_0)\cdot \partial_xw_j - ih\partial_x^2w_j + nw_j $, and thus, by an iteration,
$$
A^N(\phi_0(x)w_j(x;h)) = h^N\phi_0(x)B^Nw_j+{\mathcal  O}(h^\infty),
$$
for any $N\geq 1$. In particular, due to the form of $B$, and since $\Vert (x-x_0)^\alpha\phi_0\Vert ={\mathcal  O}(1)$ for any $\alpha\in\Z_+^n$ (actually, ${\mathcal  O}(h^{|\alpha |/2})$), we obtain, 
$$
\Vert A^N(U_j\varphi_0)\Vert_{L^2(\Omega'_j,{\mathcal  H})} ={\mathcal  O}(h^N),
$$
for any $\Omega'_j\subset\subset\Omega_j$. Now, if $(x_1,\xi_1)\in T^*\Omega_j$ is different from $(x_0,\xi_0)$, then $A^N$ is elliptic at $(x_1,\xi_1)$ and thus, given any $\i\in C_0^\infty (T^*\Omega_j) $ with $\i(x_1,\xi_1)=1$, the standard construction of a microlocal parametrix (see, e.g., \cite{DiSj1}) gives an uniformly bounded  operator $A'_N$ , such that,
$$
A'_N\circ A^N =\i (x,hD_x)+{\mathcal  O}(h^\infty ).
$$
As a consequence, we obtain,
$$
\Vert\i (x,hD_x)(U_j\varphi_0)\Vert_{L^2(\Omega'_j,{\mathcal  H})} ={\mathcal  O}(h^N),
$$
for all $N\geq 1$. Therefore $(x_1,\xi_1)\notin FS (\phi_0(x)v_j(x))$, and thus, we have proved,
$$
FS(U_j\varphi_0) \subset \{ (x_0,\xi_0 )\}\cap T^*\Omega_j.
$$
This means that $FS(U_j\varphi_0)$ consists in at most one point. Conversely, if $x_0\in\Omega_j$ and $FS(U_j\varphi_0)=\emptyset$, by the ellipticity of $A^N$ as $\vert\xi\vert\to \infty$, we would have (see, e.g., \cite{Ma2} Prop. 2.9.7),
 $$
 \Vert U_j\varphi_0\Vert_{\Omega_j'} ={\mathcal  O}(h^\infty ),
 $$
 for any $\Omega_j'\subset\subset \Omega_j$. But this contradicts the fact that $ \Vert U_j\varphi_0\Vert_{\Omega_j'} =  \Vert \varphi_0\Vert_{\Omega_j'} =1+{\mathcal  O}(h)$ if $x_0\in \Omega_j'$.
\hfill$\bullet$
\vskip 0.3cm
Now, applying Theorem \ref{MAINTH} and Corollary \ref{CORT} (or rather Remark \ref{esttemps2}), we obtain,
\be
\label{mainrepcs}
e^{itP/h}\varphi_0 ={\mathcal  W}^*e^{-itA/h}{\mathcal  W}\varphi_0 +{\mathcal  O}(\la t\ra h^\infty ),
\ee
uniformly for $t\in [0, T_{\Omega'}(x_0,\xi_0))$,where $\Omega'\subset\subset\Omega$ is the same as the one used to define $\tilde P$ in Chapter \ref{modop}, and 
\be
\label{defthamilt}
T_{\Omega'}(x_0,\xi_0):=\sup \{ T>0\, ;\, \pi_x(\cup_{t\in [0,T]} \exp tH_{a_0}(x_0,\xi_0))\subset \Omega'\}.
\ee
Moreover, by Lemma \ref{FSW} and Proposition \ref{FSCS}, we see that, 
\be
FS({\mathcal  W}\varphi_0) =\{(x_0,\xi_0)\}.
\ee
Assuming, e.g., that $x_0\in \Omega_1$, and taking $\i_1\in C_0^\infty (\Omega_1)$ such that $\i_1 =1$ in a neighborhood of $x_0$, we also have,
$$
{\mathcal  W}\varphi_0 = {\mathcal  W}\i_1^2\varphi_0 +{\mathcal  O}(h^\infty ) = {\mathcal  W}U_1^{-1}\i_1 U_1\i_1\varphi_0+{\mathcal  O}(h^\infty ),
$$
and therefore, using (\ref{asexpphi2}),  (\ref{defW}), and the fact that $ {\mathcal  W}U_1^{-1}\i_1$ is an $h$-admissible operator from $L^2(\R^n;{\mathcal  H})$ to $L^2(\R^n)$ (see Theorem \ref{th2}), we obtain as before (by a stationary phase expansion),
\be
{\mathcal  W}\varphi_0(x;h)\sim (\pi h)^{-n/4}e^{ix\xi_0/h - (x-x_0)^2/2h} \sum_{k=0}^\infty h^k w_k(x)  +{\mathcal  O}(h^\infty), 
\ee
with $w_k\in C_b^\infty(\R^n)$, $w_0 (x) =\la \tilde u_{L'+1}(x),\tilde u_{L'+1}(x)\ra +{\mathcal  O}(|x-x_0|)= 1+{\mathcal  O}(|x-x_0|)$, and where the asymptotic expansion takes place in $C_b^\infty(\R^n)$.
\vskip 0.3cm
This means that ${\mathcal  W}\varphi_0$ is a coherent state in $L^2(\R^n)$, centered at $(x_0,\xi_0)$, and from this point we can apply all the standard (and less standard) results of semiclassical analysis for scalar operators, in order to compute $e^{-itA/h}{\mathcal  W}\varphi_0$ (see, e.g., \cite{CoRo, Ha1, Ro1, Ro2} and references therein). In particular, we learn from \cite{CoRo} Theorem 3.1 (see also \cite{Ro2}), that, for any $N\geq 1$,
\be
\label{estcoro}
e^{-itA/h}{\mathcal  W}\varphi_0 = e^{i\delta_t/h}\sum_{k=0}^{3(N-1)}c_k(t;h)\Phi_{k,t} +{\mathcal  O}(e^{NC_0t}h^{N/2}),
\ee
where $\Phi_{k,t}$ is a (generalized) coherent state centered at $(x_t,\xi_t):=\exp tH_{a_0}(x_0,\xi_0)$, $\delta_t:= \int_0^t(\dot x_s\xi_s -a_0(x_s,\xi_s))ds + (x_0\xi_0 - x_t\xi_t)/2$, $C_0>0$ is a constant,   the coefficients $c_k(t;h)$'s are of the form,
\be
\label{ckth}
c_k(t;h)=\sum_{\ell =0}^{N_k}h^\ell c_{k,\ell}(t),
\ee
with $c_{k,\ell}$ universal polynomial with respect to $(\partial^\gamma a_0(x_t, \xi_t))_{|\gamma|\leq M_k}$, and where the estimate is uniform with respect to $(t,h)$ such that $0\leq t <T_{\Omega'}(x_0,\xi_0)$ and $he^{C_0t}$ remains bounded ($h>0$ small enough). In particular, (\ref{estcoro}) supplies an asymptotic expansion of $e^{-itA/h}{\mathcal  W}\varphi_0$ if one restricts to the values of $t$ such that $0\leq  t <<\ln\frac1{h}$.
\vskip 0.3cm
Now, applying ${\mathcal  W}^*$ to (\ref{estcoro}), and observing that ${\mathcal  W}^*\Phi_{k,t} = {\mathcal  V}^*(\Phi_{k,t}\tilde u_{L'+1}) = U_j^{-1}{\mathcal  V}^*_j(\Phi_{k,t}u_{L'+1,j})$, where $j=j(t)$ is chosen in such a way that $\exp tH_{a_0}(x_0,\xi_0)\in \Omega_j$, and where ${\mathcal  V}^*_j := U_j{\mathcal  V}^*U_j^{-1}$ is an $h$-admissible operator on $L^2(\Omega_j ;{\mathcal  H})$ (that is, becomes an $h$-admissible operator on $L^2(\R^n ;{\mathcal  H})$ once sandwiched by cutoff functions supported in $\Omega_j$), we deduce from (\ref{mainrepcs}),
\begin{theorem}\sl 
\label{WP}
Let $\varphi_0$ be as in (\ref{cohststart}), and let $T_{\Omega'}(x_0,\xi_0)$ defined in (\ref{defthamilt}). Then, there exists $C>0$ such that, for any $N\geq 1$, one has,
$$
e^{-itP/h}\varphi_0 = e^{i\delta_t/h}\sum_{k=0}^{3(N-1)} c_k(t ;h)\Phi_{k,t} U_{j(t)}^{-1}\tilde v_{k, j(t)}( x)+{\mathcal  O}(h^{N/4}),
$$
where $\Phi_{k,t}$ is a  coherent state centered at $(x_t,\xi_t):=\exp tH_{a_0}(x_0,\xi_0)$,  $j(t)\in \{1,\dots,r\}$ is such that $\exp tH_{a_0}(x_0,\xi_0)\in\Omega_{j(t)}$, $\tilde v_{k, j(t)} \in C^\infty (\Omega_{j(t)} ;{\mathcal  H})$, $c_k(t;h)$ is as in (\ref{ckth}), $\delta_t:= \int_0^t(\dot x_s\xi_s -a_0(x_s,\xi_s))ds + (x_0\xi_0 - x_t\xi_t)/2$, and where the estimate is uniform with respect to $(t,h)$ such that $h>0$ is small enough and $t\in [0,\min(T_{\Omega'}(x_0,\xi_0), C^{-1} \ln\frac1{h}))$.
\end{theorem} 
\begin{remark}\sl  Actually, the coherent state $\Phi_{k,t}$ is of the form,
$$
\Phi_{k,t}= c_k(t)f_k(x,\sqrt {h})h^{-n/4}e^{ix\xi_t/h -q_t(x-x_t)/h},
$$
where $c_k(t)$ is a normalizing factor, $f_k$ is polynomial in 2 variables, and $q_t$ is a $t$-dependent quadratic form with positive-definite real part, that can be explicitly computed by using a classical evolution involving the Hessian of $a_0$ at $(x_t,\xi_t)$  (see \cite{CoRo}). More precisely, one has $q_t(x)=-i\la \Gamma_tx , x\ra/2$ with $\Gamma_t=(C_t+iD_t)(A_t+iB_t)^{-1}$, where the $2n\times 2n$ matrix,
$$
F_t=\left(\begin{array}{cc}
A_t & B_t\\
C_t & D_t
\end{array}\right)
$$
is, by definition, the solution of the classical problem,
$$
\dot F_t = J{\rm Hess}\hskip 1pt a_0(x_t,\xi_t)F_t \quad ;\quad F(0) = I_{2n}.
$$
Here,  $J:=\left(\begin{array}{cc}
0& I_n\\
-I_n & 0
\end{array}\right)$, and ${\rm Hess}\hskip 1pt a_0$ stands for the Hessian of $a_0$.
(We are grateful to M. Combescure and D. Robert for having explained to us this construction and the main result of \cite{CoRo}.)
\end{remark}
\begin{remark}\sl  As in \cite{CoRo}, one can also consider more  general initial states, of the form,
$$
\varphi_0 (x) = e^{i(\xi_0 \cdot x - x_0\cdot  hD_x)/h}f\left(\frac{x}{\sqrt h}\right),
$$
where $f\in{\mathcal  S}(\R^n)$ (we refer to \cite{CoRo} Theorem 3.5 for more details). In the same way, a similar result can also be obtained  for oscillating initial states  of the form,
$$
\varphi_0 (x) = f(x)e^{iS(x) /h},
$$
where $f\in C_0^\infty (\R^n)$ and $S\in C^\infty (\R^n ;\R)$ (see \cite{CoRo} Remark 3.9).
\end{remark}
\begin{remark}\sl  In principle, all the terms of the asymptotic series can be computed explicitly by an inductive procedure (although, in practical, this task may result harder than expected since the simplifications are sometimes quite tricky). Indeed, all our constructions mainly rely on symbolic pseudodifferential calculus, that provides very explicit inductive formulas.  
\end{remark}
\chapter{Application to Polyatomic Molecules}
\label{PolMol}
\setcounter{equation}{0}%
\setcounter{theorem}{0}%
In this chapter, we apply all the previous results to the particular case of a polyatomic molecule with Coulomb-type interactions, imbedded in an electromagnetic field. Denoting by $x=(x_1,\dots ,x_n)\in \R^{3n}$ the position of the $n$ nuclei, and by $y=(y_1,\dots ,y_p)\in\R^{3p}$ the position of the $p$ electrons, the corresponding Hamiltonian takes the form,
\be
\label{physmolsys}
H=\sum_{j=1}^n \frac1{2M_j}(D_{x_j} - A(x_j))^2 + \sum_{k=1}^p \frac1{2m_k}(D_{y_k} - A(y_k))^2 + V(x,y),
\ee
where the magnetic potential $A$ is assumed to be in $C_b^\infty (\R^3)$,  and where the electric potential $V$ can be written as,
\be
\label{potmol}
V(x,y) =  V_{\rm nu} (x)+V_{\rm el}(y) + V_{\rm el\mbox{-}nu}(x,y) + V_{\rm ext} (x,y)= V_{\rm int}(x,y) + V_{\rm  ext}(x,y).
\ee
Here, $V_{\rm nu}$ ( resp. $V_{\rm el}$, resp. $V_{\rm el\mbox{-}nu}$) stands for sum of the nucleus-nucleus  (resp. electron-electron, resp. electron-nucleus) interactions, and $ V_{\rm  ext}$ stands for the external electric potential. Actually, our techniques can be applied to a slightly more general form of Hamiltonian (also allowing, somehow, a strong action of the magnetic field upon the nuclei), namely,
\be
H=\sum_{j=1}^n \frac1{2M_j}(D_{x_j} - a_jA_j(x))^2 + \sum_{k=1}^p \frac1{2m_k}(D_{y_k} - B_k(x,y))^2 + V(x,y),
\ee
where  $A_1,\dots, A_n$ (respectively  $B_1,\dots, B_p$) are assumed to be in $C_b^\infty (\R^n;\R)$ (respectively $C_b^\infty (\R^{n+p};\R)$), the $a_j$'s are  extra parameters, and $V$ is as in (\ref{potmol}) with,
\begin{eqnarray}
V_{\rm nu} (x)= \sum_{1\leq j < j'\leq n}\frac{\alpha_{j,j'}}{\vert x_j -x_{j'}\vert}\,\, &;&\,\,V_{\rm el}(y) = \sum_{1\leq k < k'\leq p}\frac{\beta_{k,k'}}{\vert y_k -y_{k'}\vert} \,\,;\nonumber\\
\label{coulomb}
V_{\rm el\mbox{-}nu}(x,y)= \sum_{\genfrac{}{}{0pt}{}{1\leq j\leq n}{1\leq k\leq p}} \frac{-\gamma_{j,k}}{\vert x_j -y_k\vert}\,\,& ;&\,\, V_{ext} \in C_b^\infty (\R^{n+p};\R),
\end{eqnarray}
$\alpha_{j,j'}, \beta_{k,k'}, \gamma_{j,k}> 0$ constant.  In fact, as in \cite{KMSW}, more general forms can be allowed for the interaction potentials, e.g., by replacing any function of the type $\vert z_j -z'_k\vert^{-1}$ (where the letters $z$ and $z'$ stand for $x$ or $y$ indifferently) by some $V_{j,k}(z_j -z'_k)$, where $V_{j,k}$ is assumed to be $\Delta$-compact on $L^2(\R^3)$ and to verify some estimates on its derivatives (see \cite{KMSW} Section 2). In the same way, one could also have admitted singularities of the same kind for the exterior potentials. However, here we keep the form (\ref{coulomb}) since it is more concrete and corresponds to the usual physical situation.
\vskip 0.2cm
Then, we consider the Born-Oppenheimer limit in the following sense: We set,
\be
\label{limsemicl}
M_j = h^{-2}b_j \, ;\, a_j = h^{-1} c_j + d_j,
\ee
and we consider the limit $h\to 0_+$ for some fix $b_j, m_k >0$, $c_j, d_j\in \R$. By scaling the time variable, too, the quantum evolution of the molecule is described by the Schr\"odinger equation,
$$
ih\frac{\partial\varphi}{\partial t} = P(h)\varphi,
$$
where,
\be
\label{pdehmol}
P(h):=\sum_{j=1}^n \frac1{2b_j}(hD_{x_j} - (c_j+hd_j)A_j(x))^2 + \sum_{k=1}^p \frac1{2m_k}(D_{y_k} - B_k(x,y))^2 + V(x,y).
\ee
In particular, we see that $P(h)$ satisfies to Assumptions (H1) and (H2), with,
\begin{eqnarray*}
&&{\boldsymbol\omega} = \sum_{j=1}^n \frac1{2b_j}(hD_{x_j} - (c_j+hd_j)A_j(x))^2,\\
&& \omega (x,\xi ;h ) = \sum_{j=1}^n \frac1{2b_j}\left[(\xi_j - (c_j+hd_j)A_j(x))^2 + ih(c_j + hd_j)(\partial_{x_j}A_j)(x)\right],\\
&& Q(x) =  \sum_{k=1}^p \frac1{2m_k}(D_{y_k} - B_k(x,y))^2 + V_{\rm el}(y) + V_{\rm el\mbox{-}nu}(x,y) + V_{\rm ext} (x,y),\\
&& W(x) = V_{\rm nu} (x).
\end{eqnarray*}
Now, following the terminology of \cite{KMSW}, we denote by
$$
{\mathcal  C}:= \bigcup_{\genfrac{}{}{0pt}{}{1\leq j,k\leq n}{j\not= k}} \{ x=(x_1,\dots,x_n)\in\R^{3n}\, ;\, x_j =x_k\}
$$
the so-called {\it collision set} of nuclei, and
we make on $Q(x)$ the following gap condition:
\\\\
$(H3')$ There exists a contractible bounded  open set $\Omega \subset \R^{3n}$ such that $\overline\Omega\cap{\mathcal  C} =\emptyset$, and, for all $x\in  \overline\Omega$, the $L'+L$ first values $\lambda_1(x), \dots ,\lambda_{L'+L}(x)$, given by the Mini-Max principle for $Q(x)$  on $L^2(\R^{3p})$, are discrete eigenvalues of $Q(x)$, and verify,
$$
\inf_{x\in \Omega} \dist \left(\sigma(Q(x))\backslash \{\lambda_{L'+1}(x), \dots ,\lambda_{L'+L}(x)\} , \{\lambda_{L'+1}(x), \dots ,\lambda_{L'+L}(x)\}\right) >0.
$$
\vskip 0.3cm
As it is well known (see \cite{CoSe}), under these assumptions, the  two spectral projections $\Pi_0^-(x)$ and $\Pi_0(x)$ of $Q(x)$, corresponding to $\{\lambda_{1}(x), \dots ,\lambda_{L'}(x)\}$ and $\{\lambda_{L'+1}(x), \dots ,\lambda_{L'+L}(x)\}$ respectively, are twice differentiable with respect to $x\in\Omega$. In particular, the whole assumption (H3) is indeed  satisfied in that case (and even with a slightly larger open subset of $\R^{3n}$).
\vskip 0.2cm
Now, in order to be able to apply the results of the previous chapters to this molecular Hamiltonian, it remains to construct a family $(\Omega_j, U_j(x))_{1\leq j\leq r}$ that verifies Assumption (H4). We do it by following \cite{KMSW}. 
\vskip 0.2cm
More precisely, for any fixed $x_0=(x_1^0,\dots,x_n^0)\in\R^{3n}\backslash {\mathcal  C}$, we choose $n$ functions $f_1, \dots, f_n\in C_0^\infty (\R^3; \R)$, such that,
$$
f_j(x_k^0) = \delta_{j,k}\,\,\, (1\leq j,k\leq n),
$$
and, for $x\in \R^{3n}$, $s\in \R^3$, and $y=(y_1,\dots,y_p)\in\R^{3p}$, we set,
\begin{eqnarray*}
&& F_{x_0}(x,s) := s + \sum_{k=1}^n (x_k -x_k^0)f_k(s)\in \R^3,\\
&& G_{x_0}(x,y):= (F_{x_0}(x,y_1), \dots, F_{x_0}(x,y_p))\in \R^{3p}.
\end{eqnarray*}
Then, by the implicit function theorem, for $x$ in a sufficiently small neighborhood $\Omega_{x_0}$ of $x_0$, the application $y\mapsto G_{x_0}(x,y)$ is a diffeomorphism of $\R^{3p}$, and we have,
\begin{eqnarray*}
&& x_k=F_{x_0}(x,x^0_k),\\
&& G_{x_0}(x,y)= y \mbox{ for } \vert y\vert\mbox{ large enough}.
\end{eqnarray*}
Now, for $v\in L^2(\R^{3p})$ and $x\in\Omega_{x_0}$, we define,
$$
U_{x_0}(x)v(y):= \vert {\rm det\hskip 1pt} d_yG_{x_0}(x,y)\vert^{\frac12} v(G_{x_0}(x,y))\vert,
$$
and we see that $U_{x_0}(x)$ is a unitary operator on $L^2(\R^{3p})$ that preserves both ${\mathcal  D}_Q =H^2(\R^{3p})$ and $C_0^\infty (\R^{3p})$. Moreover, denoting by $U_{x_0}$ the operator on $L^2(\Omega_{x_0}\times \R^{3p})$ induced by $U_{x_0}(x)$, we have the following identities:
\begin{eqnarray}
&& U_{x_0}hD_xU_{x_0}^{-1}=hD_x + hJ_1(x,y)D_y + hJ_2(x,y),\nonumber \\
&& U_{x_0}D_yU_{x_0}^{-1}= J_3(x,y)D_y + J_4(x,y),\nonumber\\
&& U_{x_0}\frac1{\vert y_k -y_k'\vert}U_{x_0}^{-1}=\frac1{\vert F_{x_0}(x,y_k) -F_{x_0}(x,y_k'),\vert}\nonumber\\
\label{singenlev}
&& U_{x_0}\frac1{\vert x_j -y_k\vert}U_{x_0}^{-1}= \frac1{\vert F_{x_0}(x,x_j^0) -F_{x_0}(x,y_k)\vert},
\end{eqnarray}
where the (matrix or operator-valued) functions $J_\nu$'s ($1\leq \nu\leq 4$) are all smooth on $\Omega_{x_0}\times \R^{3p}$. Indeed, denoting by $\tilde G_{x_0}(x,\cdot )$ the inverse diffeomorphism of $G_{x_0}(x,\cdot )$, one finds, 
\begin{eqnarray*}
J_1(x,y) &=& ({}^td_x\tilde G_{x_0})(x, y'=G_{x_0}(x,y)),\\
J_2(x,y)&=&\vert {\rm det\hskip 1pt} d_yG_{x_0}(x,y)\vert^{\frac12}D_x\left(\vert{\rm det\hskip 1pt} d_{y'}\tilde G_{x_0}(x, y')\vert^{\frac12}\right)\left\vert_{y'=G_{x_0}(x,y))}\right.,\\
J_3(x,y) &=& ({}^td_{y'}\tilde G_{x_0})(x, y'=G_{x_0}(x,y)),\\
J_4(x,y)&=&\vert {\rm det\hskip 1pt} d_yG_{x_0}(x,y)\vert^{\frac12}D_{y'}\left(\vert{\rm det\hskip 1pt} d_{y'}\tilde G_{x_0}(x, y')\vert^{\frac12}\right)\left\vert_{y'=G_{x_0}(x,y))}\right..
\end{eqnarray*}
The key-point in (\ref{singenlev}) is that the ($x$-dependent) singularity at $y_k=x_j$ has been replaced by the (fix) singularity at $y_k=x_j^0$.
Then, as in \cite{KMSW}, one can easily deduce that  the map $x\mapsto U_{x_0}Q(x)U_{x_0}^{-1}$ is in $C^\infty (\Omega_{x_0}; {\mathcal  L}(H^2(\R^{3p}), L^2(\R^{3p}))$. Moreover, so is the map  $x\mapsto U_{x_0}\Delta_yU_{x_0}^{-1}$, and we also see that $U_{x_0}{\boldsymbol\omega}U_{x_0}^{-1}$ can be written as in (\ref{conjomega}) (with $\Omega_{x_0}$ instead of $\Omega_j$, $m=2$, and $Q_0 = -\Delta_y + C_0$, $C_0 >0$ large enough). Indeed,  with the notations of (\ref{singenlev}), and setting ${\mathcal  J}(x) =( {\mathcal  J}_1(x),\dots, {\mathcal  J}_n(x)):= J_1(x,y)D_y + J_2(x,y)$,  we have,
\begin{eqnarray}
U_{x_0}{\boldsymbol{\omega}}U_{x_0}^{-1} &=& \sum_{k=1}^n\frac1{2b_k}(hD_{x_k} +h{\mathcal  J}_k(x) -(c_k+hd_k)A_k(x))^2\nonumber\\
\label{omegconj}
&=& {\boldsymbol{\omega}}+h\sum_{k=1}^n\frac1{b_k}{\mathcal  J}_k(hD_{x_k}-c_kA_k)\\
&& \hskip 1cm+ h^2\sum_{k=1}^n\frac1{2b_k}({\mathcal  J}_k^2-i(\nabla_x{\mathcal  J}_k) -2d_kA_k{\mathcal  J}_k).\nonumber
\end{eqnarray}
\vskip 0.2cm
To complete the argument, we just observe that the previous construction can be made around any point $x_0$ of $\overline\Omega$, and since this set is compact, we can cover it by a finite family $\tilde\Omega_1,\dots,\tilde\Omega_r$ of open sets such that each one corresponds to some $\Omega_{x_0}$ as before. Denoting also $U_1(x),\dots,U_r(x)$ the corresponding operators $U_{x_0}(x)$, and setting $\Omega_j =\tilde\Omega_j\cap \Omega$, we can conclude that the family $(\Omega_j ,U_j(x))_{1\leq j\leq r}$ verifies (H4) with ${\mathcal  H}_\infty = C_0^\infty (\R^{3p})$. As a consequence, we can apply to this model all the results of the previous chapters, and thus, we have proved,
\begin{theorem}\sl 
\label{MAINTHmol}
Let $P(h)$ be as in (\ref{pdehmol}) with $V$ given by (\ref{potmol}) and (\ref{coulomb}), $A_1,\dots,A_n\in C_b^\infty (\R^n ;\R)$, and $B_1,\dots, B_p\in C_b^\infty (\R^{n+p};\R)$.  Assume also (H3'). Then, the conclusions of Theorem \ref{MAINTH} are valid for $P=P(h)$.
\end{theorem}
We also observe that, in this case, we have,
$$
\omega (x,\xi ;h)= \omega_0(x,\xi ) + h\omega_1(x,\xi ) + h^2\omega_2(x),
$$
with,
\begin{eqnarray}
\omega_0(x,\xi )&=&\sum_{k=1}^n \frac1{2b_k}(\xi_k - c_kA_k(x))^2 \nonumber\\
\label{decomeg}
\omega_1(x,\xi )&=&\sum_{k=1}^n \frac1{2b_k}\left[2d_kA_k(x)(c_kA_k(x)-\xi_k)+ic_k(\partial_{x_k}A_k)(x) \right]\\
\omega_2(x)&=& \sum_{k=1}^n \frac1{2b_k}\left[d_k^2A_k(x)^2+id_k(\partial_{x_k}A_k)(x) \right].\nonumber
\end{eqnarray}
In particular, the conditions (\ref{coeffprinc}) and (\ref{coefclas}) are satisfied, and thus, we also have,
\begin{theorem}\sl 
Let $P(h)$ be as in (\ref{pdehmol}) with $V$ given by (\ref{potmol}) and (\ref{coulomb}), $A_1,\dots,A_n\in C_b^\infty (\R^n ;\R)$, and $B_1,\dots, B_p\in C_b^\infty (\R^{n+p};\R)$. Assume also (H3') and $L=1$. Then, the conclusions of Corollary \ref{CORT} and Theorem \ref{WP} are valid for $P=P(h)$.
\end{theorem}
Moreover, concerning the symbol of the effective Hamiltonian, in that case we have,
\begin{theorem}\sl 
\label{remsymbol}
Let $P(h)$ be as in (\ref{pdehmol}) with $V$ given by (\ref{potmol}) and (\ref{coulomb}), $A_1,\dots,A_n\in C_b^\infty (\R^n ;\R)$, and $B_1,\dots, B_p\in C_b^\infty (\R^{n+p};\R)$. Assume also (H3') and $L=1$. Then, the symbol $a(x,\xi ;h)$ of the effective Hamiltonian verifies,
$$
a(x,\xi ;h)= a_0(x,\xi ) +ha_1(x,\xi) + h^2a_2(x,\xi) +{\mathcal  O}(h^3),
$$
with, for $(x,\xi)\in T^*(\Omega)$,
\begin{eqnarray*}
a_0(x,\xi)&=& \omega_0(x,\xi ) +\lambda_{L'+1} (x) +W(x);\\
a_1(x,\xi)&=& \omega_1(x,\xi )  -i \nabla_\xi\omega_0(x,\xi)\la \nabla_{x}u(x),u(x)\ra\\
a_2(x,\xi)&=&  \sum_{k=1}^n\frac1{2b_k}\la (\xi_k-d_kA_k(x))^2u(x), u(x)\ra \\
&& +\sum_{k,\ell =1}^n \frac1{b_kb_\ell}(\xi_k - c_kA_k)(\xi_\ell -c_\ell A_\ell)\la R'(x,\lambda (x))\nabla_{x_k}u,\nabla_{x_\ell}u\ra,
\end{eqnarray*}
where $\omega_0$ and $\omega_1$ are defined in (\ref{decomeg}), and 
$$
R'(x,\lambda (x)):=\Pi_0^\perp (x) (\lambda (x) - Q(x))^{-1}\Pi_0^\perp (x),
$$
 is the reduced resolvent of $Q(x)$.
\end{theorem}
{\em Proof -- } A possible proof may consist in using Proposition \ref{symhameff}. Then,  observing (with the notations of (\ref{omegconj})) that, by definition, 
\be
\label{defJ}
 {\mathcal  J}= U_{x_0}D_xU_{x_0}^{-1} - D_x,
\ee
and, exploiting the fact that  
the $(L'+1)$-th normalized eigenstate $u(x)$ of $Q(x)$ is a twice differentiable function of $x$ with values in $L^2(\R^{2p})$ (see , e.g., \cite{CoSe}, but this is also an easy consequence of (\ref{defJ}) and the fact that $x\mapsto U_{x_0}(x)u(x)$ is smooth), and setting $v(x)=U_{x_0}(x)u(x)$, one can write,
$$
\la {\mathcal  J}v,v\ra_{\mathcal  H} = \la D_xu,u\ra_{\mathcal  H} - \la D_xv,v\ra_{\mathcal  H}.
$$
As a consequence, one also finds,
$$
\sum_{k=1}^n\frac1{b_k}(\xi_k-c_kA_k)\la {\mathcal  J}_kv,v\ra_{\mathcal  H}-i\la \nabla_\xi\omega_0\nabla_{x}v,v\ra_{\mathcal  H} = -i\la \nabla_\xi\omega_0 \nabla_{x}u,u\ra_{\mathcal  H} .
$$
where $\omega_\ell$ ($0\leq\ell\leq 2$) are defined in (\ref{decomeg}), and this permits to make appear many cancellations in the expression of $a(x,\xi ;h)$ given in Proposition \ref{symhameff}, leading to the required formulas. 
\vskip 0.3cm
However, there is a much simpler way to prove it, using directly the expressions (\ref{hameffl1})  given in  Proposition \ref{symhameff} for the operator $A$. Indeed, since in our case $x\mapsto u(x)$ is twice differentiable, for all $w\in L^2(\R^{n+p})$, we can write,
$$
[D_x, \tilde\Pi_0]w = -i\la w, \nabla_xu(x)\ra u(x) -i \la w, u(x)\ra\nabla_xu(x),
$$
and, for all $w\in C^1(\R^{3n};L^2(\R^{3p} ))$,
\begin{eqnarray*}
[D_x^2, \tilde\Pi_0]w&=& [D_x, \tilde\Pi_0] \cdot D_xw+D_x\cdot [D_x, \tilde\Pi_0]w\\
&=& -2i\la D_xw, \nabla_xu(x)\ra u(x) -2i \la D_xw, u(x)\ra\cdot \nabla_xu(x)\\
&& -\la w, \nabla_xu(x)\ra \cdot \nabla_x u(x) - \la w, u(x)\ra \nabla_x\cdot\nabla_x u(x).
\end{eqnarray*}
This permits to write explicitly the operator $[\tilde\Pi_0, \tilde P] = [\tilde\Pi_0, {\boldsymbol{\omega}}]$ as,
\begin{eqnarray*}
[\tilde\Pi_0, \tilde P]  w &=&  ih \sum_{k=1}^n\frac1{b_k}\la (hD_{x_k}-(c_k+hd_k)A_k)w, \nabla_{x_k}u(x)\ra u(x) \\
&&+ih \sum_{k=1}^n\frac1{b_k} \la (hD_{x_k}-(c_k+hd_k)A_k)w, u(x)\ra \cdot \nabla_{x_k}u(x)\\
&&+h^2 \sum_{k=1}^n\frac1{2b_k}\left( \la w, \nabla_{x_k}u(x)\ra \cdot \nabla_{x_k} u(x) +\la w, u(x)\ra \nabla_{x_k}^2 u(x)\right).
\end{eqnarray*}
In particular, taking $w=Z_1^* \alpha (x) = \alpha (x)u(x)$, $\alpha \in H^1(\R^{3n})$), and using the fact that $R'(x,\lambda (x))u(x) =0$, one finds,
\begin{eqnarray*}
R'(x,\lambda (x)) [\tilde\Pi_0, \tilde P]  Z_1^*\alpha  &=& ih  \sum_{k=1}^n\frac1{b_k}\left((hD_{x_k}-c_kA_k)\alpha\right) R'(x,\lambda (x))\nabla_{x_k}u(x)\\
&& \hskip 2cm +{\mathcal  O}(h^2\Vert\alpha\Vert),
\end{eqnarray*}
and then,
\begin{eqnarray*}
&& Z_1[\tilde P, \tilde\Pi_0] R'(x,\lambda (x)) [\tilde\Pi_0, \tilde P] Z_1^*\alpha  \\
&& = h^2  \sum_{k,\ell =1}^n\frac1{b_kb_\ell}\left((hD_{x_k}-c_kA_k)(hD_{x_\ell}-c_\ell A_\ell)\alpha\right) \times \\
&& \hskip 1cm \times \la R'(x,\lambda (x))\nabla_{x_k}u(x),\nabla_{x_\ell}u(x)\ra  +{\mathcal  O}(h^3\Vert\alpha\Vert),
\end{eqnarray*}
This obviously permits to compute the principal symbol of the partial differential operator $A_2$ appearing in (\ref{hameffl1}). The (full) symbol of $A_1=Z_1\tilde PZ_1^*$ is even easier to compute, and the result follows.
\hfill$\bullet$
\begin{remark}\sl
The smoothness with respect to $x$ of all the coefficients appearing in $a(x,\xi ;h)$ is {\it a priori} known, but can also be recovered directly by using (\ref{defJ}). For instance, writing $\la \nabla_xu(x), u(x)\ra$ as, $$
\la \nabla_xu(x), u(x)\ra = \la \nabla_xU_{x_0}u(x), U_{x_0}u(x)\ra + i\la {\mathcal  J}(x)U_{x_0}u(x), U_{x_0}u(x)\ra,
$$
 permits to see its smoothness near $x_0$.
\end{remark}
\begin{remark}\sl
Using the expression of $A_3$ appearing in (\ref{hameffl1}), one could also compute the next term (i.e., the $h^3$-term) in $a(x,\xi ;h)$.
\end{remark}
\begin{remark}\sl
Analogous formulas can be obtained in a very similar way in the case where $L$ is arbitrary but $\lambda_{L'+1}=\dots =\lambda_{L'+L}$.
\end{remark}
\begin{remark}\sl
Although we did not do it here, we can also treat the case of unbounded magnetic potential (e.g., constant magnetic field). Then, the estimates on the coefficients $c_\alpha$'s in Assumption (H1) are not satisfied anymore, but, since we mainly work in a compact region of the $x$-space, it is  clear that an adaptation of our arguments  lead to the same results.
\end{remark}
\begin{remark}\sl  In the case of a free molecule (or, more generally, if the external electromagnetic field is invariant under the translations of the type $(x,y)\mapsto (x_1+\alpha, \dots, x_n+\alpha, y_1+\alpha, \dots ,y_p+\alpha)$ for any $\alpha\in\R^3$), one can factorize the quantum motion, e.g., by using the so-called {\it center of mass of the nuclei} coordinate system, as in \cite{KMSW}. Then, denoting by $R$ the position of the center of mass of the nuclei, the operator takes the form,
$$
P(h) = H_0(D_R) + P'(h) + h^2p(D_y),
$$
where $H_0(D_R)$ stands for the quantum-kinetic energy of the center of mass of the nuclei, $P'(h)$ has a form similar to that of $P(h)$ in (\ref{pdehmol}) (but now, with $x\in\R^{3(n-1)}$ denoting the {\it relative} positions of the nuclei), and $p(D_y)$ is a PDO of order 2 with respect to $y$, with constant coefficients (the so-called {\it isotopic term}). Therefore, one obtains the factorization,
\be
\label{factor}
e^{-itP(h)/h} = e^{-itH_0(D_R)/h}e^{-it(P'(h)+h^2p(D_y))/h},
\ee
and it is easy to verify that our previous constructions can be performed with $Q(x)$ replaced by $Q(x) + h^2p(D_y)$. In particular, under the same assumptions as in Theorem \ref{MAINTHmol}, the quantum evolution under $P'(h) + h^2p(D_y)$ of an initial state $\varphi_0$ verifying (\ref{condloc}) with $P$ replaced by $P'(h)$ (that is, a much weaker assumption) can be expressed in terms of the quantum evolution associated to a $L\times L$ matrix of $h$-admissible operators on $L^2(\R^{3(n-1)})$. In that case, (\ref{factor}) provides a way to reduce the evolution of $\varphi_0$ under $P(h)$, too.
\end{remark}

\appendix
\chapter{Smooth Peudodifferential Calculus with Operator-Valued Symbol}
\label{stpdo}
\setcounter{equation}{0}%
\setcounter{theorem}{0}%
We recall the usual definition of $h$-admissible operator  with operator-valued symbol. In some sense, this corresponds to a simple case of the more general definitions given in \cite{Ba, GMS}. For $m\in\R$ and $\mathcal  H$ a Hilbert space,  we denote by $H^m(\R^n;{\mathcal  H})$ the standard $m$-th order Sobolev space on $\R^n$ with values in ${\mathcal  H}$.
\begin{definition}\sl Let $m\in\R$ and let  ${\mathcal  H}_1$ and  ${\mathcal  H}_2$ be two Hilbert space. An operator 
$A=A(h)\, :\, H^m({\R^n}; {\mathcal  H}_1)\rightarrow L^2({\R^n};{\mathcal  H}_2)$ 
with $h\in (0,h_0 ]$  is called $h$-admissible (of degree $m$) if,  for any $N\geq 1$,
\be
\label{defadmiss}
A(h)=\sum_{j=0}^Nh^j{\rm Op}_h(a_j(x,\xi ;h))+h^NR_N(h),
\ee
where $R_N$ is uniformly bounded from $H^m({\R^n}; {\mathcal  H}_1)$ to $L^2({\R^n};{\mathcal  H}_2)$ for $h\in (0,h_0 ]$, and, for all $h>0$ small enough, $a_j\in  C^{\infty}(T^*{\R}^n; {\mathcal  L}({\mathcal  H}_1;{\mathcal  H}_2))$, with 
\be
\label{estsymb}
\Vert \partial^{\alpha}a_j(x,\xi ;h)\Vert_{{\mathcal  L}({\mathcal  H}_1;{\mathcal  H}_2)}\leq C_{\alpha}\la \xi\ra^m
\ee 
for all $\alpha\in  {\Z_+}^{2n}$ and some positive constant $C_\alpha$, uniformly for 
$(x,\xi)\in  T^*{\R}^n$ and $h>0$ small enough. In that case, the formal series,
\be
\label{formalsymb}
a(x,\xi ;h) = \sum_{j\geq 0} h^ja_j(x,\xi ;h),
\ee
is called the {\it symbol} of $A$ (it can be resummed up to a remainder in ${\mathcal  O}(h^\infty \la \xi\ra^m)$ together with all its derivatives).
Moreover, in the case $m=0$ and ${\mathcal  H}_2 ={\mathcal  H}_1$, $A$ is called a (bounded) $h$-admissible operator on $L^2(\R^n ;{\mathcal  H}_1)$.
\end{definition}
Here, we have denoted by ${\rm Op}_h(a)$ the standard quantization of a symbol $a$, defined by the following formula:
\be
\label{quantfiorm}
{\rm Op}_h(a)u (x) :=\frac1{(2\pi h)^n}\int e^{i(x-y)\xi /h}a\left(x ,\xi\right) u(y)dyd\xi,
\ee
valid for any tempered distribution $u$, and where the integral has to be interpreted as an oscillatory one. Actually, by the Calder\' on-Vaillancourt Theorem (see, e.g., \cite{GMS,DiSj1, Ma2, Ro1}, and below), the estimate (\ref{estsymb}) together with the quantization formula (\ref{quantfiorm}), permit to define ${\rm Op}_h(a)$ as a bounded operator $H^m({\R^n}; {\mathcal  H}_1)\rightarrow L^2({\R^n};{\mathcal  H}_2)$. Let us also observe that, very often, the formal series (\ref{formalsymb}) are indeed identified with one of their resummations (and thus, the symbol is considered as a function, rather than a formal series). Indeed, since the various resummations (together with all their derivatives) differ by uniformly ${\mathcal  O}(h^\infty \la \xi \ra^m)$ terms, in view of (\ref{defadmiss}) and the Calder\' on-Vaillancourt Theorem, it is clear that this has no real importance.
\vskip 0.2cm
As it is well known (see, e.g., \cite{Ba, DiSj1, GMS, Ma2}), with such a type of quantization is associated a full  and explicit  symbolic calculus that permits to handle these operators in a very easy and pleasant way. In particular, we have the following results:
\begin{proposition}[Composition]\sl 
\label{sharp}
Let $A$ and $B$ be two bounded $h$-admissible operators on $L^2(\R^n ;{\mathcal  H}_1)$, with respective symbols $a$ and $b$. Then, the composition $A\circ B$ is an $h$-admissible operators on $L^2(\R^n ;{\mathcal  H}_1)$, too, and its symbol $a\sharp b$ is given by the formal series,
$$
a\sharp b (x,\xi ;h) =\sum_{\alpha \in\Z_+^n} \frac{h^{\vert\alpha\vert}}{i^{\vert\alpha\vert}\alpha !}\partial_\xi^\alpha a(x,\xi ;h)\partial_x^\alpha b(x,\xi ;h).
$$
\end{proposition}
\begin{remark}\sl There is a similar result for the composition of unbounded $h$-admissible operators, but it requires more conditions on the remainder $R_N(h)$ appearing in (\ref{defadmiss}) (see \cite{Ba, GMS}).
\end{remark}
\begin{proposition}[Parametrix]\sl Let $A$ be a bounded $h$-admissible operator on $L^2(\R^n ;{\mathcal  H}_1)$, such that any resummation $a$ of its symbol is elliptic,  in the sense that $a(x,\xi ;h)$ is invertible on $\mathcal  H_1$ for any $(x,\xi ;h)$, and its inverse verifies,
$$
\Vert a(x,\xi ;h)^{-1}\Vert_{{\mathcal  L}({\mathcal  H}_1)}={\mathcal  O}(1),
$$
uniformly for 
$(x,\xi)\in  T^*{\R}^n$ and $h>0$ small enough. Then, $A$ is invertible on $L^2(\R^n ;{\mathcal  H}_1)$, its inverse $A^{-1}$ is $h$-admissible, and its symbol $b$ verifies,
$$
b=a^{-1} + hr,
$$
with $r=\sum_{j\geq 0}h^jr_j$, $\Vert \partial^\alpha r_j\Vert_{{\mathcal  L}({\mathcal  H}_1)}={\mathcal  O}(1)$ uniformly.
 \end{proposition}
 \begin{remark}\sl
It is easy to see that the ellipticity of any resummation of the symbol is equivalent to the ellipticity of the function $a_0(x,\xi ;h)$ appearing in (\ref{defadmiss}) (and thus, to the ellipticity of at least one resummation).
 \end{remark}
 \begin{remark}\sl
Of course, the $r_j$'s can actually be all determined recursively, by using the identity $a\sharp b =1$ (this gives a possible choice for them, but this choice is not unique since we have allowed them to depend on $h$). 
 \end{remark}
\begin{proposition}[Functional Calculus]\sl Let $A$ be a self-adjoint $h$-admissible operator on $L^2(\R^n ;{\mathcal  H}_1)$, and let $f\in C_0^\infty (\R)$. Then, $f(A)$ is $h$-admissible, and its symbol $b$ verifies,
$$
b=f(\re a)+hr,
$$
where $\re a := (a+a^*)/2$, and $r=\sum_{j\geq 0}h^jr_j$, $\Vert \partial^\alpha r_j\Vert_{{\mathcal  L}({\mathcal  H}_1)}={\mathcal  O}(1)$ uniformly.
 \end{proposition}
 \begin{proposition}[Calder\'on-Vaillancourt Theorem]\sl Let $a=a(x,\xi)$ be in $ C^{\infty}(T^*{\R}^n; {\mathcal  L}({\mathcal  H}_1;{\mathcal  H}_2))$, such that, for all $\alpha\in\Z_+^{2n}$,
$\Vert \partial^{\alpha}a(x,\xi )\Vert_{{\mathcal  L}({\mathcal  H}_1;{\mathcal  H}_2)}$ is uniformly bounded on $T^*{\R}^n$. Then, ${\rm Op}_h(a)$ (defined, e.g., on ${\mathcal S}(\R^n; {\mathcal  H}_1)$) extends to a bounded operator : $L^2(\R^n; {\mathcal  H}_1)\rightarrow L^2(\R^n; {\mathcal  H}_2)$, and  there exist two constants $C_n$ and $M_n$, depending only on the dimension $n$, such that,
$$
\Vert {\rm Op}_h(a)\Vert_{{\mathcal L}(L^2(\R^n; {\mathcal  H}_1);L^2(\R^n; {\mathcal  H}_2))}
\leq C_n\sum_{|\alpha|\leq M_n}\sup_{T^*\R^n}|\partial^\alpha a(x,\xi)|.
$$
 \end{proposition}
\chapter{Propagation of the Support}
\label{psupp}
\setcounter{equation}{0}%
\setcounter{theorem}{0}%
\begin{theorem}\sl 
\label{th:appA} Let $P$ be as in (\ref{P}) with (H1)-(H2), and
let $K_0$ be a compact subset of ${\R}^n_x$, $f\in C_0^\infty (\R)$ and $\varphi_0\in  L^2({\R}^{n};{\mathcal  H})$, such that $\Vert\varphi_0\Vert =1$, and,
$$
\Vert (1-f(P))\varphi_0\Vert_{L^2(\R^n ;{\mathcal  H})} + \Vert \varphi_0\Vert_{L^2(K_0^c ;{\mathcal  H})}={\mathcal  O}(h^{\infty}).
$$
Then,  for any $\varepsilon >0$, any $T>0$, and any $g\in C_0^\infty (\R)$ such that $gf=f$, the compact set defined by,
$$
K_{T,\varepsilon }:=\{ x\in \R^n\, ;\, \dist (x,K_0)\leq \varepsilon +C_1T\},
$$
with
$$
C_1:=\frac12\Vert \nabla_\xi\omega (x,hD_x)g(P)\Vert,
$$
verifies,
$$
\sup_{t\in [0,T]}\Vert e^{-itP/h}\varphi_0\Vert_{L^2(K_{T,\varepsilon}^c ;{\mathcal  H})}={\mathcal  O}(h^{\infty}),
$$
as $h\to 0$.
\end{theorem}
{\em Proof -- }
First, we need the following lemma:
\begin{lemma}\sl 
\label{premlemappB}
For any $\i\in  C_b^{\infty}({\R}^n)$, such that ${\rm supp}\i\subset K_0^c$, and
for any $g\in  C_0^{\infty}({\R})$, one has,
$$\Vert \i(x) g(P)\varphi_0\Vert={\mathcal  O}(h^{\infty}).$$
\end{lemma}
{\em Proof -- }  Consider a sequence $(\i_j)_{j\in\N}\subset 
C_b^{\infty}({\R}^n)$,  ${\rm supp}\i_j\subset K_0^c$ and such that
$$\i_{j+1}\i_j=\i_j,\quad \i_j\i =\i.$$
Then, in view of (\ref{calfonc}), it is sufficient to show that, for any $N\geq 0$, 
$$
\Vert \i_j(x) (P-\lambda)^{-1}\varphi_0\Vert={\mathcal  O}(h^{N}\vert
\im\lambda\vert^{-(N+1)}),
$$
uniformly as $h,\vert \im\lambda\vert\to 0_+$.
\vskip 0.3cm
We set,
$u_j=\i_j(x) (P-\lambda)^{-1}\varphi_0$,
and we observe that,  for all $j\in  {\N}$, one has $\Vert u_j\Vert={\mathcal  O}(\vert
\im\lambda\vert^{-1})$. By induction on $N$, 
let us suppose, for all $j\in \N$,
$$
\Vert \i_j(x) (P-\lambda)^{-1}\varphi_0\Vert={\mathcal  O}(h^{N}\vert
{\im}\lambda\vert^{-(N+1)}).
$$
Since $\i_{j+1}=1$ on $\supp \i_j$, and $P$ is differential in $x$, we have,
$$(P-\lambda)u_j=\i_j\varphi_0+[P, \i_j]\i_{j+1}(P-\lambda)^{-1}\varphi_0,$$
and  thus,
$$u_j=(P-\lambda)^{-1}\i_j\varphi_0+(P-\lambda)^{-1}[\boldsymbol{\omega},
\i_j]u_{j+1}.$$
Now, by assumption, we have $\Vert \i_j\varphi_0\Vert={\mathcal  O}(h^{\infty})$,
and  therefore, $\Vert (P-\lambda)^{-1}\i_j\varphi_0\Vert={\mathcal  O}(h^{\infty}\vert
{\im}\lambda\vert^{-1})$. Moreover, using (H1)-(H2), it is easy to see that the operator $\vert\im\lambda \vert  h^{-1}(P-\lambda)^{-1}[\boldsymbol{\omega},
\i_j]$ is uniformly bounded on $L^2(\R^n;{\mathcal  H})$.
Hence, using the induction hypothesis, we obtain,
$$\Vert u_j\Vert={\mathcal  O}(h^{\infty}\vert
{\im}\lambda\vert^{-1})+{\mathcal  O}(h^{N+1}\vert
{\im}\lambda\vert^{-(N+2)})={\mathcal  O}(h^{N+1}\vert
{\im}\lambda\vert^{-(N+2)})$$
for any $j\in  {\N}$, and the lemma follows.
\hfill$\bullet$
\vskip 0.2cm
Now, for any $F\in  C^{\infty}({\R}_+\times {\R}_x^n ;\R)$,
let us compute the quantity, 
\begin{eqnarray}
&&\partial_t\la F(t,x)f(P)e^{-itP/h}\varphi_0, f(P)e^{-itP/h}\varphi_0\ra\nonumber\\
&&\hskip 1.5cm={\re}\langle
(\partial_t F -ih^{-1}FP)f(P)e^{-itP/h}\varphi_0, f(P)e^{-itP/h}\varphi_0\rangle\nonumber\\
&&\hskip 1.5cm=\langle (\partial_t F
-\frac{i}{2h}[F,P])f(P)e^{-itP/h}\varphi_0, f(P)e^{-itP/h}\varphi_0\rangle\nonumber\\
\label{derivt}
&&\hskip 1.5cm=\langle (\partial_t F
+\frac{i}{2h}[\boldsymbol{\omega}, F])f(P)e^{-itP/h}\varphi_0, f(P)e^{-itP/h}\varphi_0\rangle.
\end{eqnarray}
Then, we fix $g\in C_0^\infty (\R)$ such that $gf=f$, and, for $j\in \N$, we set,
\be
\label{defFj}
F_j (t,x):= \varphi_j(\dist (x,K_0) -C_1t),
\ee
where $C_1=\frac12\Vert \nabla_\xi\omega (x,hD_x)g(P)\Vert$, and the $\varphi_j$'s are in $C_b^\infty (\R;\R_+)$ with support in $[\varepsilon ,+\infty)$, verify $\varphi_j(s) =1$ for $s\geq \varepsilon + \frac1{j}$, $\varphi_{j+1} =1$ near $\supp\varphi_j$, and are such that, 
$$
\varphi_j':=\phi_j^2\geq 0\,\, {\rm with}\,\, \phi_j\in C_b^\infty (\R;\R).
$$
In particular,  $F_j\in C_b^\infty({\R}_+\times {\R}_x^n ;\R_+)$, and, setting $d(x):=\dist (x,K_0)$, we have,
$$
\nabla_x F_j=\varphi_j'(d(x) -C_1t)\nabla d(x) ,
\quad\partial_t F_j=-C_1\varphi_j'(d(x) -C_1t).
$$ 
Moreover, since $\boldsymbol{\omega}=\omega (x,hD_x)$ is a differential operator  with respect $x$, of degree $m$, we see that,
\be
\label{commFomeg}
\frac{i}{h}[\boldsymbol{\omega}, F_j]=\nabla_x F_j\cdot\nabla_\xi\omega (x,hD_x)+ hR_j,
\ee
where $R_j=R_j(t, x,hD_x)$ is a differential operator of degree $m-2$ in $x$, with coefficients in $C_b^\infty({\R}_+\times {\R}_x^n)$ and supported in $\{ F_{j+1}=1\}$. 
\begin{lemma}\sl  
\label{lemmappendB}
For any $N\geq 1$,
$$
\Vert R_j f(P)u\Vert ={\mathcal  O}(\sum_{k=0}^Nh^{k}\Vert F_{j+k+1}f(P)u\Vert  + h^{N+1}\Vert u\Vert).
$$
\end{lemma}
{\em Proof -- }  We write,
$$
R_jf(P)=R_jF_{j+1}f(P)
= R_jg(P)F_{j+1}f(P)+ R_j[F_{j+1},g(P)]f(P).
$$
Then, using (\ref{calfonc})
and the fact that $[P, F_{j+1}]=[\boldsymbol{\omega},F_{j+1}]$, we obtain,
\begin{eqnarray*}
&&R_j[F_{j+1},g(P)]\\
 &&\hskip 0.5cm =\frac1{\pi}\int{\overline{\partial}}\tilde g(z)
R_j(P-z)^{-1}[\boldsymbol{\omega},F_{j+1}](P-z)^{-1} dz\; d\bar z\\
&&\hskip 0.5cm =\frac1{\pi}\int{\overline{\partial}}\tilde g(z)
R_j(P-z)^{-1}[\boldsymbol{\omega},F_{j+1}]F_{j+2}(P-z)^{-1} dz\; d\bar z\\
&&\hskip 0.5cm=\frac1{\pi}\int{\overline{\partial}}\tilde g(z)
R_j(P-z)^{-1}[\boldsymbol{\omega},F_{j+1}](P-z)^{-1}F_{j+2} dz\; d\bar z \\
&&\hskip 1cm +\frac1{\pi}\int{\overline{\partial}}\tilde g(z)
R_j(P-z)^{-1}[\boldsymbol{\omega},F_{j+1}](P-z)^{-1}[\boldsymbol{\omega},F_{j+2}](P-z)^{-1} dz\; d\bar z,
\end{eqnarray*}
and thus, by iteration,
\begin{eqnarray*}
&&R_j[F_{j+1},g(P)] \\
&&\hskip 0.5cm = \sum_{k=1}^N\frac1{\pi}\int {\overline{\partial}}\tilde g(z)
R_j(P-z)^{-1} \left( \prod_{\ell =1}^k \left([\boldsymbol{\omega},F_{j+\ell}](P-z)^{-1}\right)\right) F_{j +k+1} dz\; d\bar z\\
&&\hskip 2cm +\frac1{\pi}\int {\overline{\partial}}\tilde g(z)
R_j(P-z)^{-1} \prod_{\ell =1}^{N+1} \left([\boldsymbol{\omega},F_{j+\ell}](P-z)^{-1}\right)  dz\; d\bar z.
\end{eqnarray*}
Since $\Vert R_j(P-z)^{-1}\Vert ={\mathcal  O}(1)$ and $\Vert [\boldsymbol{\omega},F_{j+\ell}](P-z)^{-1}\Vert ={\mathcal  O}(h)$, the result follows. {}$\,$ \hfill$\bullet$
\vskip 0.3cm
As a consequence, we deduce from (\ref{commFomeg}),
\begin{eqnarray*}
&&\frac{i}{h}[\boldsymbol{\omega}, F_j]f(P)e^{-itP/h}\varphi_0\\
&&\hskip 0.2cm =\varphi_j'(d(x)
 -C_1t) \nabla d(x)\nabla_\xi\omega (x,hD_x)f(P)e^{-itP/h}\varphi_0\\
&&\hskip 4cm+{\mathcal  O}(\sum_{k=0}^Nh^{k+1}\Vert F_{j+k+1}f(P)e^{-itP/h}\varphi_0\Vert  + h^{N+2})\\
&&\hskip 0.2cm =\phi_j(d(x)
-C_1t) \nabla d(x)\nabla_\xi\omega (x,hD_x)g(P) \phi_j(d(x)
-C_1t) f(P)e^{-itP/h}\varphi_0\\
&& \hskip 0.5cm+ \phi_j(d(x)
-C_1t) )\nabla d(x)[ \phi_j(d(x)
-C_1t),\nabla_\xi\omega (x,hD_x)] f(P)e^{-itP/h}\varphi_0\\
&& \hskip 0.5cm+ \phi_j(d(x)
-C_1t) )\nabla d(x)\nabla_\xi\omega (x,hD_x)[ \phi_j(d(x)
-C_1t),g(P)] f(P)e^{-itP/h}\varphi_0\\
&&\hskip 4cm+{\mathcal  O}(\sum_{k=0}^Nh^{k+1}\Vert F_{j+k+1}f(P)e^{-itP/h}\varphi_0\Vert  + h^{N+2}),
\end{eqnarray*}
and thus, since $\phi_j$ is supported in $\{ F_{j+1}=1\}$, as in the proof of  Lemma \ref{lemmappendB}, we obtain,
\begin{eqnarray*}
&&\frac{i}{h}[\boldsymbol{\omega}, F_j]f(P)e^{-itP/h}\varphi_0\\
&&\hskip 0.2cm =\phi_j(d(x)
-C_1t) \nabla d(x)\nabla_\xi\omega (x,hD_x)g(P) \phi_j(d(x)
-C_1t) f(P)e^{-itP/h}\varphi_0\\
&&\hskip 4cm+{\mathcal  O}(\sum_{k=0}^Nh^{k+1}\Vert F_{j+k+1}f(P)e^{-itP/h}\varphi_0\Vert  + h^{N+2}),
\end{eqnarray*}
for any fixed  $N\geq 1$.
\vskip 0.3cm
Going back to (\ref{derivt}), and using the fact that $\Vert \nabla d(x)\nabla_\xi\omega (x,hD_x)g(P)\Vert\leq C_1 $, this gives,
\begin{eqnarray*}
&&\partial_t\la F_j(t,x)f(P)e^{-itP/h}\varphi_0,f(P)e^{-itP/h}\varphi_0\ra\\
&&\hskip 1.5cm\leq {\mathcal  O}(\sum_{k=0}^Nh^{k+1}\Vert F_{j+k+1}f(P)e^{-itP/h}\varphi_0\Vert^2 + h^{N+2}),
\end{eqnarray*}
and  therefore, integrating 
between $0$ and $t$, and using Lemma \ref{premlemappB}, 
\begin{eqnarray*}
&&\la F_j(t,x)f(P)e^{-itP/h}\varphi_0,f(P)e^{-itP/h}\varphi_0\ra\\
&&\hskip 1.5cm= {\mathcal  O}(\sum_{k=0}^Nh^{k+1}\int_0^t\Vert F_{j+k+1}f(P)e^{-isP/h}\varphi_0\Vert^2ds  + th^{N+2}),
\end{eqnarray*}
In particular, since 
$$\Vert F_j(t,x)f(P)e^{-itP/h}\varphi_0\Vert^2\leq \la F_j(t,x)f(P)e^{-itP/h}\varphi_0,f(P)e^{-itP/h}\varphi_0\ra,$$ 
we have $\Vert F_j(t,x)f(P)e^{-itP/h}\varphi_0\Vert^2={\mathcal  O}(h)$ for any $j\in  {\N},$ and then,  by
induction,  $\Vert F_j(t,x)f(P)e^{-itP/h}\varphi_0\Vert^2={\mathcal  O}(h^N)$
for all $N\in  {\N}$. Due to the definition (\ref{defFj}) of $F_j$, this proves the theorem.
\hfill$\bullet$
\chapter{Two Technical Lemmas}
\label{technr}
\setcounter{equation}{0}%
\setcounter{theorem}{0}%
\begin{lemma}\sl 
\label{tec2}
Let $\psi_j, \i_j\in C_0^\infty (\R^n)$, such that  $\i_j =1$ near $\supp \psi_j$. Then, for any $f\in C_0^\infty (\R)$, one has,
$$
U_j\psi_j f(\tilde P)U_j^{-1}\i_j = \psi_j f (U_j\i_j\tilde P U_j^{-1}\i_j) +{\mathcal  O}(h^\infty).
$$
\end{lemma}
{\em Proof -- } By (\ref{calfonc}), and taking the adjoint, it is enough to prove, for any $N\geq 1$,
$$
U_j\i_j (\tilde P -z)^{-1}U_j^{-1}\psi_j =  (U_j\i_j\tilde P U_j^{-1}\i_j -z)^{-1}\psi_j +{\mathcal  O}(h^N|\im z|^{-N'}),
$$
locally uniformly for $z\in\C$, and with some $N'=N'(N)<+\infty$. Let $v\in L^2(\R^n)$ and set $u:=(\tilde P -z)^{-1}U_j^{-1}\psi_jv$. By  Lemma \ref{suppdisj} (and its proof), we know that, 
\be
\label{tecc3}
u = \i_j u + {\mathcal  O}(h^N |\im z|^{-N'}\Vert v\Vert ),
\ee
for some $N'=N'(N)<+\infty$. On the other hand, we have,
\begin{eqnarray*}
(U_j\i_j\tilde P U_j^{-1}\i_j -z)U_j\i_j u&=& U_j\i_j\tilde Pu - zU_j\i_j u + U_j\i_j\tilde P (\i_j^2 -1)u\\
&=& U_j\i_j (zu + U_j^{-1}\psi_jv)- zU_j\i_j u + U_j\i_j\tilde P (\i_j^2 -1)u\\
&=& \psi_j v+ U_j\i_j\tilde P (\i_j^2 -1)u,
\end{eqnarray*}
and thus, using (\ref{tecc3}), 
\begin{eqnarray*}
U_j\i_j u &=& (U_j\i_j\tilde P U_j^{-1}\i_j -z)^{-1}(\psi_j v + U_j\i_j\tilde P (\i_j^2 -1)u)\\
&=& (U_j\i_j\tilde P U_j^{-1}\i_j -z)^{-1}\psi_j v +{\mathcal  O}(h^N |\im z|^{-N''}\Vert v\Vert ),
\end{eqnarray*}
for some other $N''=N''(N)<+\infty$. Then, the result follows.\hfill$\bullet$
\begin{lemma}\sl 
\label{tec1}
Let $\psi, \i\in C_0^\infty (\R^n)$, such that  $\i =1$ near $\supp \psi$. Then, for any $\rho\in C_0^\infty (\R)$, one has,
$$
\rho (\i \boldsymbol{\omega} \i)\psi = \rho (\boldsymbol{\omega})\psi +{\mathcal  O}(h^\infty).
$$
\end{lemma}
{\em Proof -- } The proof is very similar to (but simpler than) the one of Lemma \ref{tec2}, and we omit it. \hfill$\bullet$

\backmatter
\bibliographystyle{amsalpha}
\bibliography{}

\printindex

\end{document}